\documentclass[preprint,a4paper,preprint,12pt]{elsarticle}
\usepackage[top=3cm,right=2cm,bottom=2cm,left=3cm]{geometry}

\usepackage{graphicx}
\usepackage{amssymb}
\usepackage{amsmath}
\usepackage{amsthm}
\usepackage{bm}
\usepackage{enumerate}
\usepackage{color}
\usepackage{multirow}
\usepackage{booktabs}
\usepackage[table,xcdraw]{xcolor}
\usepackage{ulem}

\usepackage{lineno}

\theoremstyle{definition}

\theoremstyle{proposition}
\newtheorem{proposition}{Proposition}[section]
\theoremstyle{lemma}

\theoremstyle{theorem}
\newtheorem{theorem}{Theorem}[section]
\theoremstyle{corollary}

\newcommand{\vect}[1]{{\underline{#1}}}
\newcommand{\vetor}[1]{{\vect{#1}}}
\newcommand{\tensor}[1]{{\underline{\underline{#1}}}}

\begin{document}

\begin{frontmatter}


\title{ Error estimates for the Scaled Boundary Finite Element Method}



\author[fec]{Karolinne O. Coelho \corref{cor1}}
\ead{karolinneoc@gmail.com}

\author[fec]{ Philippe R. B. Devloo }

\ead{phil@fec.unicamp.br}

\author[IMECC]{ S\^onia M. Gomes }

\ead{soniag@unicamp.br}

\address[fec]{ FEC - Universidade Estadual de  Campinas, R. Josiah Willard Gibbs 85 - Cidade Universitária, Campinas, SP, CEP 13083-839, Brazil}

\address[IMECC]{ IMECC -   Universidade Estadual de  Campinas, Campinas, SP, Brazil}

\cortext[cor1]{Corresponding author - Phone +55 19 3521-1149}

\begin{abstract}
The Scaled Boundary Finite Element Method (SBFEM) is a technique  in which  approximation spaces are constructed using a semi-analytical approach. They are based on  partitions of the computational domain by  polygonal/polyhedral subregions, where  the shape functions  
approximate local Dirichlet problems with piecewise polynomial trace data. Using this operator adaptation approach, and by imposing a  starlike scaling requirement  on the subregions, the   representation of  local SBFEM  shape functions in  radial and surface directions are obtained  from eigenvalues and eigenfunctions of  an ODE system, whose coefficients are determined by the element geometry and the trace polynomial spaces.  The aim of this paper is to derive a priori error estimates  for SBFEM's  solutions of harmonic test problems. 
 For that, the SBFEM spaces are characterized in the context  of Duffy's approximations for which  a gradient-orthogonality  constraint is imposed.  As a consequence,  the scaled boundary functions are
gradient-orthogonal  to any  function in Duffy's spaces  vanishing at the mesh skeleton, a mimetic version of a well-known property valid for harmonic functions.  This orthogonality property is applied to provide  a priori SBFEM error estimates  in terms of known finite element interpolant errors of the exact solution. Similarities  with  virtual harmonic approximations are also explored  for the understanding of SBFEM  convergence properties.  Numerical experiments with  2D and 3D polytopal meshes confirm optimal   SBFEM convergence rates for two test problems with smooth solutions. Attention is  also paid to the approximation of a  point singular solution by  using SBFEM close to the singularity and finite element approximations elsewhere,  revealing optimal accuracy rates of  standard regular contexts.
\end{abstract}


\end{frontmatter}


\section{Introduction }
\label{S:1}

The Scaled Boundary Finite Element Method (SBFEM) is a Galerkin method  in which the approximation spaces are constructed using a semi-analytical approach \cite{Song1997, Song1998,Wolf2003,Song2018a}. They  are based on general  partitions of the computational domain by polygonal/polyhedral subregions $S$ (called $S$-elements), which are supposed to verify  the  starlike scaling requirement  such that any point on the boundary  $\partial S$ can be directly
 visible from a  center point (scaling center). The shape functions are computed  by the application of the scaled boundary technique,  involving  a specific parametrization of the  $S$-elements,  which is possible thanks to  their scaling property.
In classical FE methods, the local approximations are (mapped) polynomials, 
 and these are known to fail or 
  have very low convergence rates when the exact solutions can not be 
 properly represented by polynomials. 
 In SBFEM,  discretization by piecewise polynomials only takes place  at   $\partial S$, whilst   the  functions are  
 constructed by approximating  local
 Dirichlet problems  internally to $S$. The method is discussed in the books   \cite{Wolf2003,Song2018a}, and articles therein cited.

This  incorporation of  analytic knowledge about the local behavior of the exact solution in the  approximation spaces is the main  property  of SBFEM. 
 Therefore, it can be viewed as an operator adapted method.  As discussed in   \cite{Babuska1997,Melenk1999}, in the context of  the Partition of Unity Method,   these methods can be expected to perform better when   compared with standard polynomial based FE approximations. 
They reduce the number of  degrees of freedom significantly 
and hence the computational cost, while improving the quality of the solutions.  
For the case of SBFEM, where only  boundary  values of the subdomains are discretized by local surface polynomials, its operator adapted approach revealed  itself to be particularly efficient to approximate problems with stress singularities, such as crack tips, v-notches, and re-entrant corners to name a few applications in elasticity \cite{Yang2006,Song2018,Pramod2019,Guo2019,Bulling2019}. More recently, the method has been applied to highly irregular and heterogeneous domains due to the flexibility in generating SBFEM meshes \cite{Liu2017,Natarajan2020}. For instance, the SBFEM has been applied in quadtree and octree meshes since hanging nodes  
 can be  avoided due to the flexible topology of SBFEM polygonal/polyhedral subregions  \cite{Saputra2017,Chen2018,Guo2019}.
 
 The aim of this paper is to derive a priori error estimates  of SBFEM  
approximations for the case of Laplace’s equation. Although numerical experiments in the literature point that optimal rates of convergence are obtained using SBFEM approximations \cite{Gravenkamp2019,Gravenkamp2020}, the mathematical demonstrations that give support to the observed numerical results are new contributions of the current  work. For that, we explore two different aspects of SBFEM  spaces, shared with Duffy's approximations \cite{Duffy1982} or  with virtual harmonic spaces \cite{Chernov2019}.

Taking advantage of  the scaling property,   functions can be represented in the $S$-elements  by coordinates in  radial and surface directions. Their values on the boundary $\partial S$  live in piecewise polynomial trace spaces, which are radially extended to the  interior of the subdomain. Therefore, this property puts  SBFEM's spaces in the context of Duffy's approximations \cite{Duffy1982},  whose definitions are summarized in Section \ref{sec:duffys}. 
Partitions of  $S$  are 
obtained by a geometric transformation collapsing a reference quadrilateral, hexahedron or prism on  triangular, pyramidal or tetrahedral elements $K\subset S$, each one sharing the scaling center as a vertex
(see  Section  \ref{sec:sbfemin}). 
 
 For the model Laplace  problem under consideration,  the focus of the SBFEM  operator adapted approach is   the  approximations  inside $S$-elements   by  “radial harmonic extensions”  of   surface components. 
 It is shown that  SBFEM's  spaces are  Duffy's approximations 
constructed to solve  Laplace problems with piecewise polynomial Dirichlet data over $\partial S$.  
SBFEM spaces are characterized by the enforcement of a gradient-orthogonality constraint with respect to Duffy’s approximations vanishing on $\partial S$ and at the center point, as demonstrated in Section  \ref{SBFEM-S}. 
This perspective on SBFEM approximations reveals 
that the local scaled boundary shape functions are constructed based on an orthogonality condition. 
By enforcing these intrinsic orthogonality constraints, their parametrization in   radial and surface directions emerge from the eigenvalues and eigenfunctions of an ODE system, whose coefficients are determined by the element geometry and the trace polynomial spaces. 
The scaled boundary functions are gradient-orthogonal to  an
extended class of  Duffy’s  functions  that vanish at the mesh skeleton.  It can be viewed as  a mimetic version of a well-known property valid for harmonic functions. These aspects are stated  in  Proposition \ref{ExtendedOrtho} and used as a key tool to the development  of energy error estimates for the SBFEM in terms of FE interpolation errors in Section \ref{sec:SBFEMGalerkin}, as  shown in Theorem \ref{apriori}, one of the main contributions of this study. 
    
SBFEM also has close similarities with virtual harmonic approximation spaces recently introduced in \cite{Chernov2019},  as explored in Section \ref{Comments} and  summarized in Theorem \ref{apriori2}.   In both  cases the   trace functions are piecewise polynomials   defined over subregion boundaries $\partial S$, which are extended to the interior of $S$ by solving local Dirichlet Laplace problems: 
whilst 
the  functions in the local virtual spaces   are strongly harmonic, in  SBFEM spaces this property is enforced  in a reduced extent.  Thus, SBFEM approximation errors may come from the   trace polynomial interpolation  or by their deviation of being harmonic. 
However, unlike for the virtual harmonic subspaces,  it is possible   to explore  the  radial Duffy's structure  to  explicitly compute  SBFEM shape functions.

In Section \ref{sec:numericaltests}, we present  results of  SBFEM computational simulations for some harmonic test problems confirming the predicted theoretical convergence results of Section \ref{sec:SBFEMGalerkin}. We consider  2D and 3D cases with  smooth solutions, and discretizations based on different $S$-partition geometry, which  are formed by internally  collapsed triangular, pyramidal, or tetrahedral elements. In the same section, we present $p$-convergence histories  verifying  asymptotic exponential convergence rates in terms of degrees of freedom (DOF), and compare results with respect to the ones given by usual  FE methods based on partitions obtained by the conglomeration of the internal collapsed elements. We also pay  attention to the approximation of a singular problem where the singularity occurs  by the  change of boundary condition  and observe that optimal rates of convergence holds using few DOF, using SBFEM to resolve the singularity. We draw some concluding remarks  in Section \ref{sec:conclusions}.

\section{Duffy's approximations  in  triangles, pyramids and tetrahedra}\label{sec:duffys}

 Duffy's transformations \citep{Duffy1982} (also referred to as  collapsed coordinate systems)  are invertible maps of  a rectangle into a triangle, a hexahedron to a pyramid,  or a prism to a tetrahedron. These maps were originally proposed for integration of vertex singularities and they  are widely applied to define integration quadrature formulae in  triangles \citep{Lyness1994, Blyth2006}. Duffy's transformations are also the basic tools for the construction of spectral methods  on simplices  (triangles, tetrahedra) \citep{Karniadakis1999}. Collapsed isoparametric elements parametrized by Duffy's transformations also have  applications in crack problems \cite{Wu1993,Pu1978,Raju1987}. 

\subsection{Duffy's geometric transformations}
The  master elements to be considered   have  the general form $\hat{K} = [0,1] \times \hat{L}\subset \mathbb{R}^d$, where $\hat{L} \subset \mathbb{R}^{d-1}$, $d=2,3$.
In the parametric coordinates $\hat{\mathbf{x}}=(\xi,\bm{\eta}) \in \hat{K}$,   $\xi$ plays the role of radial variable,  and  $\bm{\eta}$ refers to surface coordinates. The geometry of the master elements may be one of the following  kinds:
\begin{itemize}
\item Rectangle $\hat{K}$, where $\hat{L}= \hat{I}$ is the interval $\hat{I}=[-1,1]$.
\item Hexahedron $\hat{K}$, where   $\hat{L}=\hat{Q}$ is the rectangle $\hat{Q}=[-1,1] \times [-1,1]$.
\item Prism  $\hat{K}$, where   $\hat{L}= \hat{T}$ is the triangle $\hat{T}=\{\boldsymbol{\eta}=(\eta_1,\eta_2); 0 \leq \eta_i \leq 1, \eta_1+\eta_2\leq 1\}.$
\end{itemize}
The key aspect of  geometric  Duffy's transformations  $F_K : \hat{K} \rightarrow K$  is  the collapse  of one facet in $\hat{K}$ on a single vertex of the deformed element $K$. These maps are also referred  in the literature as  collapsed coordinate systems  \citep{Karniadakis1999}.   If  $\mathbf{x}$ denotes   the Cartesian coordinate in $K$,  the mapped points 
 $\mathbf{x} = F_K(\xi,\bm{\eta}) \in K$  are generically defined by
\begin{equation}
	F_K(\xi,\bm{\eta}) = \xi \left( F_L(\bm{\eta}) - \mathbf{a}_0\right) + \mathbf{a}_0,  \label{eq:Duffy}
\end{equation}
where     $\mathbf{a}_0$ is  a vertex in $K$,  and $L \subset \partial K$ refers to a facet  opposite to $\mathbf{a}_0$, which is supposed to be mapped by the  geometric transformation  $F_L: \hat{L}\rightarrow L$. Notice that the whole facet $\{ (0,\bm{\eta}), \bm{\eta} \in \hat{L}\} \subset \hat{K}$ is collapsed over the vertex  $\mathbf{a}_0 \in K$, so that $K$ can be regarded as a quadrilateral with two identical vertices,  a hexahedron with  four equal vertices, or a prism with three identical vertices. That is why $\mathbf{a}_0$ is called the collapsed vertex.
The mapping  $F_K$ can also be seen as a scaling from a point $F_L(\boldsymbol{\eta})\in L$ to the vertex $\mathbf{a}_0$. This process  generates  radial lines  $[ \mathbf{a}_0, F_L(\boldsymbol{\eta})] = \mathbf{a}_0+ \xi \mathbf{r}(\boldsymbol{\eta})$,  where $\mathbf{r}(\boldsymbol{\eta})=F_L(\boldsymbol{\eta})-\mathbf{a}_0$.

The Jacobian matrix $\mathbf{J}_K = \nabla_{\hat{\mathbf{x}}} F_K$ of the transformation \eqref{eq:Duffy}  is
\begin{align}
\mathbf{J}_K (\xi,\bm{\eta}) &= \begin{bmatrix} F_L(\bm{\eta}) - \mathbf{a}_0 & \xi\nabla_{\bm{\eta}} F_L(\bm{\eta}) \end{bmatrix} = \mathbf{J}_K(1,\bm{\eta}) \begin{bmatrix} 1 & 0 \\ 0 & \xi\mathbf{I}_{d-1} \end{bmatrix}, 
\label{eq:jmatrix}
\end{align}
where $\mathbf{I}_{d-1}$ is the $d-1\times d-1$ identity matrix, and 
$
\mathbf{J}_{K}(1,\bm{\eta})=\begin{bmatrix} F_L(\bm{\eta}) - \mathbf{a}_0 & \nabla_{\bm{\eta}}F_L(\bm{\eta})  \end{bmatrix}$
is the Jacobian matrix at the surface points where $F_K(1,\bm{\eta}) = F_L(\bm{\eta})$.  Thus
\begin{align}
	\mathbf{J}_K^{-1} = \begin{pmatrix} 1 & 0 \\ 0 & \frac{1}{\xi} \mathbf{I}_{d-1}
	\end{pmatrix}\mathbf{J}_K(1,\bm{\eta})^{-1}.
\label{eq:jacinv}
\end{align}

In the following, the geometric transformation \eqref{eq:Duffy} is illustrated for the three different element geometries 
considered in the current study.

\subsubsection*{Case 1:  quadrilateral $\hat{K}$ to triangular $K$}
\begin{figure}[h]
    \begin{centering}
    \includegraphics[scale=0.3]{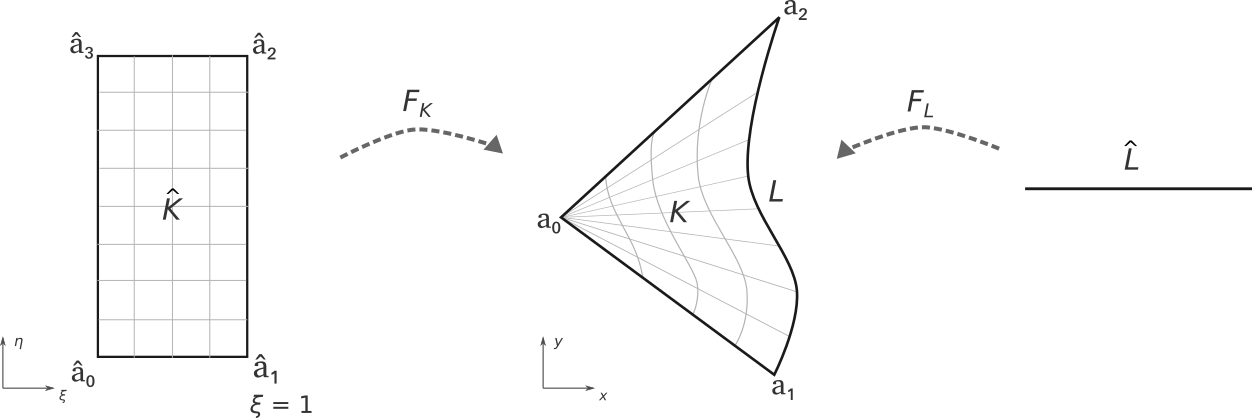}
    \par\end{centering}
    \caption{Geometric illustration of Duffy's transformation over a triangle as a collapsed quadrilateral. \label{fig:DuffyTriangle}}
\end{figure}

Let $\hat{K}$ be the rectangular master element with vertices  listed in the next table
\setlength{\tabcolsep}{0.1pt}
\begin{center}
\begin{tabular}{|c|c|c|c|}
\hline 
$\hat{\mathbf{a}}_{0}$ & $\hat{\mathbf{a}}_{1}$ & $\hat{\mathbf{a}}_{2}$ & $\hat{\mathbf{a}}_{3}$\tabularnewline
\hline 
$(0,-1)$ & $(1,-1)$ & $(1,1)$ & $(0,1)$\tabularnewline
\hline 
\end{tabular}
\end{center}
 and consider a  general  triangular element, with vertices  $\mathbf{a}_0 =  F_K(\hat{\mathbf{a}}_0)$, $\mathbf{a}_1 =  F_K(\hat{\mathbf{a}}_1)$, and $\mathbf{a}_2 =  F_K(\hat{\mathbf{a}}_2)$,  as   illustrated in Figure \ref{fig:DuffyTriangle}. Notice that the edge  $[\hat{\mathbf{a}}_0, \hat{\mathbf{a}}_3]$  collapses onto the vertex $\mathbf{a}_0 = F_K(\hat{\mathbf{a}}_0) \in K$, whilst   $\mathbf{a}_1$ and  $\mathbf{a}_2$ are the  vertices of the opposite edge   $L=F_K(1,\eta) = F_L(\eta)$.

\subsubsection*{Case 2: hexahedral $\hat{K}$ to pyramidal $K$}
\begin{figure}[h]
    \begin{centering}
    \includegraphics[scale=0.3]{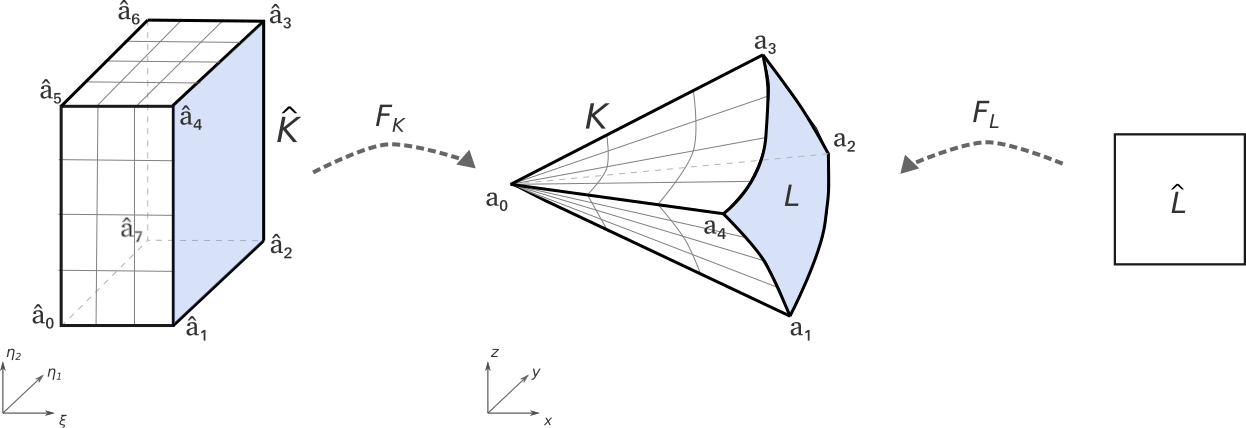}
    \par\end{centering}
    \caption{Geometric illustration of a Duffy's transformation over a  pyramid  as a collapsed hexahedron. 
    \label{fig:DuffyPyramid}}
\end{figure}

The master element is the hexahedron $\hat{K}$ whose vertices are listed bellow.
\setlength{\tabcolsep}{0.1pt}
\begin{center}
\begin{tabular}{|c|c|c|c|c|c|c|c|}
\hline 
$\hat{\mathbf{a}}_{0}$ & $\hat{\mathbf{a}}_{1}$ & $\hat{\mathbf{a}}_{2}$ & $\hat{\mathbf{a}}_{3}$ & $\hat{\mathbf{a}}_{4}$ & $\hat{\mathbf{a}}_{5}$ & $\hat{\mathbf{a}}_{6}$ & $\hat{\mathbf{a}}_{7}$\tabularnewline
\hline 
$(0,-1,-1)$ & $(1,-1,-1)$ & $(1,1,-1)$ & $(1,1,1)$ & $(1,-1,1)$ & $(0,-1,1)$ & $(0,1,1)$ & $(0,1,-1)$\tabularnewline
\hline 
\end{tabular}
\end{center}
Figure \ref{fig:DuffyPyramid}  illustrates a  mapped pyramid  with vertices $\mathbf{a}_i =  F_K(\hat{\mathbf{a}}_i)$, $i=0, \cdots 4$,  $\mathbf{a}_0$ being the collapsed vertex  with opposite quadrilateral face  $L=[\mathbf{a}_1, \mathbf{a}_2, \mathbf{a}_3, \mathbf{a}_4]$.
Observe that: 
\begin{enumerate}
\item The rectangular face $\text{\ensuremath{\left[\hat{\mathbf{a}}_{0},\hat{\mathbf{a}}_{5},\hat{\mathbf{a}}_{6},\hat{\mathbf{a}}_{7}\right]}}$
collapses onto $\mathbf{a}_{0}$;
\item The face $\text{\ensuremath{\left[\hat{\mathbf{a}}_{0},\hat{\mathbf{a}}_{1},\hat{\mathbf{a}}_{4},\hat{\mathbf{a}}_{5}\right]}}$
 collapses onto the triangle $\text{\ensuremath{\left[{\mathbf{a}}_{0},{\mathbf{a}}_{1},{\mathbf{a}}_{4}\right]}}$;
\item The face $\text{\ensuremath{\left[\hat{\mathbf{a}}_{0},\hat{\mathbf{a}}_{1},\hat{\mathbf{a}}_{2},\hat{\mathbf{a}}_{7}\right]}}$
 collapses onto the triangle $\text{\ensuremath{\left[{\mathbf{a}}_{0},{\mathbf{a}}_{1},{\mathbf{a}}_{2}\right]}}$; 
\item The face $\text{\ensuremath{\left[\hat{\mathbf{a}}_{2},\hat{\mathbf{a}}_{3},\hat{\mathbf{a}}_{6},\hat{\mathbf{a}}_{7}\right]}}$
 collapses onto the triangle $\text{\ensuremath{\left[{\mathbf{a}}_{0}, {\mathbf{a}}_{2},{\mathbf{a}}_{3}\right]}}$;
 \item The face $\text{\ensuremath{\left[\hat{\mathbf{a}}_{3},\hat{\mathbf{a}}_{6},\hat{\mathbf{a}}_{5},\hat{\mathbf{a}}_{4}\right]}}$
 collapses onto the triangle $\text{\ensuremath{\left[{\mathbf{a}}_{0}, {\mathbf{a}}_{4},{\mathbf{a}}_{3}\right]}}$.\\
\end{enumerate}

\subsubsection*{Case 3:  prismatic $\hat{K}$ to tetrahedral $K$}
\begin{figure}[h]
    \begin{centering}
    \includegraphics[scale=0.3]{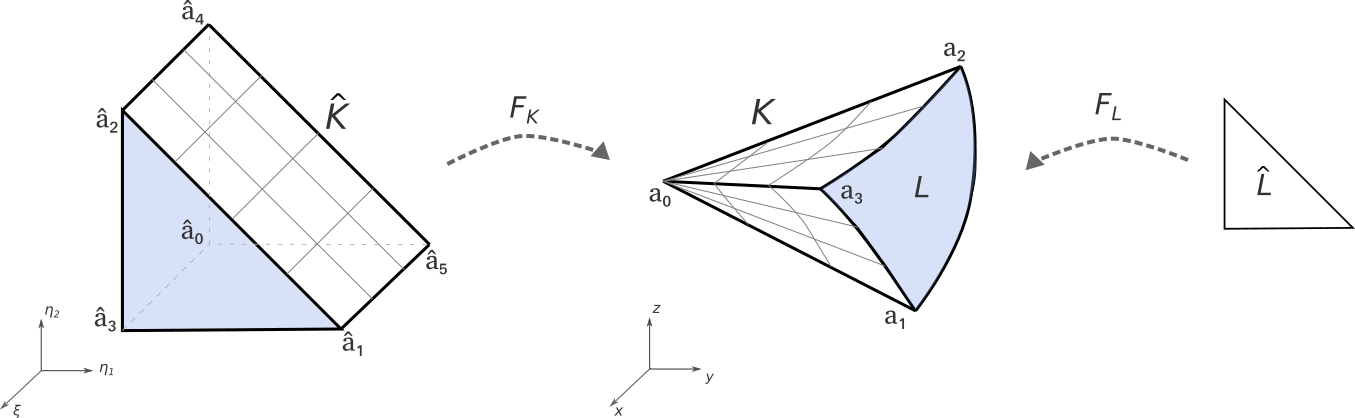}
    \par\end{centering}
    \caption{Geometric illustration of a Duffy's transformation over a  tetrahedron  as a collapsed prism. 
    \label{fig:DuffyThetra}}
\end{figure}

The  master element is the prism $\hat{K}$ whose vertices are listed bellow.
\setlength{\tabcolsep}{0.1pt}
\begin{center}
\begin{tabular}{|c|c|c|c|c|c|}
\hline 
$\hat{\mathbf{a}}_{0}$ & $\hat{\mathbf{a}}_{1}$ & $\hat{\mathbf{a}}_{2}$ & $\hat{\mathbf{a}}_{3}$ & $\hat{\mathbf{a}}_{4}$ & $\hat{\mathbf{a}}_{5}$ \tabularnewline
\hline 
$(0,0,0)$ & $(1,1,0)$ & $(0,1,0)$ & $(0,1,1)$ & $(1,0,0)$ & $(0,0,1)$ \tabularnewline
\hline 
\end{tabular}
\end{center}
In the tetrahedron   shown in Figure \ref{fig:DuffyThetra},  with vertices $\mathbf{a}_i =  F_K(\hat{\mathbf{a}}_i)$, $i=0, \cdots 3$,   the collapsed vertex is $\mathbf{a}_0$ and the opposite quadrilateral face  is $L=[\mathbf{a}_1, \mathbf{a}_2, \mathbf{a}_3]$.
Note  that: 
\begin{enumerate}
\item The triangular face $\text{\ensuremath{\left[\hat{\mathbf{a}}_{0},\hat{\mathbf{a}}_{4},\hat{\mathbf{a}}_{5}\right]}}$
 collapses onto the vertex $\mathbf{a}_{0}$;
\item The  quadrilateral face $\left[\hat{\mathbf{a}}_{0},\hat{\mathbf{a}}_{4},\hat{\mathbf{a}}_{2},\hat{\mathbf{a}}_{3}\right]$
 collapses onto the triangle $\left[{\mathbf{a}}_{0},{\mathbf{a}}_{2}, {\mathbf{a}}_{3}\right]$;
\item The   quadrilateral face $\left[\hat{\mathbf{a}}_{0},\hat{\mathbf{a}}_{3}, \hat{\mathbf{a}}_{1}, \hat{\mathbf{a}}_{5}\right]$
 collapses onto the triangle $\left[{\mathbf{a}}_{0},{\mathbf{a}}_{1},{\mathbf{a}}_{3}\right]$;
\item The  quadrilateral face $\left[\hat{\mathbf{a}}_{1},\hat{\mathbf{a}}_{2},\hat{\mathbf{a}}_{4},\hat{\mathbf{a}}_{5}\right]$
collapses onto the triangle $\left[{\mathbf{a}}_{0},{\mathbf{a}}_{1},{\mathbf{a}}_{2}\right]$.
\end{enumerate}
We recall  that a hexahedron to tetrahedron  Duffy's transformation can also be derived, as adopted in \citep{Karniadakis1999},  first via  a preliminary step hexahedron to prism, and then the prism to tetrahedron described above.

\subsection{Duffy's approximations }\label{sec:duffy}
 Duffy's approximations  refer to   functions  $\phi=\mathbb{F}_{K}(\hat{\phi})$ defined in  $K$ and  obtained  backtracking  functions $\hat{\phi}(\xi,\bm{\eta})$ defined in $\hat{K}$, meaning that
\[\phi(\mathbf{x})=\hat{\phi}(\xi,\boldsymbol{\eta}), \quad \mbox{for} \; \mathbf{x} ={F}_{K}(\xi,\boldsymbol{\eta}) \in K.\]
The focus of this paper is on functions  ${\phi}$ obtained by separating variables in  $\hat{\phi}(\xi,\boldsymbol{\eta})=\hat{\rho}(\xi)\hat{\alpha}(\boldsymbol{\eta})$, where $\hat{\rho}(\xi)$ is called the radial component, and  $\hat{\alpha}(\boldsymbol{\eta)}$ is the surface  component.  
 It is clear that   constant functions in $\hat{K}$ are  mapped to constant functions in $K$. It should  also be noted that for the cases where  $\hat{\alpha}(\bm{\eta})$ is  not a constant function, the well definition of  $\phi$ at the collapsed point $\mathbf{a}_0$ requires that $\hat{\rho}(0)=0$, so that $\phi(\mathbf{a}_0)=0$.

We consider function spaces  $${\mathcal{D}}_k (\hat{K})= \{ \hat{\phi}(\xi,\boldsymbol{\eta})=\hat{\rho}(\xi)\hat{\alpha}(\boldsymbol{\eta}); \; \hat{\alpha}(\boldsymbol{\eta}) \in  V_k(\hat{L}),$$
where the surface components $ \hat{\alpha}(\boldsymbol{\eta})\in  V_k(\hat{L})$  used to define FE approximation spaces $V_k(L)=\mathbb{F}_L(V_k(\hat{L}))$, are   finite dimensional  polynomial spaces $V_k(\hat{L})$. The following cases shall be studied:
\begin{enumerate}
\item $V_k(\hat{L})=\mathbb{P}_{k}(\hat{L})$,  polynomials  of total degree not greater than $k$, for the interval $\hat{L}=[-1,1]$ or for the triangle $\hat{L}=\hat{T}$.
\item $V_k(\hat{L})=\mathbb{Q}_{k,k}(\hat{L})$,   polynomials of  degree  not greater than $k$ on each coordinate $\eta_1, \eta_2$, for the quadrilateral $\hat{L}= \hat{Q}$.
\end{enumerate}

 \subsubsection*{Gradient operation in $\mathcal{D}_k ({K})$}

We  restrict the study to  mapped spaces $\mathcal{D}_k ({K})=\mathbb{F}_K({\mathcal{D}}_k (\hat{K}))\subset {H}^1({K})$. 
For instance, as already observed in \citep{Shen2008} for the case of triangular elements $K$,  $H^1(K)$ corresponds to ${H}^1_\omega(\hat{K})$ where
${H}^1_\omega(\hat{K}):=\{ \hat{\phi} \in L^2_\omega(\hat{K}):\;\partial_{\bm{\eta}} \,\hat{\phi} \in L^2_{\omega^{-1}}(\hat{K}) \; \mbox{and} \; \partial_\xi \,\hat{\phi} \in L^2_{\omega}(\hat{K})\}$,
where $\omega(\xi,\bm{\eta})=\xi |\mathbf{J}_{K}(1,\bm{\eta}) |$. Particularly, $\partial_{\bm{\eta}} \, \hat{\phi}(0,\bm{\eta})=0$ for bounded $\partial_y \phi(x,y)$. 

The chain rule implies that 
{\small
 \begin{align}
\nabla_{\mathbf{x}}\phi(\mathbf{x}) & =[\mathbf{J}_{K}(1,\boldsymbol{\eta})]^{-T}\left[\begin{array}{cc}
1 & 0\\
0 & \frac{1}{\xi} \mathbf{I}_{d-1}
\end{array}\right]\left[\begin{array}{c}
\hat{\rho}'(\xi)\hat{\alpha}(\boldsymbol{\eta}) \nonumber\\
\hat{\rho}(\xi)\nabla_{\mathbf{\boldsymbol{\eta}}}\hat{\alpha}(\boldsymbol{\eta})
\end{array}\right]\\
 & =[\mathbf{J}_K(1,\boldsymbol{\eta})]^{-T}\left[\begin{array}{c}
\hat{\rho}'(\xi)\hat{\alpha}(\boldsymbol{\eta})\\
\frac{1}{\xi}\hat{\rho}(\xi)\nabla_{\boldsymbol{\eta}}\hat{\alpha}(\boldsymbol{\eta})
\end{array}\right] 
  =[\mathbf{J}_K(1,\boldsymbol{\eta})]^{-T}\left[\begin{array}{cc}
\hat{\alpha}(\boldsymbol{\eta}) & 0\\
0 & \nabla_{\boldsymbol{\eta}}\hat{\alpha}(\boldsymbol{\eta})
\end{array}\right]\left[\begin{array}{c}
\hat{\rho}'(\xi)\\
\frac{1}{\xi}\hat{\rho}(\xi)
\end{array}\right]. \label{eq:Grad}
\end{align}}
If $\hat{\alpha}(\boldsymbol{\eta})=\sum_{l}\alpha^{l}\hat{N}_k^{l}(\boldsymbol{\eta})$
is a linear combination of FE shape functions $\hat{N}_k^{l}(\boldsymbol{\eta})$
forming a basis for $V_k(\hat{L})$, then 
\begin{equation}
\nabla_{\mathbf{x}}\phi(\mathbf{x})=\sum_{l}\alpha^{l}\left[\begin{array}{cc}
\underline{B}_{1l}(\boldsymbol{\eta}) & \underline{B}_{2l}(\boldsymbol{\eta})\end{array}\right]\left[\begin{array}{c}
\hat{\rho}'(\xi)\\
\frac{1}{\xi}\hat{\rho}(\xi)
\end{array}\right], \label{eq:Grad2}
\end{equation}
where 
\begin{align}
\underline{B}_{1l}(\boldsymbol{\eta})  =[\mathbf{J}_K(1,\boldsymbol{\eta})]^{-T}\left[\begin{array}{c}
\hat{N}_k^{l}(\boldsymbol{\eta})\\
0
\end{array}\right], \; \mbox{and}  \; \;
\underline{B}_{2l}(\boldsymbol{\eta})  =[\mathbf{J}_K(1,\boldsymbol{\eta})]^{-T}\left[\begin{array}{c}
0\\
\nabla_{\boldsymbol{\eta}}\hat{N}_k^{l}(\boldsymbol{\eta})
\end{array}\right].\label{eq:BVect}
\end{align}

 \subsubsection*{Special case: $\hat{\alpha}(\bm{\eta}) \equiv 1$ }\label{Example}
 Let us  consider the particular cases of  $\phi(\mathbf{x})  \in \mathcal{D}_k ({K})$, for which  $\hat{\phi}(\xi,\bm{\eta})=\hat{\rho}(\xi)$, meaning that $\hat{\alpha}(\bm{\eta}) \equiv 1$.  
  A closer look on formula \eqref{eq:Grad} reveals that  
 \begin{align}
\nabla_{\mathbf{x}} \phi (\mathbf{x}) & =  [\mathbf{J}_K(1,\boldsymbol{\eta})]^{-T}    \hat{\rho}^{'}(\xi).   \label{eq:GradC} 
\end{align}
 For   affine elements $K$ and   $\hat{\rho}(\xi)=\xi$,  the mapped  function has constant gradient normal to $L$, so that $\phi   \in H^1(K)$  is an affine function vanishing at the collapsed vertex  $\mathbf{a}_1$,  and constant unitary values $\phi|_{L}\equiv 1$ over the  facet $L$  opposite to  $\mathbf{a}_1$. 

\section{SBFEM  spaces in the context of Duffy's approximations} 
\label{sec:scaledboundary}
Our purpose in this section  is to summarize the main aspects  of  SBFEM  approximation spaces under the point of view of Duffy's approximations and to prove  some  of their orthogonality properties  to a large range of $H^1$-conforming  functions.   

 \subsection{S-elements}\label{sec:Selements}

\begin{figure}[!htb]. 
    \centering
    \includegraphics[scale=0.3]{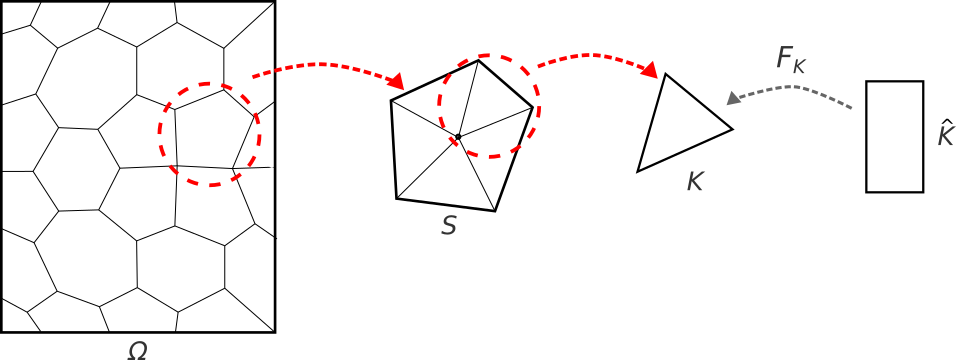}\vspace{0.4cm}\\
    \includegraphics[scale=0.25]{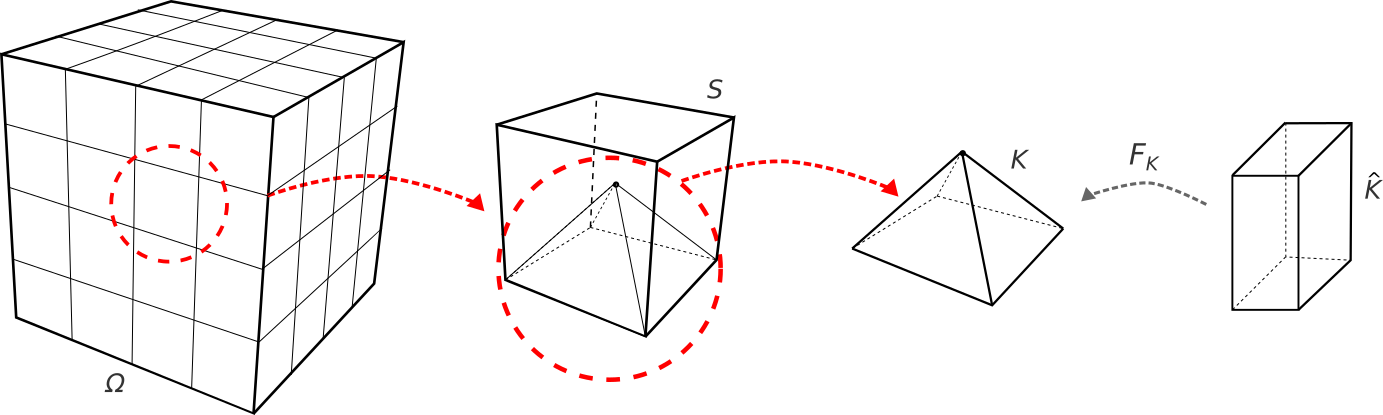}\vspace{0.4cm}\\
    \includegraphics[scale=0.25]{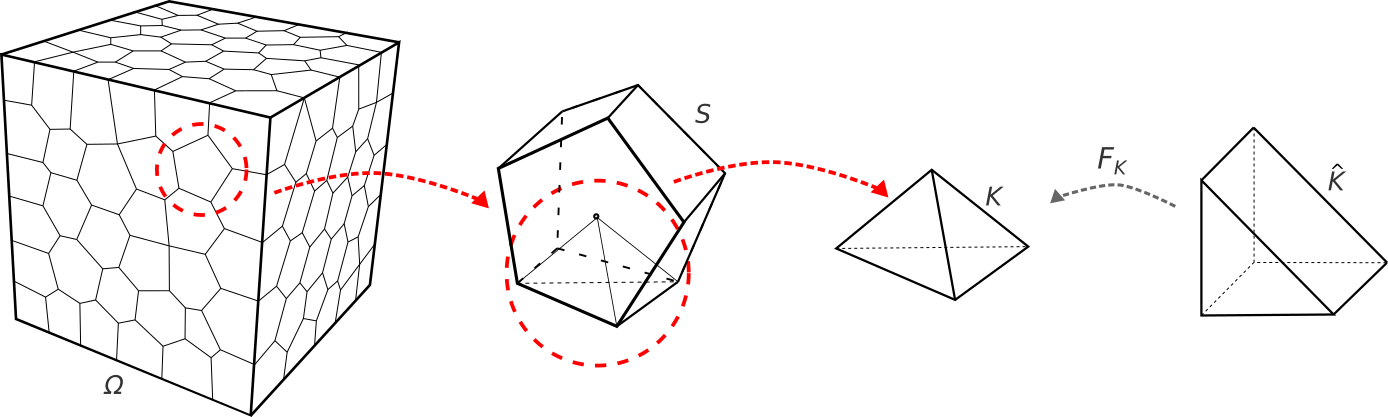}
    \caption{Illustration of  macro  partitions  $\mathcal{T}= \{S\}$,   with focus on a  sector  $ K \in \mathcal{T}^S$, with corresponding Duffy's transformation, for  triangular, pyramidal and  tetrahedral $K$. 
    \label{fig:partition}}
\end{figure}

The SBFEM  adopts  macro partitions $\mathcal{T}=\{S\}$ of the computational domain $\Omega  \subset \mathbb{R}^d$ by subregions  $S$ verifying the  starlike scaling requirement  that any point on the  boundary  of $S$ should be directly visible from a  point  $\mathbf{O} \in S$, called the scaling center. We restrict  the study  to   convex polytopal  $S$-elements   (polygonal or  polyhedral with flat facets $L^e$).  In the literature covering this method, the set  $\Gamma^S=\cup_e L^e$,   $e=1,\cdots, N^{\Gamma^S}$  is known as the scaled boundary element.   A  conformal  sub-partition  $\mathcal{T}^S =\{K^e\}$  of $S$ is formed by  sectors $K^e$ sharing  the  scaling center $\mathbf{O}$ as one of their vertices,   $L^e$  being  the  facet of  $K^e$  opposite to the scaling center.   As illustrated in  Figure \ref{fig:partition}, the sectors $K^e$  may have different geometry: triangular in 2D, pyramidal, or tetrahedral in 3D, the facets $L^e$  being a line segment, a quadrilateral or a triangular  element, respectively.  Moreover, we notice that a three-dimensional  $S$-element   may also be partitioned by hybrid tetrahedral-pyramidal meshes, combining  elements of different geometry,  with scaled boundary $\Gamma^S$  formed by triangular-quadrilateral facets.  For simplicity, we shall  restrict the analysis to  partitions   $\mathcal{T}^S$ where all elements $K^e$ have the same geometry.

  This scaled geometry of $S$  implies  that  the points $\mathbf{x} \in S$ can be uniquely represented by a radial coordinate $0\leq \xi \leq 1$  and a surface coordinate  $ \mathbf{x}_b$. The radial  coordinate   (or scaling factor)  points from the scaling center ($\xi = 0$) to a point $\mathbf{x}_b\in \Gamma^S$ (where $\xi=1$).   
The geometry of   $S$ may also be defined in each   sector $K^e \in  \mathcal{T}^S$   by  a transformation from the cartesian  coordinates $\mathbf{x} \in K^e$ to parametric Duffy's coordinates $(\xi, \boldsymbol{\eta}) \in \hat{K}= [0,1]\times \hat{L}$.   This   correspondence  defines  a  geometric mapping $F_{K^e}: \hat{K} \rightarrow K^e$ in the class of Duffy's transformations described in  the previous section, where $K^e$ is interpreted as a collapsed quadrilateral, hexahedral or prismatic geometric element  for which the facet $ F_{K^e}(0,\boldsymbol{\eta})$ is collapsed on top of its  vertex $\mathbf{x}_0$ in  the scaling center $\mathbf{O}$. The points $\mathbf{x}_b$ in  the  opposed facet $L^e$ are expressed as  $F_{K^e}(1,\boldsymbol{\eta})=F_{L^e}(\boldsymbol{\eta})$, $\boldsymbol{\eta}\in \hat{L}$.  For hexahedral or prismatic reference  elements $\hat{K}$, the  lateral quadrilateral faces are collapsed on triangular faces to form a pyramid or a tetrahedron, respectively. These maps are illustrated in Figure \ref{fig:partition}.
 
 \subsection{Duffy's spaces in  S-elements}\label{sec:sbfemin}
There are two stages in the construction of approximations  on  polytopal elements $S$:
\begin{itemize}
\item[ 1)] Definition of a trace space over the  boundary   $\Gamma^S$.
\item[ 2)] Extension of the traces  to the interior of  $S$. 
\end{itemize}
The first stage is typical of FE contexts,  but for specific scaled $S$-elements  the extension  to the interior can be performed in the radial direction.
\subsubsection*{Trace FE space over the scaled boundary $\Gamma^S$}
Let  $\Lambda_k(\Gamma^S) = C(\Gamma^S) \cap \prod_{L^e \subset \Gamma^S} V_k(L^e)$ be a FE  space defined over $\Gamma^S$. 
Recall that  $V_k(L^e)=\mathbb{F}_{L^e}(V_k(\hat{L}))$, where $V_k(\hat{L})$ is the polynomial space considered in $\hat{L}$. 
Let  $ N_k^{l,e}=\mathbb{F}_{L^e}(\hat{N}^l_k)$ be shape functions for the local FE spaces $V_k(L^e)$ over the facets $L^e \subset \Gamma^S$  obtained backtracking polynomial shape functions  $\hat{N}^l_k$ for the reference polynomial space $V_k(\hat{L})$. Thus, if  $\alpha \in \Lambda_k(\Gamma^S)$  and $\mathbf{x}_b=F_{L^e}(\boldsymbol{\eta}) \in L^e$, then $\alpha(\mathbf{x}_b)  = \hat{\alpha}^e(\boldsymbol{\eta}) =\sum_l {\alpha}^{l,e} \hat{N}_k^l (\boldsymbol{\eta})$. As usual, shape functions $N_k^{n,S}(\mathbf{x})$  for  $\Lambda_k(\Gamma^S)$ (say, of cardinality $\mathcal{N}^S$) can be obtained by the assembly of the local  shape functions $N_k^{l,e}$, and   the functions $\alpha \in \Lambda_k(\Gamma^S)$ can globally represented by  linear combinations
$\alpha(\mathbf{x}_b)  = \sum_{n=1}^{\mathcal{N}^S} {\alpha}^n {N}_{k}^{n,S} (\mathbf{x}_b), \mathbf{x}_b\in \Gamma^S$.  By collecting  the shape functions and  multiplying coefficient  in $\mathcal{N}^S$-vectors $\vetor{N}^S=[N^{n,S}_{k}]$ and   $\vetor{\alpha}=[\alpha^n]$, we may use the alternative expression $\alpha = \vetor{N}^S \cdot  \vetor{\alpha}$.

 \subsubsection*{  Radial  extensions: Duffy's space over $S$}
Given a trace function $\alpha \in   \Lambda_k(\Gamma^S)$, take a radial function $\hat{\rho}(\xi)$, $0\leq \xi \leq 1$, to  induce the  definition of a  function $\phi(\mathbf{x})$  by radial extension to the interior of $S$.  Inside each  sector   $K^e \in \mathcal{T}^S$ and  over $L^e$,  consider the parametrizations  $\mathbf{x}=F_{K^e}(\xi,\boldsymbol{\eta}) \in K^e$ and   $\mathbf{x}_b=F_{L^e}(\boldsymbol{\eta})$.  Recall the representation $\alpha(\mathbf{x}_b)  = \hat{\alpha}^e(\boldsymbol{\eta})$ to  define the radial extension $$\phi(\mathbf{x})=\hat{\phi}^e(\xi,\boldsymbol{\eta}):=\hat{\rho}(\xi) \hat{\alpha}^e(\boldsymbol{\eta}).$$
Notice that    the surface component  $\hat{\alpha}^e(\boldsymbol{\eta})$ varies over the partition $\mathcal{T}^S$, whilst the radial component $\hat{\rho}(\xi)$ is the same in all sectors $K^e$.  

Thus, we are in the following context of Duffy's approximation spaces
\begin{equation}
\mathcal{D}_k(S)=\left\{ \phi \in H^1(S);  \exists\; \hat{\phi} \in  \mathcal{D}_k(\hat{K})\; \mbox{such that} \; \phi|_{K^e} =\mathbb{F}_{K^e}(\hat{\phi}^e), \forall K^e \in \mathcal{T}^S\right\}, \label{eqDS}
\end{equation}
where $ \mathcal{\hat{D}}_k(\hat{K})$  is  a given   reference Duffy's approximation space   in the master element $\hat{K}$  described in Section \ref{sec:duffy}.   
For instance,  $\mathcal{D}_0(S)$  corresponds to the class of functions  in association with  $\hat{\phi}^e(\xi,\boldsymbol{\eta})= C \hat{\rho}(\xi), \forall K^e$,   obtained from  constant  trace functions $\alpha \equiv C \in \Lambda_0(\Gamma^S) $,  where $\Lambda_0(\Gamma^S)$ are the functions with constant value on $\Gamma^S$. It is clear that  $\mathcal{D}_0(S) \subset \mathcal{D}_k(S), \forall k \geq 0$. Particularly, let us also consider the subspace $\mathcal{D}^0_0(S)\subset \mathcal{D}_0(S)$ associated to radial  functions $\hat{\rho}(\xi)$ vanishing at $\xi=1$.
 
So far, $\mathcal{D}_k(S)$ is a functional space of infinite dimension, for  discretization only happens  for the surface component, living in a finite dimensional trace FE space $\Lambda_k(\Gamma^S)$, whilst the radial component can be chosen  arbitrarily. The SBFEM spaces   to be considered in Section \ref{SBFEM-S} are examples of  finite dimensional subspaces of $\mathcal{D}_k(S)$.  Other finite dimensional subspaces  $\mathcal{D}_{k,m}(S) \subset  \mathcal{D}_{k}(S)$ are also  of   interest:   functions  $\phi$  having local components   $\phi|_{K^e} =\mathbb{F}_{K^e}(\hat{\phi}^e)$,  where $ \hat{\phi}^e(\xi,\boldsymbol{\eta})=\hat{\rho}(\xi)\hat{\alpha}^e(\boldsymbol{\eta})$ with $\hat{\rho}\in \mathbb{P}_m[0,1]$ and $\hat{\alpha}^e\in V_k(\hat{L})$.

\subsubsection*{Gradient inner product  in  $\mathcal{D}_k(S)$ }\label{Inner}
Let  a pair of functions  $\phi, \psi \in \mathcal{D}_k(S)$ with local components   $\phi|_{K^e}= \mathbb{F}_{K^e}(\hat{\phi}^e)$, and  $\psi|_{K^e}=\mathbb{F}_{K^e}(\hat{\psi}^e)$,     $ \hat{\phi}^e(\xi,\boldsymbol{\eta})=\hat{\rho}(\xi)\hat{\alpha}^e(\boldsymbol{\eta})$,  and   $ \hat{\psi}^e(\xi,\boldsymbol{\eta})=\hat{\sigma}(\xi)\hat{\mu}^e(\boldsymbol{\eta})$ being associated  with  radial  $\hat{\rho}(\xi), \hat{\sigma}(\xi)$ and  surface $\hat{\alpha}(\boldsymbol{\eta}), \hat{\mu}(\boldsymbol{\eta})$ components.
Recalling the trace representation $\alpha(\mathbf{x}_b)  = \hat{\alpha}^e(\boldsymbol{\eta}) =\sum_l {\alpha}_l^e \hat{N}_k^l (\boldsymbol{\eta})$  for $\mathbf{x}_b \in L^e$, then formula  \eqref{eq:Grad2} becomes 
\[
\nabla_{\mathbf{x}}\phi(\mathbf{x})=\sum_{l}\alpha^{l,e}\left[\begin{array}{cc}
\underline{B}_{1l}^e(\boldsymbol{\eta}) & \underline{B}_{2l}^e(\boldsymbol{\eta})\end{array}\right]\left[\begin{array}{c}
\hat{\rho}'(\xi)\\
\frac{1}{\xi}\hat{\rho}(\xi)
\end{array}\right], \quad  \mbox{for} \; \mathbf{x} \in K^e,
\]
 both   $d \times 1$ matrices
\[
\underline{B}_{1l}^e(\boldsymbol{\eta})  =\mathbf{J}_{K^e}(1,\boldsymbol{\eta})^{-T}\left[\begin{array}{c}
\hat{N}_k^{l}(\boldsymbol{\eta})\\
0
\end{array}\right],  \quad 
\underline{B}_{2l}^e(\boldsymbol{\eta})  =\mathbf{J}_{K^e}(1,\boldsymbol{\eta})^{-T}\left[\begin{array}{c}
0\\
\nabla_{\boldsymbol{\eta}}\hat{N}_k^{l}(\boldsymbol{\eta})
\end{array}\right]
\] 
depending  on the geometry of the element at the boundary, and  on the surface component, but being independent of  the radial coordinate $\xi$
(see \cite{Song2018a} for the occurrence of these matrices in the  formulation of SBFEM methods). 
Analogous formula holds for $\psi$:
\[
\nabla_{\mathbf{x}}\psi(\mathbf{x})=\sum_{m}\mu^{m,e}\left[\begin{array}{cc}
\underline{B}_{1m}^e(\boldsymbol{\eta}) & \underline{B}_{2m}^e(\boldsymbol{\eta})\end{array}\right]\left[\begin{array}{c}
\hat{\sigma}'(\xi)\\
\frac{1}{\xi}\hat{\sigma}(\xi)
\end{array}\right].
\]
Thus, if $ \langle \phi, \psi\rangle_{\nabla,K^e}:= \int_{K^e}  \nabla_{\mathbf{x}} \phi(\mathbf{x}) \cdot \nabla_{\mathbf{x}} \psi(\mathbf{x}), \mbox{d}K^e$, then
{\small
 \begin{align}
 \langle \phi, \psi\rangle_{\nabla,K^e}& =\sum_{l,m}   \mu^{m,e} \, \alpha^{l,e} \int_{0}^{1} \int_{-1}^1  \begin{bmatrix} \vetor{B}_{1l}^e& \vetor{B}_{2l}^e\end{bmatrix} \begin{bmatrix}   \hat{\rho}'(\xi) \\  \frac{1}{\xi} \hat{\rho}(\xi) \end{bmatrix}  \cdot  \begin{bmatrix} \vetor{B}_{1m}^e& \vetor{B}_{2m}^e \end{bmatrix} \begin{bmatrix}   \hat{\sigma}^{'}(\xi) \\  \frac{1}{\xi}  \hat{\sigma}(\xi) \end{bmatrix}  \xi^{d-1} |\mathbf{J}_{K^e}(1,\boldsymbol{\eta})|\, d\boldsymbol{\eta} d\xi \nonumber\\
  & =\sum_{l,m}   \mu^{m,e} \, \alpha^{l,e}\int_{0}^1\int_{-1}^1\begin{bmatrix} \hat{\rho}'(\xi) &   \frac{1}{\xi} \hat{\rho}(\xi) \end{bmatrix} \cdot  \left(
     \begin{bmatrix} \vetor{B}_{1l}^{eT}\\ \vetor{B}_{2l}^{eT}\end{bmatrix} 
     \begin{bmatrix} \vetor{B}_{1m}^e& \vetor{B}_{2m}^e \end{bmatrix}
     |\mathbf{J}_{K^e}(1,\boldsymbol{\eta}) | \right) \begin{bmatrix}  \hat{\sigma}^{'}(\xi)\\   \frac{1}{\xi}\hat{\sigma}(\xi) \end{bmatrix}\,   \xi^{d-1} \,d\boldsymbol{\eta} d\xi  \nonumber\\  
   &  = \sum_{l,m}   \mu^{m,e} \, \alpha^{l,e}  \int_{0}^1 \begin{bmatrix} \hat{\rho}'(\xi) &   \frac{1}{\xi} \hat{\rho}(\xi) \end{bmatrix}  \tensor{E}^e_{ml}\begin{bmatrix}  \hat{\sigma}^{'}(\xi)\\   \frac{1}{\xi}\hat{\sigma}(\xi) \end{bmatrix}\,   \xi^{d-1}  d\xi,  \label{eq:InnerProdu}
  \end{align}}
  where the entries in   the    matrix $\tensor{E}_{ml}^e= \begin{bmatrix}  E_{11,ml}^e & E_{12,ml}^e \\ E_{21,ml}^e & E_{22,ml}^e\end{bmatrix}$
  are
\begin{align*}
    E_{11,ml}^e = \int_{-1}^1\vetor{B}_{1l}^{eT}(\boldsymbol{\eta})\vetor{B}^e_{1m}(\boldsymbol{\eta}) |\mathbf{J}_{K^e}(1,\boldsymbol{\eta})|\ \mbox{d} \boldsymbol{\eta}, &  \quad   E_{12,lm}^e  = \int_{-1}^1\vetor{B}_{1l}^{eT}(\boldsymbol{\eta})\vetor{B}^e_{2,m}(\boldsymbol{\eta})|\mathbf{J}_{K^e}(1,\boldsymbol{\eta})|\ \mbox{d}\boldsymbol{\eta}.\\
    E_{21,ml}^e = \int_{-1}^1\vetor{B}_{2l}^{eT}(\boldsymbol{\eta})\vetor{B}^e_{1m}(\boldsymbol{\eta}) |\mathbf{J}_{K^e}(1,\boldsymbol{\eta})|\ \mbox{d}\boldsymbol{\eta}, &   \quad   E_{22,lm}^e  = \int_{-1}^1\vetor{B}_{2l}^{eT}(\boldsymbol{\eta})\vetor{B}^e_{2,m}(\boldsymbol{\eta})|\mathbf{J}_{K^e}(1,\boldsymbol{\eta})|\ \mbox{d}\boldsymbol{\eta}.\
\end{align*}

\subsection{SBFEM spaces in $S$-elements}\label{SBFEM-S}
 There are two stages in the construction of local SBFEM approximation  spaces in $S$-elements, that we shall denote by $\mathbb{S}_k( S)$:   the restriction of a function in $\mathbb{S}_k( S)$ over the scaled boundary   $\Gamma^S$ is set in the  FE trace space $\Lambda_k (\Gamma^S)$, and  in the radial direction, it is  obtained analytically in terms of  eigenvectors and eigenfunctions of an ODE system, known SBFEM equation. 
Our purpose is to highlight  the  main aspects of SBFEM spaces in the context of Duffy's approximations    $\mathcal{D}_k(S)$ for $S$-elements, and to  show that   a paramount for the derivation of the SBFEM equation   is the enforcement of a gradient orthogonality constraint. 

Precisely,  having in mind that our goal is the solution of harmonic model problems, let us define the    subspace 
{\small
\begin{equation} \hspace{-.15cm}\mathbb{S}_k( S)=\left\{ \phi \in \mathcal{D}_k(S);  \langle \phi, \psi\rangle_{\nabla,S}:= \int_S \nabla_{\mathbf{x}} \phi(\mathbf{x}) \cdot \nabla_{\mathbf{x}} \psi(\mathbf{x}) \, dS =0, \;  \forall \psi \in \mathcal{D}^0_0(S), \,  \psi (\mathbf{O})=0\right\}. \label{SBFEM-space}\end{equation} }
 This definition suggests  that  the  functions $ \phi \in  \mathbb{S}_k( S) \subset  \mathcal{D}_k(S)$  have  boundary values  $\phi|_{\Gamma^S }= \alpha \in \Lambda_k (\Gamma^S)$, and  they are ``weak solutions'' of the harmonic equation $\Delta \Phi=0$ in $S$ with Dirichlet data $\alpha$.   Thus,  in some extent, $\mathbb{S}_k( S)$  can be interpreted as   ``radial harmonic extensions'' of the trace FE space $ \Lambda_k(\Gamma^S)$ to  the interior of $S$.

Notice that $\phi_0(\mathbf{x})\equiv 1$ is clearly in  $\mathbb{S}_k( S)$. The goal is to construct  linearly independent shape functions $\phi_i \in  \mathbb{S}_k( S)$ such that
\[\mathbb{S}_k( S)=  \mbox{span} \,\{\phi_i \} .\]
It is known that  the radial components and boundary values for the  SBFEM shape functions $\phi_i  $  are determined by a particular  family of  exact  eigenvalues and eigenfunctions solving  an  ODE system \cite{Song2018a}.  Next, we recover this representation of  $\phi_i$ using the current approach of Duffy's approximations constrained by the gradient orthogonality property expressed in \eqref{SBFEM-space}. 

Recall that, as a function in  $\mathcal{D}_k(S)$,  the shape function  $\phi_i \in \mathbb{S}_k( S)$ must be  obtained as  $\phi_i|_{K^e}=\mathbb{F}_{K^e}(\hat{\phi}_i^e)$, backtracking  a function $\hat{\phi}_i^e(\xi,\boldsymbol{\eta}) = \hat{\rho}_i(\xi)\hat{\alpha}_i^e(\boldsymbol{\eta}) \in \hat{\mathcal{D}}_k(\hat{K})$.  Moreover, we are assuming that   the local surface components $\hat{\alpha}_i^e(\boldsymbol{\eta})$ have  expressions   $\hat{\alpha}_i^e(\boldsymbol{\eta}) = \sum_l \alpha_{i}^{l,e} \hat{N}_k^l(\boldsymbol{\eta})$,  as linear combinations of  shape functions $\hat{N}_k^l(\boldsymbol{\eta}) \in {V}_k(\hat{L})$.   
 Thus,   it is necessary to  characterize  the radial functions $\hat{\rho}_i(\xi)$ and   the   multiplying coefficients  $\alpha_{i}^{l,e}$ allowing the verification of the gradient orthogonality property  stated in definition \eqref{SBFEM-space}.

\subsubsection*{Derivation of the SBFEM equation}\label{sec:ODE}
Let  $\psi (\mathbf{x}) \in \mathcal{D}_k(S)$ be a  general function locally defined as  $\psi|_{K^e}=\mathbb{F}_{K^e}(\hat{\psi}^e)$, where  $\hat{\psi}^e(\xi,\boldsymbol{\eta})=\hat{\sigma}(\xi)\hat{\mu}^e(\boldsymbol{\eta}) \in \hat{\mathcal{D}}_k(\hat{K})$ and  consider its  gradient inner product 
\[   \langle \phi_i , \psi\rangle_{\nabla,S} =\int_{S} \nabla_{\mathbf{x}} \phi_i(\mathbf{x})\cdot \nabla_{\mathbf{x}} \psi(\mathbf{x}) \,dS=\sum_e    \langle \phi_i , \psi\rangle_{\nabla,K^e},\] 
with a (searched)  shape function $\phi_i \in \mathbb{S}_k( S)$, where the   terms $  \langle \phi_i , \psi\rangle_{\nabla,K^e}$  are expressed  as in \eqref{eq:InnerProdu}. In fact, this formula can be rewritten as:
\begin{align}
  \langle \phi_i , \psi\rangle_{\nabla,K^e}= 
\sum_{m,l} \mu^{m,e} \, \alpha^{l,e} 
 \int_{0}^1 & \left( 
   \xi^{d-1}\hat{\rho}_i^{'}(\xi)E_{11,ml}^{e} \hat{\sigma}^{'}(\xi)+
   \xi^{d-2}\hat{\rho}_i^{'}(\xi)E_{12,ml}^e\hat{\sigma}(\xi)\right. \nonumber \\
 &  \left.   +  \xi^{d-2} \hat{\rho}_i(\xi)E_{21,ml}^e\hat{\sigma}^{'}(\xi)+
   \xi^{d-3}\hat{\rho}_i(\xi)E_{22,ml}^e\hat{\sigma}(\xi)
  \right) \, d\xi .\label{eq:SBFemInnerProd}
\end{align}
Let us denote  by $ \tensor{E}_{rs}$, $r,s \in\{1,2\}$,  the $\mathcal{N}^S\times \mathcal{N}^S$ matrices obtained by assembling the  matrices $\tensor{E}_{rs,ml}^e$,  element-by-element, according to the interelement connectivity. The process is similar to matrix assembly for FE discretizations of 
 boundary problems in $\mathbb{R}^{d-1}$.  
Moreover, consider  the vector functions   $\hat{\vetor{\Phi}}_i(\xi)= \hat{\rho}_i(\xi) \vetor{\alpha}_{i}$,  and   $\hat{\vetor{\Psi}}(\xi)= \hat{\sigma}(\xi) \vetor{\mu}$  collecting both radial and trace information of the shape functions $\phi_i(\mathbf{x})$ and of test functions  $\psi(\mathbf{x})$. Applying this notation, and summing up the contributions in \eqref{eq:SBFemInnerProd}, we obtain
\begin{align}
	 \langle \phi_i , \psi\rangle_{\nabla,S}  =   \int_{0}^1  \hat{\vetor{\Psi}}^{'}(\xi)  \cdot  \left[  \xi^{d-1}  \tensor{E}_{11}\hat{\vetor{\Phi}}_i^{'}(\xi) +  \xi^{d-2} \tensor{E}_{21} \hat{\tensor{\Phi}}_i(\xi)  \right] + \nonumber \\
    \hat{\vetor{\Psi}}(\xi) \cdot \left[  \xi^{d-2} \tensor{E}_{12}\hat{\vetor{\Phi}}_i^{'}(\xi)  +   \xi^{d-3} \tensor{E}_{22} \hat{\vetor{\Phi}}_i(\xi)\right]   d\xi. \label{eq:SBFemInnerProdTb}
\end{align} 
Consider $\hat{\vetor{Q}}_{i}(\xi)=\left[  \xi^{d-1} \tensor{E}_{11} \hat{\vetor{\Phi}}_i^{'}(\xi) +   \xi^{d-2}\tensor{E}_{21} \hat{\vetor{\Phi}}_i(\xi)  \right]$, and apply  integration by parts to obtain
\begin{align}  \int_{0}^1  \hat{\vetor{\Psi}}^{'}(\xi) \cdot  \hat{\vetor{Q}}_{i}(\xi)d\xi = \left.  \hat{\vetor{\Psi}} \cdot \  \hat{\vetor{Q}}_i\right]_0^1 -  \int_0^1 \hat{\vetor{\Psi}}(\xi) \cdot \hat{\vetor{Q}}_{i}^{'}(\xi) \, d\xi.  \label{eq:intp1}\end{align}
For $ \hat{\vetor{Q}}^{'}_{i}(\xi)=  \tensor{E}_{11}\left( \xi^{d-1}\hat{\vetor{\Phi}}^{''}_i (\xi ) + (d-1) \xi^{d-2}\hat{\vetor{\Phi}}^{'}_i (\xi )\right) + \tensor{E}_{21}\left(\xi^{d-2}\hat{\vetor{\Phi}}^{'}_i (\xi ) + (d-2)\xi^{d-3}\hat{\vetor{\Phi}}_i (\xi )\right)$, the inclusion of   formula  \eqref{eq:intp1} in \eqref{eq:SBFemInnerProdTb} gives
  \begin{align}
 \langle \phi_i , \psi\rangle_{\nabla,S}  =   \left.  \hat{\vetor{\Psi}} \cdot  \hat{\vetor{Q}}_i\right|_0^1  - \int_{0}^1  \hat{\vetor{\Psi}}(\xi)  \cdot  &\left[  \xi^{d-1} \tensor{E}_{11}\hat{\vetor{\Phi}}^{''}_i (\xi )   + \left[  (d-1) \tensor{E}_{11} - \tensor{E}_{12}+  \tensor{E}_{21}\right]   \xi^{d-2} \hat{\vetor{\Phi}}^{'}_i(\xi) \right. \nonumber\\ 
&   \left. +  \left[ (d-2) \tensor{E}_{21} - \tensor{E}_{22}\right]  \xi^{d-3}\hat{\vetor{\Phi}}_i(\xi) \right] \, d\xi.
  \label{eq:finalp} \end{align}

Recall that the purpose is to characterize the   functions ${\phi}_i (\mathbf{x})\in \mathcal{D}_k(S)$  such that the orthogonality property $M_{i}=0$ holds  for  all   functions  $\psi \in {\mathcal{D}}_0^0(S)$, i.e., vanishing  on $\Gamma^S$,  but  also vanishing on the scaling center. That is, for   $\hat{\sigma}(0)=\hat{\sigma}(1)=0$  and consequently
$\hat{\vetor{\Psi}}(0) = \hat{\vetor{\Psi}}(1)= 0$.  
These constraints on $\psi$ cancel  the boundary term   in \eqref{eq:finalp}. On the other hand, the condition for   vanishing the integral term in  \eqref{eq:finalp}  for all  $\hat{\vetor{\Psi}}(\xi)$  is equivalent to say  that $\hat{\vetor{\Phi}}(\xi)$ must solve the following equation  
  \begin{equation}
 {\small
  \xi^{d-1}\tensor{E}_{11}\hat{\vetor{\Phi}}^{''}_i (\xi )   + \left[  (d-1) \tensor{E}_{11} - \tensor{E}_{12}+  \tensor{E}_{21}\right]  \xi^{d-2}\hat{\vetor{\Phi}}^{'}_i(\xi)+  \left[  (d-2)\tensor{E}_{21} - \tensor{E}_{22}\right]  \xi^{d-3}\hat{\vetor{\Phi}}_i(\xi)}=0. \label{ODE}\end{equation}
Notice that this is the usual scaled boundary   equation  documented in \cite{Song2018a} for the SBFEM shape functions. The resolution  of \eqref{ODE} is well documented in the SBFEM literature, and it involves an  auxiliary eigenvalue problem for an ODE  system  in terms of  both $\hat{\vetor{\Phi}}_i(\xi)$ and   $\hat{\vetor{Q}}_{i}(\xi)$. For self completeness, the methodology is briefly described  in \ref{Riccati}.

In summary, the    resulting solutions $\hat{\vetor{\Phi}}_i = \hat{\rho}_i(\xi) \vetor{\alpha}_{i}$  have the form   $\hat{\rho}_i(\xi) =\xi^{\lambda_{i}}$, and $\vetor{\alpha}_i=	\vetor{A}_{i}$, 
where   ${\lambda_{i}}$  and $\vetor{A}_{i}$ refer to   positive real parts of the eigenvalues  and the associated eigenfunctions for the   ODE system equivalent to  the SBFEM equation \eqref{ODE}.  This information  is required for  the construction of the SBFEM basis functions $\phi_i$, giving the radial components  $\hat{\rho}_i(\xi)$  and  the trace surface components ${\alpha}_i  \in \Lambda_k(\Gamma^S)$ recovered from the  coefficient vectors $\vetor{\alpha}_i$. 
 Thus, the corresponding expressions  are
\begin{align} \phi_i(\mathbf{x})= {\hat{\phi}}_i^e(\xi,\boldsymbol{\eta}) &=\xi^{\lambda_{i}} \sum_{l} \vetor{A}_{i}^{l,e}\hat{N}_k^l(\boldsymbol{\eta}), \; \mbox{for}\;  \mathbf{x}=F_{K^e}(\xi,\boldsymbol{\eta}) \in K^e . \label{Phii}
\end{align}
 Analogously, associated to  $\hat{\vetor{Q}}_i(\xi)$ are  the flux functions
\begin{align} q_i(\mathbf{x})=\hat{{q}}_i^e(\xi,\boldsymbol{\eta})=\xi^{\lambda_{i}} \sum_{l} \vetor{Q}_{i}^{l,e}\hat{N}_k^l(\boldsymbol{\eta}), \; \mbox{for}\;  \mathbf{x}=F_{K^e}(\xi,\boldsymbol{\eta}) \in K^e. \label{i} \end{align}

\subsection{Orthogonality properties  of the SBFEM spaces}
\label{sec:developmentof}
In this section, we highlight two kinds of gradient orthogonality properties  held   by the SBFEM approximation spaces.

\subsubsection*{Intrinsic   gradient orthogonality property  for $\mathbb{S}_k(S)$}
The usual procedure for the  construction of SBFEM shape functions is the  determination of  analytic eigenfunctions for the  SBFEM equation  \eqref{ODE}. 
We have shown in the previous section  that there is another characterization of these shape functions that are  not  well recognized. Namely,    implicit in the  condition for a function $\phi \in \mathcal{D}_k(S) $ to solve the SBFEM equation   \eqref{ODE} is  the gradient orthogonality property, enforced from the start, in the definition of the subspaces $\mathbb{S}_k(S) $ in \eqref{SBFEM-space}. Precisely,   a function $\phi \in \mathbb{S}_k(S) \subset  \mathcal{D}_k(S) $ if  the   gradient orthogonality constraint 
\begin{equation}\langle \phi, \psi \rangle_{\nabla,S}= \int_S \nabla_{\mathbf{x}} \phi(\mathbf{x}) \cdot \nabla_{\mathbf{x}} \psi(\mathbf{x}) \,dS=0 \label{eq:orto1} \end{equation}
holds for all  $\psi  \in \mathcal{D}_0^0(S)$, with   $\psi (\mathbf{O})=0$. In such case, then $\phi$ solves equation \eqref{ODE}.

\subsubsection*{Extended   gradient orthogonality property for  $\mathbb{S}_k(S)$}
 Let $\mathcal{H}(S)$ denotes the space of  harmonic functions in $S$. Then, it is clear that $\langle \phi, \psi \rangle_{\nabla,S}=0$ for all $\phi \in \mathcal{H}(S)$ and $\psi \in H^1_0(S)$, giving  the  well-known decomposition
  \begin{equation} H^1(S) = \mathcal{H}(S) \overset{\nabla}{\oplus}  H^1_0(S),  \label{harmonic_dec} \end{equation}
  where the symbol $\overset{\nabla}{\oplus}$ denotes the orthogonality  relation  with respect to the gradient inner product $\langle \cdot, \cdot \rangle_{\nabla,S}$. Our purpose is to show a similar   relation  for Duffy's spaces $ \mathcal{D}_k(S) \subset H^1(S) $,  $\mathbb{S}_k(S)$ playing the role of the harmonic functions. For that, we need to extend the gradient orthogonality property \eqref{eq:orto1} to functions  $\psi \in \mathcal{D}_0(S)$.

 \begin{proposition}\label{ExtendedOrtho} The orthogonality property
  \begin{equation}	\langle \phi, \psi \rangle_{\nabla,S}=0, \quad \forall \phi \in \mathbb{S}_k(S) \; \mbox{ and}\;    \psi  \in \mathcal{D}_0(S) \label{eq:orto2}\end{equation}
 is valid. Thus,  
\begin{equation}\mathcal{D}_k(S) = \mathbb{S}_k(S)  \overset{\nabla}{\oplus} \mathcal{D}_0(S) \label{harmonic_dec_discr} \end{equation}
holds as  a mimetic version of \eqref{harmonic_dec}.

\end{proposition}

\begin{proof} A crucial step in the derivation of   the  SBFEM equation \eqref{ODE}  is the formula   for the gradient inner   product $ \langle \phi, \psi \rangle_{\nabla,S}$ given in  \eqref{eq:finalp},  where two terms enter into play: a boundary term and an integral term. The constraints $ \psi  \in \mathcal{D}_0^0(S)$ and  $\psi (\mathbf{O})=0$ make the boundary term to be zero, and \eqref{ODE}  derives from the assumption \eqref{eq:orto1}.  

Now let us  relax the constraints $\psi  \in \mathcal{D}_0^0(S)$ and   $\psi (\mathbf{O})=0$. Instead,   take  $\psi$ in a broader space  $\mathcal{D}_0(S)$. Clearly, the property  $	\langle \phi, \psi \rangle_{\nabla,S}=0$  holds for $\phi=\phi_0\equiv1$. Thus, it is sufficient to verify it for all  shape functions $\phi=\phi_i$ associated with eigenvalues $\lambda_{i}\neq 0$.  

Notice that the desired orthogonality  property  \eqref{eq:orto2}  is valid for $\psi  \in\mathcal{D}_0(S)$,  with    $\hat{\psi}^e(\xi,\boldsymbol{\eta})= C \hat{\sigma}(\xi)$  in the sectors $K^e$,  if and only if it holds for functions $\varphi=\psi - C\hat{\sigma}(1) \in  \mathcal{D}^0_0(S)$, i.e., for the cases where $ \hat{\varphi}(\xi)=  C (\hat{\sigma}(\xi)-  \hat{\sigma}(1))$, with $\hat{\varphi}(1)=0$.
For them, 
we apply  equation \eqref{ODE},   valid for all shape functions $\phi_i \in \mathbb{S}_k(S)$,  to reduce the equation   \eqref{eq:finalp}  to
 $$\langle \phi_i, \varphi \rangle_{\nabla,S} = \hat{\varphi}(1)  \sum_{n}\hat{\vetor{Q}}_{in}(1) -  \hat{\varphi}(0)  \sum_{n}\hat{\vetor{Q}}_{in}(0).  $$ 
Thus,  since  $\hat{\varphi}(1) =0$ and $\hat{\vetor{Q}}_{in}(0)=0$,  we obtain the orthogonality property \eqref{eq:orto2}.
\end{proof}

\section{Interpolants}\label{Interpolants}
When a Galerkin method is used to approximate a boundary value
problem, one of the most important choices is the family
of approximation spaces. For elliptic problems the achievable error
of approximation is equal to the
error obtained by approximating the solution of the partial
differential equation directly from the trial space.
 The accuracy  is accessed {\it a priori} by  bounds computed in terms of  interpolant errors using the approximation space.  In  the context of piecewise defined approximations over subregions (elements)  of the computational domain, as is the case of  FE methods,  the interpolants usually show the following characteristics: 
 \begin{itemize} 
\item Locality:  in each subregion,  a polynomial trace interpolant  over the  boundary is extended   to the interior  (a process  also called   lifting). 
\item  Global conformity: it follows  directly from the hypothesis that the trace interpolants depend exclusively on the function restriction over subregion  boundaries.
\item  Optimality: optimal interpolation error estimates  are achieved with respect to  the discretization parameters: mesh width  and polynomial order. 
\end{itemize}

In this direction,   the plan  is to construct    interpolants  in   SBFEM trial  spaces,  and to explore them to evaluate  the potential of SBFEM approximations. Firstly, let us introduce some new notation and auxiliary results already known in other contexts. 

Consider a family  of  conformal  polytopal   partitions $\mathcal{T}^h= \{S\}$     of $\Omega$  by $S$-elements, as described in Section \ref{sec:Selements}.  Define  the  mesh skeleton $\Gamma^h = \cup_{L\in \mathcal{E}^h} L$ by the  assembly  of all facets  (edges of faces) in $ \mathcal{E}^h=\{ L \subset \Gamma^{h,S},  S\in  \mathcal{T}^h\}$. The parameter  $h$ refers to  the  characteristic size of  the facets in $\Gamma^h$.  Moreover, define   the    conglomerate partitions $\mathcal{P}^h= \cup_{S\in \mathcal{T}^h} \mathcal{T}^{h,S}$  of $\Omega$. Recall that the elements $K \in \mathcal{T}^{h,S}$ may be affine triangles, pyramids, or tetrahedra inheriting the  conformal property from $\mathcal{T}^h$. In principle,  shape regularity of  $\mathcal{P}^h$ is  not a  granted property.

Based on the partitions $\Gamma^h$, $\mathcal{T}^h$ or  $\mathcal{P}^h$,  we  consider the following approximation spaces.
\begin{itemize}
\item   FE trace spaces:  $\Lambda_k(\Gamma^h) = C^0(\Gamma^h)\cap \prod_{L \in  \mathcal{E}^h} V_k(L)$, piecewise polynomial spaces,  where    $V_k(L)=\mathbb{P}_{k}({L})$, for 1D edges  and  triangular facets $L$, 
and $V_k({L})=\mathbb{Q}_{k,k}({L})$,  for  quadrilateral facets  ${L}$.

\item Duffy's spaces  $\mathcal{D}_{k}^h \subset H^1(\Omega)$:   given  the local Duffy's spaces $ \mathcal{D}_{k}^h(S),  S\in \mathcal{T}^h$ defined in Section \ref{sec:sbfemin}, set 
 \begin{align*}\mathcal{D}_{k}^h & =  \{w \in H^1(\Omega); w|_S\in  \mathcal{D}_{k}^h(S), S\in \mathcal{T}^h\},\\
 \mathcal{D}^{0,h}_0 & =\{w \in H^1(\Omega);  w|_{S}\in  \mathcal{D}^{0}_0(S), S\in\mathcal{T}^h\}. \end{align*}
Notice that   $\mathcal{D}^{0,h}_0\subset \mathcal{D}_{k}^h, \forall k \geq 0$. 

\item SBFEM  spaces $\mathbb{S}_k^h \subset H^1(\Omega)$: given  local SBFEM spaces $\mathbb{S}_k^h(S)\subset  \mathcal{D}_{k}^h(S), S\in \mathcal{T}^h$,   described in  Section \ref{SBFEM-S},  define
\[\mathbb{S}_k^h =\{w \in H^1(\Omega); w|_S\in \mathbb{S}_k^h(S), S\in \mathcal{T}^h\},\]
and set  $\mathbb{S}_{k,0}^{h}= \mathbb{S}_k^h\cap H^1_0(\Omega)$.

\item FE spaces $\mathcal{V}^{h,FE}_k \subset H^1(\Omega)$:  Consider the following FE spaces  based on  the conglomerated meshes $\mathcal{P}^h$.
\begin{enumerate} 
\item  Triangular  (2D) and tetrahedral  (3D) meshes $\mathcal{P}^h$:  $\mathcal{V}^{h,FE}_k:= \mathbb{P}_k(\mathcal{P}^h) \cap H^1(\Omega)$, where $\mathbb{P}_k(\mathcal{P}^h)$ stands for functions piecewise defined by polynomials in  $\mathbb{P}_k(K)$, $K \in \mathcal{P}^h$,  of degree not greater than $k$.  
\item Pyramidal (3D) meshes  $\mathcal{P}^h$:  let us consider $\mathcal{V}^{h,FE}_k : = \mathcal{U}^{(0),k}(\mathcal{T}^h) \cap H^1(\Omega)$, piecewise defined by a class of rational   polynomials $\mathcal{U}^{(0),k}(K)$,  for     $K \in \mathcal{P}^h$ \cite{Nigam2012}.  Traces  of functions in $\mathcal{U}^{(0),k}(K)$  are in  $\mathbb{P}_k(L)$ for  triangular faces, and  in $\mathbb{Q}_{k,k}(L)$  if $L$ is  quadrilateral.  Moreover, $\mathbb{P}_k(K) \subset \mathcal{U}^{(0),k}(K)$.  
\end{enumerate}
\end{itemize}
 \begin{proposition} \label{embbeding} (i)  For  $w\in \mathcal{V}_k^{h,FE}$, $w|_{\Gamma^h}\in \Lambda_k(\Gamma^h)$.  (ii)   $ \mathbb{P}_k(\mathcal{P}^h)  \subset \mathcal{V}^{h,FE}_k  \subset \mathcal{D}_{k}^h$.  \end{proposition}

\begin{proof}   The trace property (i)  and the  polynomial inclusion in  (ii)  are already known. To proof the  second embedding property  in (ii), let us start by considering   three particular collapsed triangular, pyramidal and tetrahedral reference  elements.  
  \begin{itemize}
  \item {\it A triangular reference element $K$}:  Let      $K $ be the reference triangle,  with collapsed vertex  $\mathbf{a}_0=(0,0)$,  and the opposed edge $L=[\mathbf{a}_1, \mathbf{a}_2]$, where  $\mathbf{a}_1=(1,0)$ and $\mathbf{a}_2=(1,1)$. Taking the  mapping  ${F}_{L}: \hat{L} \rightarrow L$, defined as   ${F}_{L}(\boldsymbol{\eta})= \left(\frac{1+\boldsymbol{\eta}}{2},  \frac{1-\boldsymbol{\eta}}{2}\right)$, 
 the   Duffy's transformation  from $\hat{K}$  over $K$ becomes
$  x = \frac{\xi}{2}(1+\boldsymbol{\eta}), \;
   y = \frac{\xi}{2}(1-\boldsymbol{\eta})$,
whose  inversion is
$ \xi = x+y,  \; \boldsymbol{\eta} =  \frac{x-y}{x+y}$.  
Let    $\psi \in \mathcal{D}_k(K)$  be the  pullback of  functions $\mathbb{F}_{K}(\hat{\psi})  \in \mathcal{D}_k(\hat{K})$, where $ \hat{\psi}(\xi,\boldsymbol{\eta})=\xi_k \hat{\alpha}(\boldsymbol{\eta})$, so that $\psi(x,y)=(x+y)_k \hat{\alpha}(\frac{x-y}{x+y})$. Thus, by varying  $\hat{\alpha}\in \mathbb{P}_k(\hat{L})$,  we conclude that all functions $\psi(x,y) \in \mathbb{P}_k(K)$ can be recovered in  $\mathcal{D}_k(K)$.
 \item {\it A pyramidal reference element}: Suppose $K$ is a  pyramid with  vertex $\mathbf{a}_0=(0,0,1)$, and opposed face $L=[\mathbf{a}_1, \mathbf{a}_2,\mathbf{a}_3, \mathbf{a}_4]$, with  vertices $\mathbf{a}_1=(0,0,0)$,  $\mathbf{a}_2=(1,0,0)$, $\mathbf{a}_3=(1,1,0)$, and $\mathbf{a}_4=(0,1,0)$.  
The  FE space $\mathcal{U}^{(0),k}(K) \subset H^1(K)$   proposed in \cite{Nigam2012} is the  first  space of an exact sequence $\mathcal{U}^{(s),k}(K)$ verifying the De Rham commuting property.  Their definition considers  the  the geometric transformation  $S_\infty: K_{\infty} \rightarrow K$ of the  ''infinite pyramid'' $K_{\infty}=\{(x,y,z)\in \mathbb{R}^3; x,y,z \geq 0, x\leq 1, y\leq 1\}\cup\{\infty\}$,    given by $S_\infty(x,y,z)=\left(\frac{x}{1+z}, \frac{y}{1+z}, \frac{z}{1+z}\right)$, $S_\infty(\infty)=\mathbf{a}_0$. The functions $ w \in \mathcal{U}^{(0),k}(K)$  are obtained by the pullback  $\mathbb{S}_\infty(u)$ of  functions $u$ in a properly chosen subspace of the  rational functions 
  $Q_{k}^{k,k,k}(\hat{K}_{\infty})=\{ \frac{q}{1+z}; q \in \mathbb{Q}_{k,k,k}(K_{\infty})\}$.  Our goal is to show that  $ \mathcal{U}^{(0),k}(K)$ can also be interpreted in the context of the Duffy's space  $\mathcal{D}_{k}(K)$. For that, consider the  hexahedron $H=[0,1] \times[0,1] \times [0,1]$, with the coordinate system $(\mu_1,\mu_2,\xi)$,  with  $(\mu_1,\mu_2)\in [0,1]\times[0,1]$ and $0\leq \xi \leq 1$.  Observe that the  geometric transformation $F_\infty: H \rightarrow K_\infty$, $F_\infty (\mu_1,\mu_2,\xi)= (\mu_1,\mu_2,\frac{\xi}{1-\xi})$ collapses the face $\xi=1$ in $H$ onto $\infty$. Moreover,   $Q_{k}^{k,k,k}(\hat{K}_{\infty})=\mathbb{F}_\infty(\mathbb{Q}_{k,k,k}(H))$.  Consequently, 
  \begin{equation}\mathcal{U}^{(0),k}(K) \subset \mathbb{S}_\infty(Q_{k}^{k,k,k}(\hat{K}_{\infty}))= \mathbb{S}_\infty(\mathbb{F}_\infty (\mathbb{Q}_{k,k,k}(H)).\label{DuffyPyramid}\end{equation} 
 On the other hand,  the transformation $F_K:H \rightarrow K$, defined by the composition $\mathbf{x}=F_K(\eta_1,\eta_2,\xi)=S_\infty (F_\infty(\eta_1,\eta_2,\xi))$  results to  be a Duffy's transformation collapsing the face $\xi=1$ in $H$  on top of the vertex  $\mathbf{a}_0 \in K$. Consequently,      $ \mathbb{F}_K (\mathbb{Q}_{k,k,k}(H)) \subset \mathcal{D}_k(K)$.  Thus, using \eqref{DuffyPyramid},  we obtain $\mathcal{U}^{(0),k}(K) \subset \mathcal{D}_{k}(K)$. 
    \item {\it A tetrahedral  reference element}: Suppose $K$ is the  reference tetrahedron with collapsed vertex $\mathbf{a}_0=(0,0,0)$,   and opposed quadrilateral face $L=[\mathbf{a}_1, \mathbf{a}_2,\mathbf{a}_3]$, with $\mathbf{a}_1=(1,0,1)$,  $\mathbf{a}_2=(1,0,0)$ and  $\mathbf{a}_3=(1,1,0)$. Notice that $L$ can be mapped by $\mathbf{x}= F_L(\boldsymbol{\eta})$, where 
$  x= 1-\eta_1-\eta_2$, $y= \eta_1$, and $z= \eta_2$.  Then, the Duffy's transformation is $F_K(\xi,\boldsymbol{\eta})= \xi F_L(\boldsymbol{\eta})$, whose inverse is 
    $\xi=x+y+z, \quad \eta_1=\frac{y}{x+y+z}, \quad \eta_2=\frac{z}{x+y+z}$
      Let   $\psi = \mathbb{F}_{K} (\hat{\psi}) \in \mathcal{D}_{k}(K)$, with $ \hat{\psi}(\xi,\boldsymbol{\eta})=\xi_k \hat{\alpha}(\boldsymbol{\eta})$,  and  $\hat{\alpha}\in \mathbb{P}_{k}(\hat{L})$. Thus,  the functions $\psi(x,y,z)=(x+y+z)_k \hat{\alpha}(\frac{y}{x+y+z}, \frac{z}{x+y+z})$ recover all functions  in $ \mathbb{P}_k(K)$.
 \end{itemize}

 Now consider  a general element $K^e  \in \mathcal{P}^h$, with collapsed vertex $\mathbf{O}$, and opposed face $L^e$ with vertices $\mathbf{a}_{l}^e$.  Notice that  $K^e$ can be seen as a geometric affine transformation of one of the reference elements $K$ described above, i.e., $K^{e} = {T}^e(K)$,  such that $\mathbf{O}={T}^e(\mathbf{a}_0)$, $\mathbf{a}_l^e={T}^e(\mathbf{a}_l)$, and thus $L^{e}={T}^e(L)$).  Since  the polynomials $\mathbb{P}_k(K)$,  for triangles and tetrahedra,  and rational polynomials $\mathbb{S}_k(K)$), for pyramids, are preserved by affine transformations, then we  conclude  that   $\mathcal{V}^{h,FE}_k  \subset \mathcal{D}_{k}^h$.  \end{proof} 

\subsection{ FE  interpolants}
Interpolant operators  $\mathcal{F}^{h,FE}_k: H^{s}(\Omega) \rightarrow   \mathcal{V}^{h,FE}_k $ have being designed as useful tools  for functions in general Sobolev spaces $H^s(\Omega), s\geq1$. As already mentioned, they are constructed  by first   defining a piecewise polynomial trace interpolant  over the  facets $L\subset \partial K$  of each element $K\in \mathcal{P}^h$, and then by extending  this trace interpolant to the interior of $K$.   Let us recall  some examples and error estimates  already available in the literature. For them, we assume the affine conglomerate triangular, pyramidal or tetrahedral partitions $\mathcal{P}^h$ are   regular (e.g. quasi-uniform and shape regular, with parameters independent of $h$). Under these circumstances, the following estimates hold.
  \begin{itemize}
   \item  There are  interpolands $\mathcal{F}^{h,FE}_k \, w$    over FE spaces $\mathcal{V}^{h,FE}_k = \mathbb{P}_k(\mathcal{P}^h)\cap H^1(\Omega)$  defined in   \cite{Babuska1987}   for  triangles and in \cite{Munoz-Sola1997} for thetrahedra.     Suppose  $w\in H^{s}(\Omega)$, $s>\frac{3}{2}$ in 2D, and $s>2$ in 3D, then the estimate
    $$|w-\mathcal{F}^{h,FE}_k \, w|_{H^{1}(\Omega)}\lesssim \frac{h^{\mu-1}}{k^{d-2}}\quad\|w\|_{H^{s}(\Omega)} $$ 
holds    for $\mu=\min(k+1,s)$, where the leading constant on the right side is independent
of $w$, $h$, and  $k$ (but depends on  $s$ and  regularity parameters of $\mathcal{P}^h$).
\item There are also the  projection-based interpolants,  proposed  by L. Demkowicz and coworkers, as expounded in  \cite{Demkowicz2008,Demkowicz2007}. They admit a general form, without requiring any specific geometric aspect, and have the flexibility to treat general local spaces, not necessarily polynomials.  Note that  such  constructions may require 
 additional regularity assumptions beyond the minimal $H^1$-conformity. Indeed,  the trace interpolants may require interpolation at   element vertices,  requiring the
 regularity $H^{1+s}$ with $s > 1/2$ in 3D FE settings.
For  FE spaces $\mathcal{V}^{h,FE}_k = \mathbb{P}_k(\mathcal{P}^h)\cap H^1(\Omega)$ based on tetrahedra,  the error estimates  stated in \cite[Theorem 2.2]{Demkowicz2007} for   projection based-interpolants  $\mathcal{F}^{h,FE}_k \, w$  have the  non-optimal form    
 \begin{align} 
   |w - \mathcal{F}^{h,FE}_k \, w|_{H^1(\Omega)} &\lesssim (\ln k)^2 \, \left(\frac{h}{k}\right)^{d-2}  |w|_{H^{s}}, \;s >3/2. \label{eq:ProjFk3D} 
 \end{align}
 The suboptimal logarithmic factor  appearing in \eqref{eq:ProjFk3D} can be dropped in the $k$-version under the more stringent regularity assumption   
 $s\geq 2$ \cite[Corollary 2.12]{Melenk2020}. 

\item For pyramidal partitions  $\mathcal{P}^h$,  projection-based interpolants  $\mathcal{F}^{h,FE}_k \, w $  over the FE spaces  $\mathcal{U}^{(0),k}(\mathcal{T}^h)$ are defined in \cite{Nigam2012}. However, to the best of our knowledge,  error estimates are still missing for them, but optimal $h$-convergence rates   have been observed in  numerical experiments presented in \cite{Bergot2010}.
 \end{itemize}

  \subsection{SBFEM interpolant}\label{interpolantSBFEM}
As for the cases of FE  spaces, we construct  interpolant operators   $\Pi_k^h:  H^{s}(\Omega)  \rightarrow  \mathbb{S}_k^h$, for sufficiently smooth functions $w\in H^{s}(\Omega)$,  following three steps:   a trace interpolant  $\mathcal{I}_k^h: H^{s}(\Gamma^h)  \rightarrow  \Lambda_k(\Gamma^h)$,  local projections  $\Pi^{h,S}_k:  H^{s}(S)  \rightarrow  \mathbb{S}_k(S)$ extending  trace functions to the interior of the element, and  assembly of local interpolants. 
\begin{enumerate}
\item Trace interpolant    $\mathcal{I}_k^h:H^{s}(\Omega) \rightarrow  \Lambda_k(\Gamma^h)$  -  it is piecewise  defined on the facets $L \in \mathcal{E}^h$,   following any of the interpolation strategies used so far for  the  FE spaces   $\mathcal{V}^{h,FE}_k$.
\item   Local projections  $ \Pi^{h,S}_k:  H^{k+1}(S)  \rightarrow  \mathbb{S}_k(S)$:    $ \Pi^{h,S}_k w \in  \mathcal{D}^{h,S}_k$  solves  the problem
\begin{align}
 \langle \Pi^{h,S}_k\, w, v\rangle_{\nabla,S}  &= 0 \quad \forall v  \in {\mathcal{D}}_{0}^{0}(S), \label{eq:LocalD}\\
 \Pi^{h,S}_k\, w|_{\Gamma^S} &=  \mathcal{I}_k^h\, w|_{\Gamma^S}. \label{eq:CC}
\end{align}
Notice that equation  \eqref{eq:LocalD} ensures that  $\Pi^{h,S}_k\, w \in \mathbb{S}_k(S)$ and the relation \eqref{eq:CC}  enforces the trace  constraint  matching  $\Pi^{h,S}_k\, w$ to the trace interpolant of $w$. It is clear from  these equations  the interpretation  of  $ \Pi^{h,S}_k$ as  "radial harmonic extension" of the trace interpolant $\mathcal{I}_k^h \, w$ to the interior of $S$.    Let      $\vetor{\omega}^{S}$ be  the coefficients in the expansion $\mathcal{I}_k^h \, w (\mathbf{x}_b)=\sum_{n=1}^{\mathcal{N}^S}  \omega^{n,S} N_k^{n,S}(\mathbf{x}_b)$, $\mathbf{x}_b\in \Gamma^S$.    We seek for coefficients $\vetor{c} =[c_i]$ such that   $ \Pi^{h,S}_k \, w   = \sum_i c_i \phi_i^S \in  \mathbb{S}_k^h( S)$. 
  According to the definition of the local spaces $\mathbb{S}_k^h(S)$,   the solution  is $\vetor{c}= \vetor{{\omega}}^{S} \tensor{A}^{-1}$, where $\tensor{A}=\tensor{A}^S$ is the eigenvector matrix associated to the traces of the  SBFEM shape-functions  $\phi^S_i$ over $\Gamma^S$. 
\item Assembly - Define $\Pi_k^h \,w$   by assembling  the  local contributions $\Pi_k^h \,w|_{S}=\Pi^{h,S}_k\,w$.  It is clear that $\Pi^{h,S}_k \,w|_{L}=\Pi^{h,S}_k \,w|_{L}$ over an interface $L=S \cap S'$ shared by two S-elements. Thus,  the  conformity property $\Pi_k^h \,w \in H^1(\Omega)$ holds.
\end{enumerate}

\subsubsection*{Remarks}
\begin{itemize}
\item[(1)]  In the same manner as  FE interpolants $\mathcal{F}^{h,FE}_k \,w$, the SBFEM interpolant   $\Pi_k^h$ satisfies the two fundamental properties:   locality and global conformity. However, they differ on the way the trace interpolant is extended to   the interior of the $S$-elements by their local projections. Recall that the  ''radial harmonic extension''   adopted in the SBFEM context is possible due to   the  particular scaled geometry of the $S$-elements. Moreover, when the SBFEM interpolant shares the trace interpolant  of  $\mathcal{F}^{h,FE}_k \,w$,  then it is clear that
\begin{equation} \Pi _k^h\, w =  \Pi _k^h\,  \mathcal{F}^{h,FE}_k \,w. \label{proj_property}\end{equation}

 \item[(2)] Since     $\mathcal{I}_k^h\, w=w|_{\Gamma^h}$ for functions   $w\in \mathcal{D}_k^{h}$, the trace constraint   \eqref{eq:CC} means that $w- \Pi_{k}^h\, w\in  {\mathcal{D}}^{0,h}_0$ for all functions $w$ in the Duffy's space $\mathcal{D}_k^{h}$. 
Consequently, Proposition \ref{ExtendedOrtho}  implies  the  orthogonality property 
  \begin{equation}  \langle w -\Pi^{h}_k\, w,  v \rangle_\nabla =\sum_{S\in \mathcal{T}^h}\langle w -\Pi^{h,S}_k\, w,  v \rangle_{\nabla,S}=0,  \quad \forall  w\in \mathcal{D}_k^{h}, \quad \forall v \in \mathbb{S}_k^h. \label{orto2D}\end{equation}
    \end{itemize}

\subsection{Comments on  the SBFEM interpolation errors}\label{Comments}

Unlike general purpose FE techniques, SBFEM  approximations are constructed to be applied for a specific type of problem. Thus, for the model Laplace  problem under consideration,  there is no interest in accessing the accuracy of  SBFEM interpolants  $\Pi _k^h \, w$  when applied to  other than   for  harmonic functions $w\in \mathcal{H}(\Omega)$.  
 For them,  the sources of SBFEM   interpolation  errors  are  two-fold:
 \begin{itemize}
\item[ (i)]  the polynomial discretization of traces $w|_{\Gamma^h} \approx  \mathcal{I}_k^h\, w \in \Lambda_k^h$.
\item[(ii)] the deviation of $\Pi^{h,S}_k\, w \in \mathbb{S}_k^h$ of being an  harmonic function.  
\end{itemize}
In this direction, let us consider the subspaces 
 $$\mathcal{V}^{h,\Delta}_k =  \{w \in \mathcal{H}(\Omega);  w|_{\Gamma^h} \in \Lambda_k(\Gamma^h)\},$$
where only trace discretization takes place. Denoted by harmonic virtual spaces, they have been used   in the context of the operator adapted virtual FE method proposed in  \cite{Chernov2019}, and   designed to solve two-dimensional harmonic problems. The term  “virtual” emphasizes that  functions in $\mathcal{V}^{h,\Delta}_k$ are not known explicitly in the interior of each subregion $S \in \mathcal{T}^h$.

The finite-dimensional spaces $\mathcal{V}^{h,\Delta}_k$   have close similarities with the SBFEM spaces  $\mathbb{S}_k^h$. In both cases, the   trace functions are in $\Lambda_k(\Gamma^h)$, which are extended to the interior of the $S$-elements by solving local Dirichlet Laplace problems: whilst the  functions in the local spaces  $V^{\Delta}_k(S)=\mathcal{V}^{h,\Delta}_k|_S$ are strongly harmonic in $S$, the ones in   $\mathbb{S}_k(S)$ are  harmonic in a weaker sense. However, unlike for the harmonic subspaces  $V^{\Delta}_k(S)$,  it is possible   to explore  the  radial Duffy's structure of  $\mathbb{S}_k^h(S)$  to  explicitly compute shape functions for them, as described in the previous section.

Let us consider  the harmonic virtual interpolant $\mathcal{F}^{h,\Delta}_k: H^{s}(\Omega) \rightarrow   \mathcal{V}^{h,\Delta }_k$ by  solving  the  local Laplace problems
\begin{align}
 \langle \mathcal{F}^{h,\Delta}_k \,w, v\rangle_{\nabla,S}  &= 0 \quad \forall v  \in H_{0}^1(S), \label{eq:LocalH}\\
\mathcal{F}^{h,\Delta}_k\, w|_{\Gamma^S} &=  \mathcal{I}_k^h\, w|_{\Gamma^S}, \label{eq:CCH}
\end{align}
where the  trace interpolant $ \mathcal{I}_k^h\, w$ is the one adopted in $ \Pi _k^h \, w$. Note that this is an analytic recovery problem for it is not directly accessible for computation, whilst the SBFEM interpolant $ \Pi _k^h\, w$ is a   computable recovery problem.  

For an harmonic  function  $u \in \mathcal{H}(\Omega)$,  let us consider the decomposition
  \begin{align} u -\Pi^{h}_k\, u =  (u- \mathcal{F}^{h,\Delta}_k \,u )+ (\mathcal{F}^{h,\Delta}_k \,u -  \Pi _k^h \, u) = (i) + (ii).  \label{dec_virtual} \end{align}
The first term $(i)= u- \mathcal{F}^{h,\Delta}_k \,u$ compares two harmonic functions  differing on the skeleton $\Gamma^h$  by the trace interpolation error $u- \mathcal{I}_k^h\, u$, meaning that only the interface errors require to be estimated. 
In fact, the application of  Neumann trace inequality (\cite[Theorem A.33]{Schwab1998})  in  each $S$-element  $S\in \mathcal{T}^h$ gives
\begin{align}
|u-\mathcal{F}_{k}^{h,\Delta}\,u|_{H^{1}(S)} \lesssim \|u-\mathcal{I}_{k}^{h}u\|_{H^{\frac{1}{2}}(\partial S)}. \label{traceerror}
\end{align}
We refer to \cite[Lemma 4.4, Lemma 4.5]{Chernov2019} for estimates of  \eqref{traceerror} in the particular Gauss-Lobatto  trace interpolation case, and  under some specific graded polygonal mesh circumstances.
On the other hand, since 
  \begin{equation} \Pi _k^h\, u =  \Pi _k^h\,  \mathcal{F}^{h,FE}_k \,u, \label{proj_propertyV}\end{equation}
 the second term becomes   $(ii) =  \mathcal{F}^{h,\Delta}_k \,u -  \Pi _k^h \, u=   \mathcal{F}^{h,\Delta}_k \,u -  \Pi _k^h \,  \mathcal{F}^{h,\Delta}_k \,u$, representing the SBFEM interpolation error for the harmonic virtual function  $\mathcal{F}^{h,\Delta}_k \,u \in \mathcal{V}_k^{h,\Delta}$.  Consequently,  according to  \eqref{eq:LocalD}  and \eqref{eq:LocalH},  we obtain
  \begin{align} \langle \mathcal{F}^{h,\Delta}_k \,u -  \Pi _k^h \,  \mathcal{F}^{h,\Delta}_k \,u, v\rangle_\nabla & 
  =\sum_{S \in \mathcal{T}^h}\langle \mathcal{F}^{h,\Delta}_k \,u -  \Pi _k^h \,  \mathcal{F}^{h,\Delta}_k \,u,  v \rangle_{\nabla,S}=0, \quad \forall v  \in {\mathcal{D}}_{0}^{0,h}.  \label{eq:orto3}\end{align}
 In other words, the second term $(ii)=\mathcal{F}^{h,\Delta}_k \,u -  \Pi _k^h \, u$, which  vanishes in $\Gamma^h$,    is orthogonal to $ {\mathcal{D}}_{0}^{0,h}$ with respect to the gradient inner product. Thus its   energy norm  is a measure of the deviation of  $ \Pi _k^h \,  \mathcal{F}^{h,\Delta}_k \,u$ of being an harmonic function. Since $\mathbb{P}_m(\mathcal{T}^S)\cap H^1_0(S) \subset {\mathcal{D}}_{0}^{0}(S) $, for polynomials of arbitrary degree $m\geq 1$, the energy norm of the second term (ii) is expected  to decay exponentially, and it is eventually   dominated by the  energy norm of the trace interpolation error  represented by the first term (i).

\subsection{Examples  of SBFEM interpolation errors in a single $S$-element}\label{ExamplesInterp}
Let us consider some examples to illustrate the accuracy capabilities of the  SBFEM interpolant, both for  smooth or boundary point singularity  harmonic functions defined in a single $S$ element,  using  refined  scaled boundary elements $\Gamma^{h,S}$, where  uniform Lagrange trace interpolation is adopted.  

\subsubsection*{Example 1 - SBFEM interpolation of a smooth harmonic function  in 2D}\label{Example_Interp2D}
\begin{figure}[!h]
	\begin{center}
	{\footnotesize
	\begin{tabular}{cccc}
\multirow{3}{*}{ $u(x,y)$}	&  \multicolumn{3}{c}{ }	 \\
 &  $h=1$ & $h=\frac{1}{2}$ & $h=\frac{1}{4}$ \\
 \includegraphics[scale=0.25]{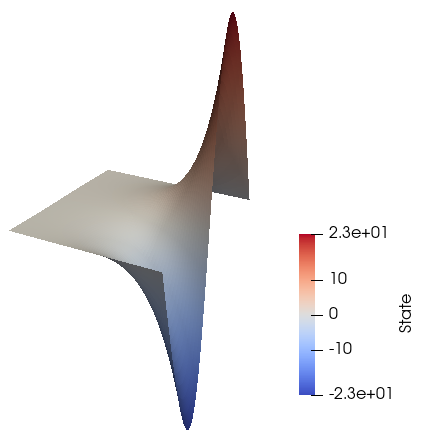}& \includegraphics[scale=0.14]{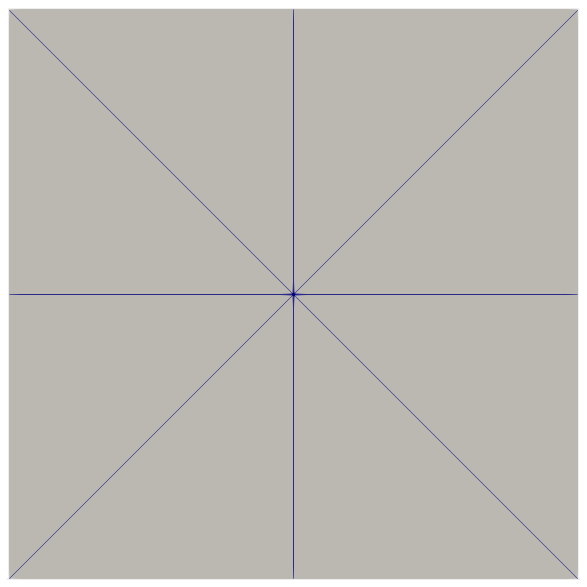} \quad & \quad  \includegraphics[scale=0.14]{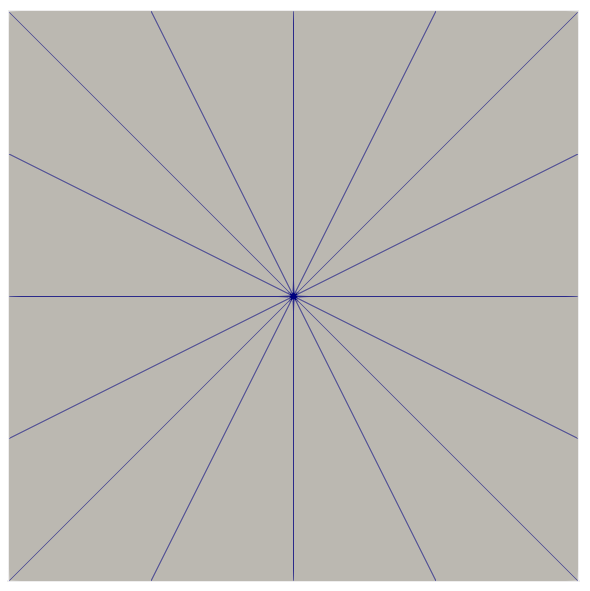}  \quad & \quad 
		\includegraphics[scale=0.14]{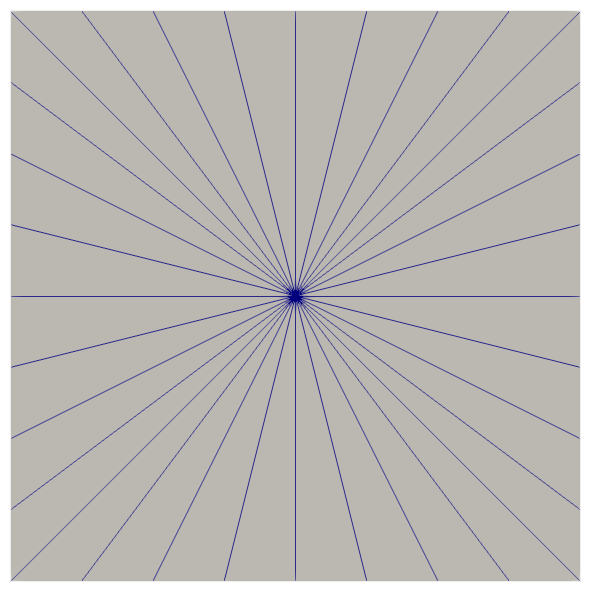} 
\end{tabular}}
	\caption{Example 1 -    Harmonic function  $u(x,y)$ and  scaled  triangular partitions $\mathcal{T}^{h,S}$ of  $S= [-1,1]\times[-1,1]$    with   scaled boundary $\Gamma^{h,S}$ formed by  $4n$ uniform facets of  width $h=\frac{2}{n}$, $n=2, 4$ and  $8$. \label{fig:example2D}}
	\end{center}
\end{figure}

In the  region  $ S=[-1,1] \times [-1, 1]$  consider the harmonic function 
$$u(x,y) = \exp{(\pi x)}\sin{(\pi y)}, $$
 and interpret $S$ as  polygonal  regions of $4n$ facets,  $n=2, 4$ and $8$,  as  illustrated in  Figure \ref{fig:example2D}.
The  scaled boundary elements $\Gamma^{h,S}$ are obtained by subdividing each   side of $\partial S$ into $n$   subintervals of width $h=\frac{2}{n}$.  In other words, $S$ is  formed by  $4n$  triangles $K^e$  sharing the scaling  center point as a vertex and   having one edge in $\Gamma^{h,S}$ as an opposite facet.     The triangles $K^e$  are mapped by Duffy's geometric transformations   described in Section \ref{sec:duffys}. 

\begin{figure}[htb]
\begin{tabular}{cc}
	\includegraphics[scale=0.3]{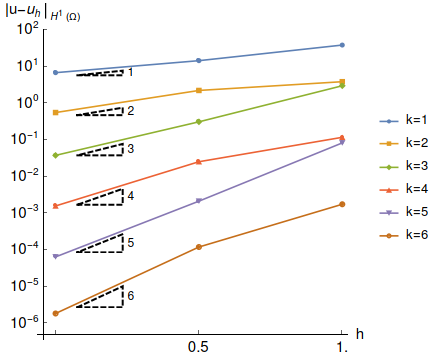} \quad & \quad \includegraphics[scale=0.3]{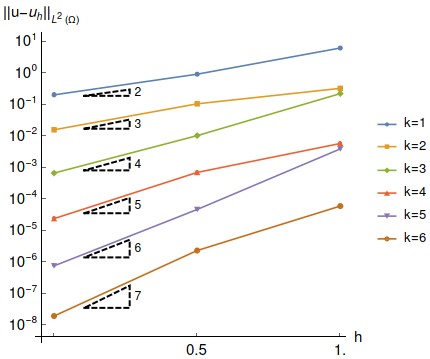}  
	\end{tabular}
	\centering
	\caption{Example 1 - Energy and  $L^2$ SBFEM interpolation errors versus $h$:   $\mathbb{S}^{h}_k(S)$ based  on the scaled  triangular partitions $\mathcal{T}^{h,S}$  of  Figure  \ref{fig:example2D}, and  trace spaces   $\Lambda^{h,S}_k$ of degree $k=1, \cdots, 6$.
	 \label{fig:convergence1element}}
\end{figure}

For these  kinds of scaled geometry, we consider the SBFEM space $\mathbb{S}^{h,S}_k$,  for $1\leq k \leq 4$,  and compute the interpolants $\Pi _k^{h,S}\, u$. The  corresponding error histories versus $h$ are plotted in Figure \ref{fig:convergence1element},  
reflecting the usual  convergence  behavior   governed by the FE trace discretizations  $\mathcal{I}_k^h\, u$ over $\partial S$, of order $k$  in the energy norm, and order $k+1$ in the $L^2$-norm.  

\subsubsection*{Example 2 - SBFEM interpolation of a smooth harmonic function  in 3D}\label{Example_Interp3D}
 The second example is for the harmonic function
\begin{equation*}
	u(x,y,z) = 4 \left(\exp{\left(\frac{\pi x}{4}\right)}\sin{\left(\frac{\pi y}{4}\right)} + \exp{\left(\frac{\pi y}{4}\right)}\sin{\left(\frac{\pi z}{4}\right)} \right) \label{eq:harmonicsolution3}
\end{equation*}
defined in   the  region  $ S=[0,1] \times [0,1] \times [0,1]$.   Let SBFEM spaces $\mathbb{S}^{h,S}_k$  obtained by considering    $S$   as polyhedral regions  with $6n^2$ facets, as  illustrated in Figure \ref{fig:mesh3D}.    The  scaled boundaries  $\Gamma^{h,S}$ are formed  by   subdividing each face in $\partial S$  into $n\times n$   quadrilaterals, and  we set the  characteristic  size   $h=\frac{1}{n}$.
Thus,  the  partitions   $\mathcal{T}^{h,S}$  are composed of $6n^2$ pyramids $K^e$  sharing the scaling  center point as a vertex, which are mapped by Duffy's geometric transformations  of the reference hexahedron, as described in Section \ref{sec:duffys}. 
\begin{figure}[!htb]
	\centering
	{\footnotesize
	\begin{tabular}{ccc}
	$h=1$ & $h=\frac{1}{2}$ & $h=\frac{1}{4}$ \\
	\includegraphics[scale=0.15]{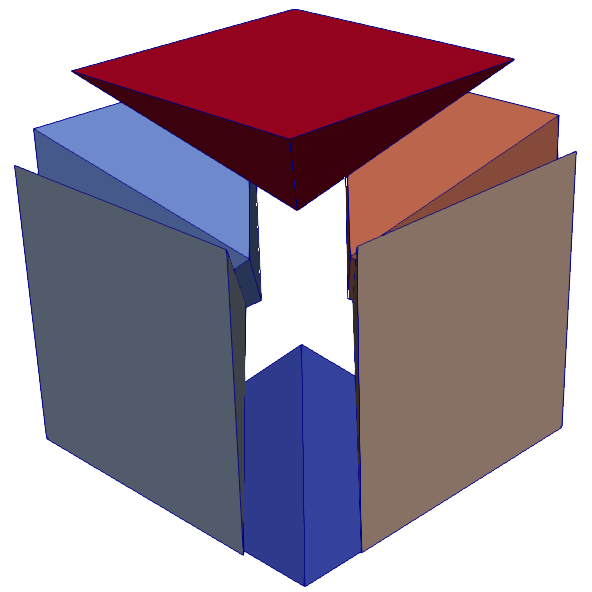} \quad & \quad \includegraphics[scale=0.14]{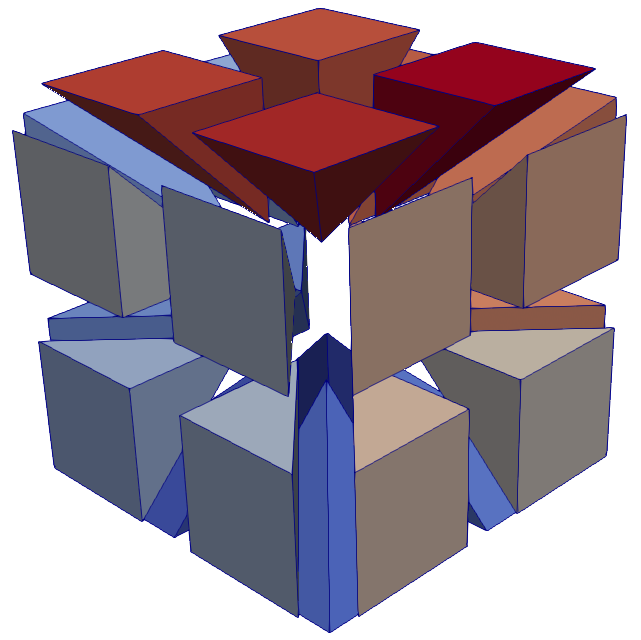} \quad & \quad \includegraphics[scale=0.14]{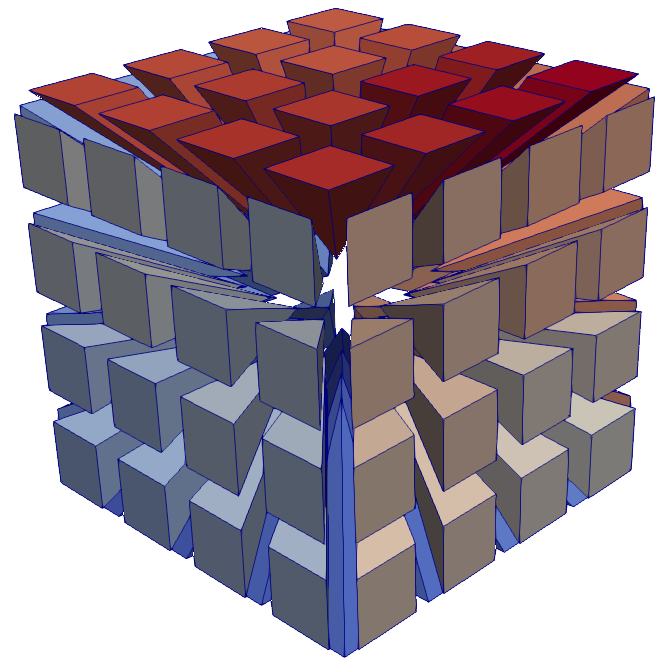} \end{tabular}}
	\caption{Example 2 -  Scaled pyramidal partitions $\mathcal{T}^{h,S}$ of   $ S= [0,1] \times [0,1] \times [0,1]$    with   scaled boundary  $\Gamma^{h,S}$ formed by $6n^2$ uniform quadrilateral  facets  of characteristic width $h=\frac{1}{n}$, $n=1, 2$ and $4$.   \label{fig:mesh3DOne}}
\end{figure}
 We approximate  $u$ by  the SBFEM interpolants   $\Pi _k^{h,S}\, u$, and the  interpolation error curves  are plotted in  Figure \ref{fig:convonecube},  
 revealing the typical optimal convergence rates of order $k$ in energy norm, and order $k+1$ in the $L^2$ norm of the trace interpolant.   
\begin{figure}[!htb]
	\centering
	\begin{tabular}{cc}
	\includegraphics[scale=0.35]{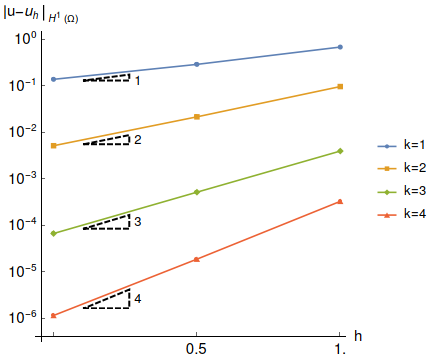} \quad & \quad 
	\includegraphics[scale=0.35]{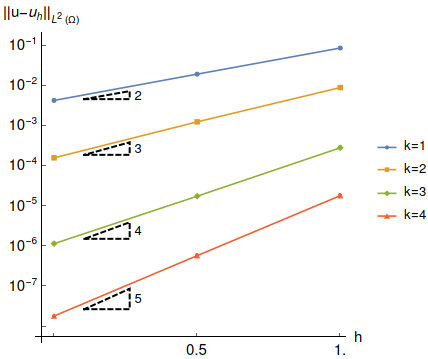}
	\end{tabular}
	\caption{Example 2  -  Energy and $L^2$ SBFEM interpolation errors   versus  $h$: $\mathbb{S}^{h}_k(S)$ based  on the scaled  pyramidal partitions $\mathcal{T}^{h,S}$ of Figure \ref{fig:mesh3DOne}, and  trace spaces   $\Lambda^{h,S}_k$ of degree $k=1,\cdots, 4$.  
	\label{fig:convonecube}}
\end{figure}

\subsubsection*{Example 3 - SBFEM interpolation of a singular harmonic function}
\begin{figure}[!h]
	\begin{center}
	{\footnotesize
	\begin{tabular}{c}
	$u(x,y)$\\
	\includegraphics[trim=2.2in 1.in 2.in 4.in, clip, 
scale=0.2]{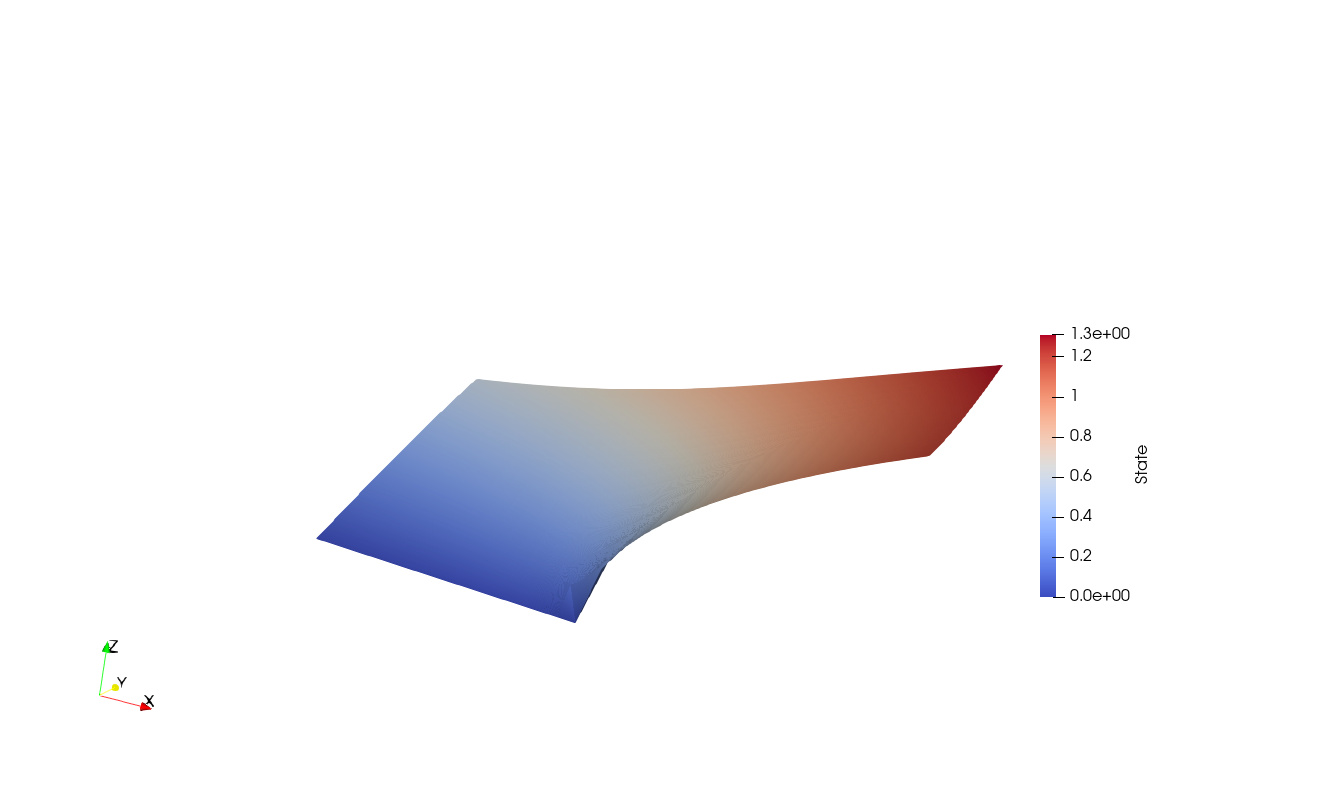}
\end{tabular}
	\begin{tabular}{cccc}
	$h=1$ & $h=\frac{1}{2}$ & $h=\frac{1}{4}$ &  $h=\frac{1}{8}$ \\
\includegraphics[scale=0.15]{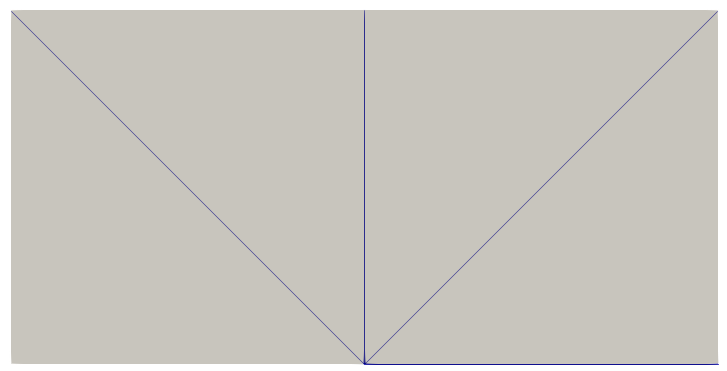} \quad & \includegraphics[scale=0.15]{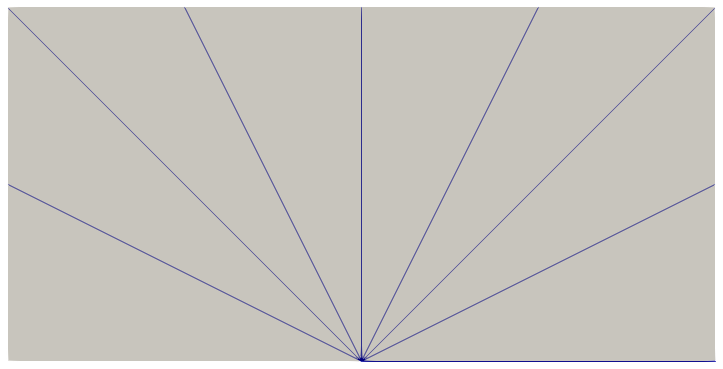} \quad &
\includegraphics[scale=0.15]{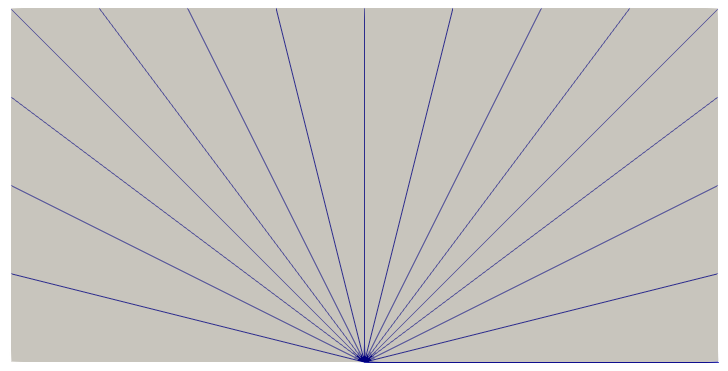} \quad &
\includegraphics[scale=0.15]{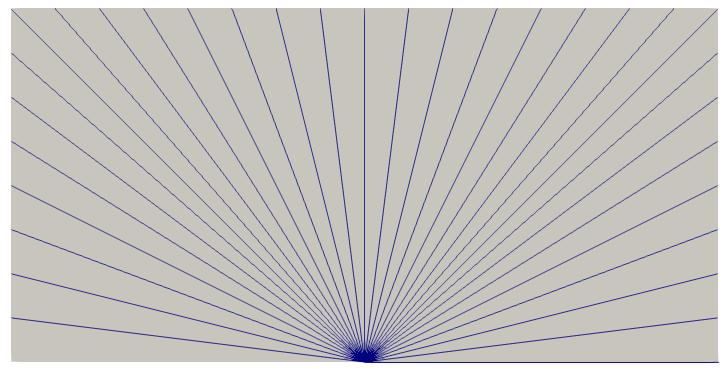}
\end{tabular}}
	\caption{Example 3 -  Singular harmonic function $u(x,y)$ and  scaled triangular partitions $\mathcal{T}^{h,S}$ of  $S=[-1,1]\times [0,1]$, with  open  scaled boundary  $\Gamma^{h,S}$, with $4n$ uniform  facets, $h=\frac{1}{n}$, $n=1, 2, 4$ and $8$. \label{fig:Singular}}
	\end{center}
\end{figure}

In the region  $S=[-1,1]\times [0,1]$ define the harmonic function
  \[u(x,y)=2^{-1/4}\sqrt{x+\sqrt{x^2+y^2}}=2^{1/4}\sqrt{r}\cos(\frac{\theta}{2}),\] 
  shown in Figure \ref{fig:Singular}, with a radial square root singularity at   the boundary point $\mathbf{O}=(0,0)$ ($r=0$), caused by boundary condition  change from  Dirichlet   $u(x,0)=0$, for $x<0$, to  Neumann   $\partial u/ \partial y (x,0) =0$, for $x>0$. This function belongs to $H^{\frac{3}{2}-\epsilon}(\Omega)$, for all $\epsilon >0$.
  
  We put the  scaling center  at  the origin  and  take an open scaled boundary  $\Gamma^{h,S}$ over  the two vertical and the top horizontal  sides of $S$,  which are uniformly     subdivided:    $n$ uniform intervals for the vertical edges, and $2n$ for the top edge,  $n=1, 2$ and  $4$.   This way,   in each refinement level, $S$ is composed of internal triangular partition $\mathcal{T}^{h,S}$  formed by  $4n$  triangles sharing the scaling center as collapsed vertex, and opposite facet width $h=\frac{1}{n}$.
Because $\Gamma^{h,S}$ is not a closed  curve,  some care must be taken in the construction  of the SBFEM  space $\mathbb{S}^{h}_k(S)$ in order to incorporate  boundary data for $u$ on the bottom boundary side of $S$. This is accomplished by enforcing in the second order SBFEM ODE system a vanishing Dirichlet boundary condition on one  side (associated with vanishing trace value at $\mathbf{x}_b=(-1,0)$), whilst a vanishing Neumann condition is assumed on the opposite side  (associated with vanishing normal trace at $\mathbf{x}_b=(1,0)$). These boundary data are  radially  extended over the  sectors $[-1,0)$ and $[0,1]$.

\begin{figure}[!htb]
	\centering
		\begin{tabular}{cc}
	\includegraphics[scale=0.35]{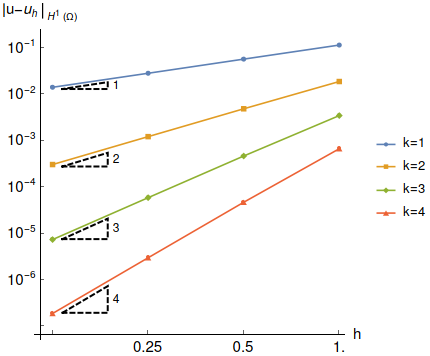} \quad & \quad 
	\includegraphics[scale=0.35]{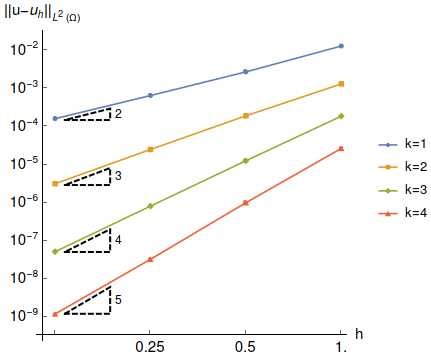}
	\end{tabular}
	\caption{Example 3  -  Energy and $L^2$ SBFEM interpolation errors    versus  $h$:  $\mathbb{S}^{h}_k(S)$ based  on the scaled triangular partitions $\mathcal{T}^{h,S}$ of Figure \ref{fig:Singular}, and  trace spaces   $\Lambda^{h,S}_k$ of degree $k=1,\cdots, 4$. 
		\label{fig:convonesingular}}
\end{figure}

SBFEM interpolation errors for this singular example are plotted in Figure \ref{fig:convonesingular}, 
revealing usual  optimal convergence rates known for trace interpolations by piecewise polynomials.   These results    reflect the role of the two terms in the decomposition  \eqref{dec_virtual},  where the  dominant contribution  is expected to come from   the virtual interpolant error, determined exclusively by the trace interpolant, which is not affected by  eventual function  singularity  not interacting with the scaled boundary $\Gamma^{h,S}$.  

\section{Galerkin SBFEM approximations}\label{sec:SBFEMGalerkin}
This section is dedicated to the Galerkin SBFEM   for   the Laplace's  model problem 
\begin{align}
    \Delta u &= 0,\ \mbox{in }\Omega, \label{eq:poisson}\\
   \gamma_0(u) &= u_D,  \mbox{on }\, \partial \Omega,   \nonumber
\end{align}
where   $u_D\in H^{1/2}(\partial \Omega)$,  and $\gamma: H^1(\Omega) \rightarrow H^{1/2}(\partial \Omega)$  is the usual trace operator.  We assume  that $u_D$ is sufficiently smooth  for the 
definition of the trace interpolant. 

Let  $\mathbb{S}^{h}_k$ be the trial   SBFEM approximation spaces based on  geometric partitions   $\mathcal{T}^h= \{S\}$ of $\Omega$ by $S$-elements,     $\Pi_k^h:  H^s(\Omega)  \rightarrow  \mathbb{S}_k^h$  being the corresponding interpolant operators, as defined in the previous section. The Galerkin SBFEM    for problem \eqref {eq:poisson} searches  approximate solutions ${u}^h \in \mathbb{S}_k^h $    satisfying: 
\begin{align}
 a( {u}_{h},v ) &= 0 \quad \forall v  \in \mathbb{S}_{k,0}^{h},
    \label{eq:discreteproblem}\\
{u}_{h}|_{\partial \Omega} &= \mathcal{I}_k^h\, {u}_{D}|_{\partial \Omega},  \label{eq:BC}
\end{align}
where  $  a( w,v):=  \int_{\Omega} \nabla u  \cdot \nabla v\, \mbox{d}\Omega$ is the usual bounded symmetric bilinear form for $u,w \in H^1(\Omega)$.  
The bilinear form $a$ is  well known to be coercive, meaning there exist $\nu >0$ such that
$a(v,v)  \geq \nu \|v\|_{H^1}^2$,   $\forall v \in H^1_0(\Omega)$.
Thus, problem \eqref{eq:discreteproblem}-\eqref{eq:BC} is well-posed  (see \cite[Proposition 3.26]{Ern2013}).

\subsection{Error analysis for the SBFEM}\label{apriori_errors}
For the error  analysis of the  Galerkin SBFEM discretization \eqref{eq:discreteproblem}-\eqref{eq:BC}, the purpose is to explore the properties \eqref{orto2D} and \eqref{eq:orto3} to estimate  energy errors  $|u - {u}^h |_{H^1}$ in approximating the harmonic exact solution $u$ from the projection errors  $|u- \mathcal{F}^{h,FE}_k \,u|_{H^1}$  on the FE spaces  $\mathcal{V}^{h,FE}_k$, or   $|u- \mathcal{F}^{h,\nabla}_k \,u|_{H^1}$  on the virtual harmonic spaces  $\mathcal{V}^{h,\Delta }_k$.  Recall that the FE interpolant errors  are available in  \cite{Babuska1987,Munoz-Sola1997,Demkowicz2007,Melenk2020} for general  functions in  Sobolev spaces,  whilst interpolant errors $|u- \mathcal{F}^{h,\Delta}_k \,u|_{H^1}$  are accessed in  \cite{Chernov2019} for harmonic functions.

  \begin{theorem}\label{apriori}  Let   $\mathcal{T}^h=\{S\}$ be a family of  polygonal partitions of $\Omega$,   $\mathbb{S}_k^h$  be the SBFEM space based 
  on $\mathcal{T}^h$, and $\mathcal{V}^{h,FE}_k$ the FE spaces based on the conglomerate meshes $\mathcal{P}^h$. Suppose the same trace interpolant is used in the definitions of  $ \Pi _k^h$  and $\mathcal{F}^{h,\Delta}_k$, and the exact solution $u \in H^1$  of  the model problem \eqref{eq:poisson} is sufficiently  regular  for them    to make sense.  If $u^h{\rm or}
  \Pi^{h,S}_k\,u \in \mathbb{S}_k^h$ is the associated Galerkin SBFEM approximation,   then 
\begin{align}
|u - {u}^h |_{H^1(\Omega)} & \leq   |u- \mathcal{F}^{h,FE}_k \,u|_{H^1(\Omega)}. \label{eq:ErProj}
\end{align}
\end{theorem}

 \begin{proof} Firstly, we observe two  orthogonality relations.
\begin{enumerate}
\item As for any Galerkin approximation,  the SBFEM  solution  verify the  orthogonality property
$a(u - {u}^h , v)=0 \quad \forall v  \in \mathbb{S}_{k,0}^{h}$, which is paramount for error estimates  for such  methods. 
\item   Proposition \ref{embbeding}  (i.e., $\mathcal{F}^{h,FE}_k \,u  \in   \mathcal{V}^{h,FE}_k  \subset {\mathcal{D}}^{h}_k$), combined with  properties \eqref{orto2D}  and  \eqref{proj_property}, implies that
\begin{equation}a(u^h, \Pi_k^h\,u - \mathcal{F}^{h,FE}_k \,u) = 0. \label{ort2}\end{equation}
\end{enumerate} 
These two orthogonality relations imply the Pythagorean equality
\begin{align*} 
| u- \mathcal{F}^{h,FE}_k \,u|_{H^1(\Omega)}^2 & = |u-{u}^h|_{H^1}^2 + |u^h- \mathcal{F}^{h,FE}_k \,u|_{H^1(\Omega)}^2.
\end{align*}
Consequently,  the estimate \eqref{eq:ErProj} holds.
\end{proof}

  \begin{theorem}\label{apriori2} Let   $\mathcal{T}^h=\{S\}$ be a family of  polygonal partitions of $\Omega$,   $\mathbb{S}_k^h$  and  $\mathcal{V}^{h,\Delta}_k$ be the SBFEM 
 and virtual spaces based on  $\mathcal{T}^h$. Suppose the same trace interpolant is used in the definitions of  $ \Pi _k^h$  and $\mathcal{F}^{h,\Delta}_k$, and the exact solution $u \in H^1$  of  the model problem \eqref{eq:poisson} is sufficiently  regular  for them    to make sense.  If $u^h \in \mathbb{S}_k^h$ is the associated Galerkin SBFEM approximation,   then
       \begin{align}
|u - {u}^h |_{H^1(\Omega)} & \leq   |u- \mathcal{F}^{h,\Delta}_k \,u|_{H^1(\Omega)}+ |\mathcal{F}^{h,\Delta}_k \,u -  \Pi _k^h \, \mathcal{F}^{h,\Delta}_k \,u|_{H^1(\Omega)}.
\label{eq:ErProj2}
\end{align}

\end{theorem}

 \begin{proof} The result is a consequence of Galerkin orthogonality property 
 \begin{align*}|u - {u}^h |_{H^1(\Omega)} = \inf_{v\in   \mathbb{S}_k^h}  |u -  v|_{H^1(\Omega)} \leq  |u -  \Pi _k^h  \,u|_{H^1(\Omega)}, \end{align*} 
 the error decomposition  \eqref{dec_virtual},  and the property 
$ \Pi _k^h  \,u =   \Pi _k^h \, \mathcal{F}^{h,\Delta}_k \,u$ remarked in \eqref{proj_propertyV}.   \end{proof}

\section{Numerical experiments}\label{sec:numericaltests}
 In this section, we  present SBFEM simulation results   for selected test problems.  First, we consider problems with smooth solutions for the verification of the predicted theoretical convergence results   of Section \ref{sec:SBFEMGalerkin}.  For a two-dimensional problem, we  explore discretizations based on quadrilateral  or polygonal $S$-elements subdivided into collapsed scaled triangles.  Then, a three-dimensional test problem is explored using SBFEM approximations based on uniform hexahedral  and  polyhedral $S$-elements subdivided into collapsed scaled pyramids,  and also on a more general geometry context of polyhedral $S$-elements subdivided by   scaled collapsed  tetrahedra. For comparison, we present results obtained by  $H^1$-conforming  FE methods based on the meshes $\mathcal{P}^h$ formed by the agglomeration of the corresponding triangles, pyramids, and  tetrahedra  partitions of the subdomains. Finally, we also evaluate the numerical performance of a coupled  FEM+SBFEM  formulation for a point singular problem, in which a traditional finite element formulation is modified by  a scaled boundary element in the vicinity of  the singularity.
 
 For the current simulations, we implemented the method in the computational framework NeoPZ\footnote{NeoPZ open-source platform: http://github.com/labmec/neopz},  which is an open source finite element library whose objective is to facilitate the development of innovative technology in finite elements \cite{Devloo1997}. Such a framework allows using a varied class of element geometries,  applying mesh refinement,  varying the approximation order, and to
  approximate partial differential equations using different approximation spaces - $H^1$, $H(\mbox{div})$, $H(\mbox{curl})$ and discontinuous -, as well as mixed and hybrid finite elements and multiscale simulations. As the NeoPZ was conceived using object-oriented concepts, with abstract classes, templates, and small blocks, it offered the required functionalities for a general coding of SBFEM  simulations. Moreover, the concept of element neighbours associated with geometric entities in NeoPZ was useful for the construction of the collapsed geometric elements and definition of the scaled boundary partitions. Two and three dimensional SBFEM approximations applied to either Laplace's equation or elasticity are implemented in a single class structure. 
 
\newcounter{example}
\refstepcounter{example}
\subsection*{Example \thesection.\theexample - smooth solution in 2D}
 The Laplace equation is approximated on the domain $\Omega = [-1;1] \times [-1;1]$, where the harmonic problem \eqref{eq:poisson} is considered with exact solution  $u(x,y) = \exp{(\pi x)}\sin{(\pi y)}$. This is the same problem of the interpolation example in Section \ref{Example_Interp2D}, illustrated in Figure \ref{fig:example2D}. 
 
\begin{figure}[h]
	\centering
	\begin{tabular}{ccc}
	 {\footnotesize Quadrilateral  $S$-elements}	\quad  \quad & \quad 	{\footnotesize  Polygonal $S$-elements - case 1}  \quad & \quad  {\footnotesize Polygonal $S$-elements - case 2}\\
	\begin{minipage}{0.2\textwidth}
		\centering
		\includegraphics[scale=0.15]{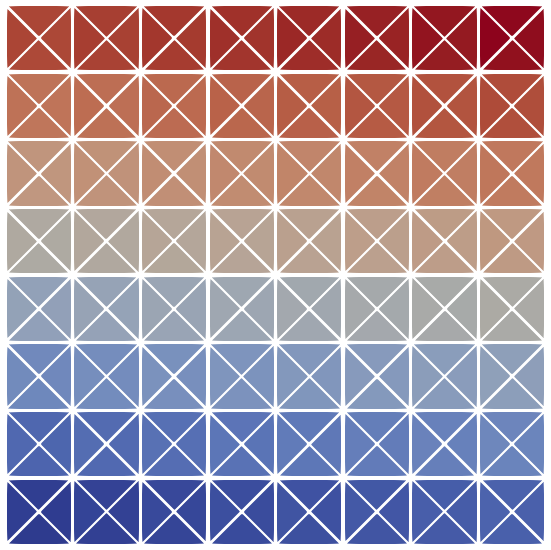}
	\end{minipage} \quad 
	&
	\quad \begin{minipage}{0.2\textwidth}
	
		\centering
		\includegraphics[scale=0.15]{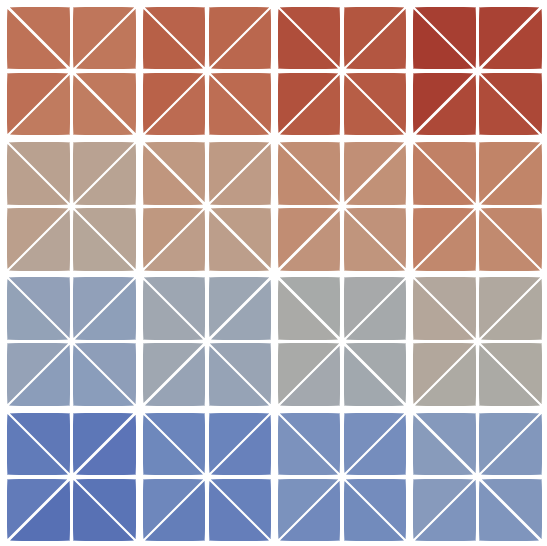}
	\end{minipage}
	&
	\begin{minipage}{0.2\textwidth}

		\centering
		\includegraphics[scale=0.23]{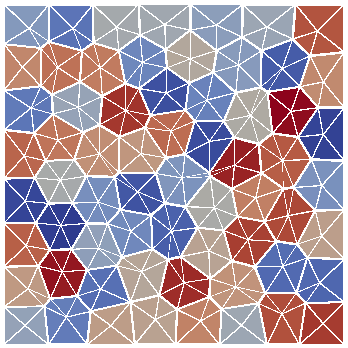}
	\end{minipage}
	\end{tabular}
	\centering
	
	\caption{Example \thesection.\theexample -  Uniform quadrilateral and  polygonal partitions $\mathcal{T}^h$, $h=\frac{1}{4}$: $S$-element  distinguished by different  colors and  subdivided into scaled  triangles.  \label{fig:geometry2D}}
\end{figure}

 We approximate  the problem     by  the Galerkin SBFEM using   sequences of partitions  $\mathcal{T}^h$ for three kinds of  $S$-elements, with  refinement levels $h=2^{-\ell}$, $\ell=1,\cdots 4$: (i)   uniform $ n\times n$,  quadrilateral     $S$-elements,   $n=2^{\ell +1}$, each one having $\Gamma^S$ formed by its  $4$ edges, (ii) polygons  with 8 edges obtained from uniform  quadrilaterals whose sides are   subdivided  once, and  (iii)  unstructured polygonal $S$-elements  constructed using the mesh generator software PolyMesher \cite{Talischi2012}, by giving as input the number of elements in $x$ and $y$ axes.  For this sequence of four irregular polygonal partitions the  scaled boundaries have average characteristic width close to  the adopt in the uniform contexts. Thus the same index parameter $h$ is adopted for them.   Recall that each $S$-element is subdivided into  triangles sharing the scaling  center point as a vertex and   having one facet in $\Gamma^{h,S}$ as  opposite edge.  Figure \ref{fig:geometry2D} illustrates the particular partitions for $h=\frac{1}{4}$.
  
 \begin{table}[!htb]
\begin{center}
	\caption{Example \thesection.\theexample -  Galerkin  SBFEM errors $E^h_{L^2}=\|u-u^h\|_{L^2(\Omega)}$  and  $E^h_{H^1}=|u-u^h|_{H^1(\Omega)}$   for  uniform   partitions  $\mathcal{T}^h$,  $h=2^{-\ell}$,  of  quadrilateral  and   polygonal  (case 1)  $S$-elements. 	} 
	\label{tab:convharmonic}
	
\vspace{.3cm}	

	\setlength{\tabcolsep}{2pt} 
	\scriptsize
	\begin{tabular}{|c|c|c|c||r|c|c||c|c|c|}
				\multicolumn{10}{c}{\footnotesize Uniform quadrilateral $S$-elements} \\
				\hline
				
{\multirow{2}{*}{ $\;\ell\;$} }		&		\multicolumn{3}{c||}{ $k$=1}    & \multicolumn{3}{c||}{ $k$=2}                                                                   & \multicolumn{3}{c|}{ $k=3$}     \\ 
		\cline{2-10}
 &		DOF  & $E^h_{L^2}$ & $E^h_{H^1}$ &  DOF & $E^h_{L^2}$ & $E^h_{H^1}$ & DOF & $E^h_{L^2}$&  $E^h_{H^1}$ \\ \hline
									
	{\color[HTML]{000000} 1}                                             & {\color[HTML]{000000} 25}   & {\color[HTML]{000000} 1.80E0}           & {\color[HTML]{000000} 1.99E1}           & {\color[HTML]{000000} 65}   & {\color[HTML]{000000} 1.31E-1}           & {\color[HTML]{000000} 2.56E0}           & {\color[HTML]{000000} 105}   & {\color[HTML]{000000} 7.78E-3}           & {\color[HTML]{000000} 2.28E-1}           \\ \hline
				
	{\color[HTML]{000000} 2}                                             & {\color[HTML]{000000} 81}   & {\color[HTML]{000000} 4.50E-1}           & {\color[HTML]{000000} 9.50E0}           & {\color[HTML]{000000} 225}   & {\color[HTML]{000000} 1.68E-2}           & {\color[HTML]{000000} 5.92E-1}           & {\color[HTML]{000000} 369}   & {\color[HTML]{000000} 4.68E-4}           & {\color[HTML]{000000} 2.62E-2}           \\ \hline
				
	{\color[HTML]{000000} 3}                                            & {\color[HTML]{000000} 289}   & {\color[HTML]{000000} 1.13E-1}           & {\color[HTML]{000000} 4.68E0}           & {\color[HTML]{000000} 833}   & {\color[HTML]{000000} 2.12E-3}           & {\color[HTML]{000000} 1.42E-1}           & {\color[HTML]{000000} 1377}  & {\color[HTML]{000000} 2.95E-5}           & {\color[HTML]{000000} 3.19E-3}           \\ \hline
	
	{\color[HTML]{000000} 4}                                            & {\color[HTML]{000000} 1089}   & {\color[HTML]{000000} 2.82E-2}           & {\color[HTML]{000000} 2.33E0}           & {\color[HTML]{000000} 3201}   & {\color[HTML]{000000} 2.65E-4}           & {\color[HTML]{000000} 3.50E-2}           & {\color[HTML]{000000} 5313}  & {\color[HTML]{000000} 1.86E-6}           & {\color[HTML]{000000} 3.96E-4}           \\ \hline
	
	& {{\color[HTML]{000000} Rate}}                                                 & {\color[HTML]{000000} 2.00}            & {\color[HTML]{000000} 1.01}            & {\color[HTML]{000000} Rate}  & {\color[HTML]{000000} 3.00}            & {\color[HTML]{000000} 2.02}            & {\color[HTML]{000000} Rate} & {\color[HTML]{000000} 3.99}            & {\color[HTML]{000000} 3.01}            \\ \hline
	\hline	
	\multicolumn{1}{|c|}{{\color[HTML]{000000} }}                    & \multicolumn{3}{c||}{{\color[HTML]{000000} $k=4$}}                        & \multicolumn{3}{c||}{{\color[HTML]{000000} $k=5$}}                                                                   & \multicolumn{3}{c|}{{\color[HTML]{000000} $k=6$}}                                                                   \\ \cline{2-10} 
	&		DOF  & $E^h_{L^2}$ & $E^h_{H^1}$ &  DOF & $E^h_{L^2}$ & $E^h_{H^1}$ & DOF & $E^h_{L^2}$&  $E^h_{H^1}$ \\ \hline
								
	{\color[HTML]{000000} 1}                                         & {\color[HTML]{000000} 145}   & {\color[HTML]{000000} 5.87E-4}           & {\color[HTML]{000000} 2.09E-2}           & {\color[HTML]{000000} 185}   & {\color[HTML]{000000} 3.89E-5}           & {\color[HTML]{000000} 1.77E-3}           & {\color[HTML]{000000} 225}   & {\color[HTML]{000000} 1.96E-6}           & {\color[HTML]{000000} 1.03E-4}           \\ \hline
				
	{\color[HTML]{000000} 2}                                         & {\color[HTML]{000000} 513}   & {\color[HTML]{000000} 1.99E-5}           & {\color[HTML]{000000} 1.29E-3}           & {\color[HTML]{000000} 657}   & {\color[HTML]{000000} 6.04E-7}           & {\color[HTML]{000000} 5.22E-5}           & {\color[HTML]{000000} 801}   & {\color[HTML]{000000} 1.53E-8}          & {\color[HTML]{000000} 1.47E-6}           \\ \hline
				
	{\color[HTML]{000000} 3}                                         & {\color[HTML]{000000} 1921}   & {\color[HTML]{000000} 6.42E-7}           & {\color[HTML]{000000} 8.14E-5}           & {\color[HTML]{000000} 2465}   & {\color[HTML]{000000} 9.47E-9}           & {\color[HTML]{000000} 1.57E-6}           & {\color[HTML]{000000} 3009}   & {\color[HTML]{000000} 1.23E-10}          & {\color[HTML]{000000} 2.21E-8}           \\ \hline
	
	{\color[HTML]{000000} 4}                                         & {\color[HTML]{000000} 7425}   & {\color[HTML]{000000} 2.03E-8}           & {\color[HTML]{000000} 5.11E-6}           & {\color[HTML]{000000} 9537}   & {\color[HTML]{000000} 1.48E-10}           & {\color[HTML]{000000} 4.81E-8}           & {\color[HTML]{000000} 11649}   & {\color[HTML]{000000} 9.46E-13}          & {\color[HTML]{000000} 3.42E-10}           \\ \hline
				
	& {{\color[HTML]{000000} Rate}}                                              & {\color[HTML]{000000} 4.98}            & {\color[HTML]{000000} 3.99}            & {\color[HTML]{000000} Rate} & {\color[HTML]{000000} 6.00}            & {\color[HTML]{000000} 5.03}            & {\color[HTML]{000000} Rate} & {\color[HTML]{000000} 6.99}            & {\color[HTML]{000000} 6.02}            \\ \hline
\multicolumn{10}{c}{}\\
				\multicolumn{10}{c}{\footnotesize Uniform polygonal $S$-elements  - case 1} \\
				\hline
			{\multirow{2}{*}{ $\;\ell\;$} }		&		\multicolumn{3}{c||}{ $k$=1}    & \multicolumn{3}{c||}{ $k$=2}                                                                   & \multicolumn{3}{c|}{ $k=3$}     \\ 
			\cline{2-10}
			&		DOF  & $E^h_{L^2}$ & $E^h_{H^1}$ &  DOF & $E^h_{L^2}$ & $E^h_{H^1}$ & DOF & $E^h_{L^2}$&  $E^h_{H^1}$ \\ \hline
			
			{\color[HTML]{000000} 1}                                             & {\color[HTML]{000000} 21}   & {\color[HTML]{000000} 8.06E-1}           & {\color[HTML]{000000} 1.23E1}           & {\color[HTML]{000000} 45}   & {\color[HTML]{000000} 8.77E-2}           & {\color[HTML]{000000} 1.95E0}           & {\color[HTML]{000000} 69}   & {\color[HTML]{000000} 8.09E-3}           & {\color[HTML]{000000} 2.53E-1}           \\ \hline
			
			{\color[HTML]{000000} 2}                                             & {\color[HTML]{000000} 65}   & {\color[HTML]{000000} 2.86E-1}           & {\color[HTML]{000000} 6.66E0}           & {\color[HTML]{000000} 145}   & {\color[HTML]{000000} 1.55E-2}           & {\color[HTML]{000000} 5.30E-1}           & {\color[HTML]{000000} 225}   & {\color[HTML]{000000} 4.98E-4}           & {\color[HTML]{000000} 2.79E-2}           \\ \hline
			
			{\color[HTML]{000000} 3}                                             & {\color[HTML]{000000} 225}   & {\color[HTML]{000000} 1.54E-2}           & {\color[HTML]{000000} 3.08E0}           & {\color[HTML]{000000} 513}   & {\color[HTML]{000000} 1.87E-3}           & {\color[HTML]{000000} 1.22E-1}           & {\color[HTML]{000000} 801}   & {\color[HTML]{000000} 3.05E-5}           & {\color[HTML]{000000} 3.24E-3}           \\ \hline
			
			{\color[HTML]{000000} 4}                                            & {\color[HTML]{000000} 833}   & {\color[HTML]{000000} 3.74E-3}           & {\color[HTML]{000000} 1.50E0}           & {\color[HTML]{000000} 1921}   & {\color[HTML]{000000} 2.30E-4}           & {\color[HTML]{000000} 2.97E-2}           & {\color[HTML]{000000} 3009}  & {\color[HTML]{000000} 1.92E-6}           & {\color[HTML]{000000} 3.94E-4}           \\ \hline
			
			& {{\color[HTML]{000000} Rate}}                                                 & {\color[HTML]{000000} 2.08}            & {\color[HTML]{000000} 1.04}            & {\color[HTML]{000000} Rate}  & {\color[HTML]{000000} 3.02}            & {\color[HTML]{000000} 2.04}            & {\color[HTML]{000000} Rate} & {\color[HTML]{000000} 3.99}            & {\color[HTML]{000000} 3.02}            \\ \hline
			\hline
			
			\multicolumn{1}{|c|}{{\color[HTML]{000000} }}                    & \multicolumn{3}{c||}{{\color[HTML]{000000} $k=4$}}                        & \multicolumn{3}{c||}{{\color[HTML]{000000} $k=5$}}                                                                   & \multicolumn{3}{c|}{{\color[HTML]{000000} $k=6$}}                                                                   \\ \cline{2-10} 
			&		DOF  & $E^h_{L^2}$ & $E^h_{H^1}$ &  DOF & $E^h_{L^2}$ & $E^h_{H^1}$ & DOF & $E^h_{L^2}$&  $E^h_{H^1}$ \\ \hline
			
			{\color[HTML]{000000} 1}                                         & {\color[HTML]{000000} 93}   & {\color[HTML]{000000} 5.49E-4}           & {\color[HTML]{000000} 2.09E-2}           & {\color[HTML]{000000} 117}   & {\color[HTML]{000000} 3.59E-5}           & {\color[HTML]{000000} 1.70E-3}           & {\color[HTML]{000000} 141}   & {\color[HTML]{000000} 1.78E-6}          & {\color[HTML]{000000} 9.68E-5}           \\ \hline
			
			{\color[HTML]{000000} 2}                                         & {\color[HTML]{000000} 305}   & {\color[HTML]{000000} 1.98E-5}           & {\color[HTML]{000000} 1.26E-3}           & {\color[HTML]{000000} 385}   & {\color[HTML]{000000} 5.76E-7}           & {\color[HTML]{000000} 4.95E-5}           & {\color[HTML]{000000} 465}   & {\color[HTML]{000000} 1.53E-8}           & {\color[HTML]{000000} 1.44E-6}           \\ \hline
			
			{\color[HTML]{000000} 3}                                         & {\color[HTML]{000000} 1089}   & {\color[HTML]{000000} 5.99E-7}           & {\color[HTML]{000000} 7.40E-5}           & {\color[HTML]{000000} 1377}   & {\color[HTML]{000000} 8.65E-9}           & {\color[HTML]{000000} 1.41E-6}           & {\color[HTML]{000000} 1665}   & {\color[HTML]{000000} 1.17E-10}          & {\color[HTML]{000000} 2.05E-8}           \\ \hline
			
			{\color[HTML]{000000} 4}                                         & {\color[HTML]{000000} 4097}   & {\color[HTML]{000000} 1.86E-8}           & {\color[HTML]{000000} 4.54E-6}           & {\color[HTML]{000000} 5185}   & {\color[HTML]{000000} 1.37E-10}           & {\color[HTML]{000000} 4.21E-8}           & {\color[HTML]{000000} 6273}   & {\color[HTML]{000000} 9.03E-13}          & {\color[HTML]{000000} 3.10E-10}           \\ \hline
						
			& {{\color[HTML]{000000} Rate}}                                              & {\color[HTML]{000000} 5.01}            & {\color[HTML]{000000} 4.03}            & {\color[HTML]{000000} Rate} & {\color[HTML]{000000} 5.98}            & {\color[HTML]{000000} 5.06}            & {\color[HTML]{000000} Rate} & {\color[HTML]{000000} 7.02}            & {\color[HTML]{000000} 6.05}            \\ \hline
		\end{tabular}
	\end{center}
\end{table}

The energy and $L^{2}$  errors  summarized in Table \ref{tab:convharmonic}  
are  for   the Galerkin SBFEM solutions in $\mathbb{S}^h_k$ based on   uniform quadrilateral  $S$-elements  and  uniform polygonal $S$-elements  of case 1,  using polynomial orders $1 \leq k \leq 6$.    The numerical results are in accordance with the predicted  rates   of order $k$   for energy errors. Optimal rates of order $k+1$ are also observed for  the errors measured by the $L^2$-norm.  

\begin{figure}[!htb]
	\centering
\begin{tabular}{ccc}
	\multicolumn{3}{c}{\footnotesize{ Polygonal  $S$-elements - case 2 } }\\
	\includegraphics[scale=0.35]{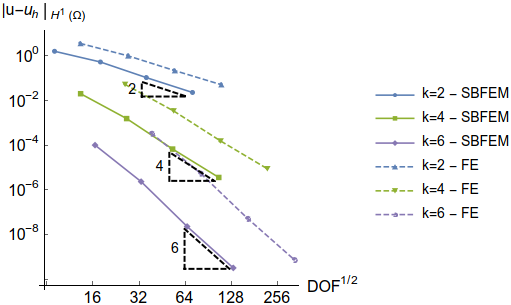} & \quad & \includegraphics[scale=0.35]{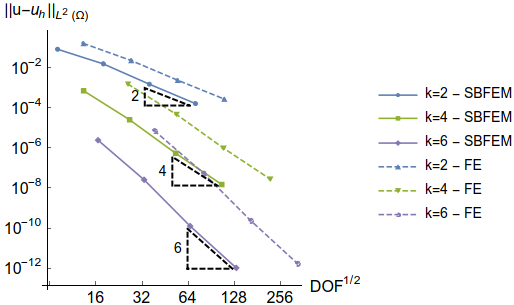} 
		\end{tabular}
	\centering
	\caption{Example \thesection.\theexample -    Energy  and $L^{2}$ errors versus  DOF for the Galerkin SBFEM  solutions in $\mathbb{S}^h_k$,   based on  the irregular  polygonal $S$-elements  of case 2  for the Galerkin FE  solutions and  for  the Galerkin FE solutions in $\mathcal{V}^{h,FE}_k$ based on the associated scaled triangular partitions $\mathcal{P}^h$, for  $k=2, 4$ and $6$.  \label{fig:convSpolygonal}}
\end{figure}

  In Figure \ref{fig:convSpolygonal}, the  energy and $L^{2}$  errors   are  plotted versus the number of DOF for Galerkin SBFEM solutions in $\mathbb{S}^h_k$ based on the polygonal meshes of case 2. For comparison, the Galerkin FE solutions in $\mathcal{V}^{h,FE}_k$ based on the associated scaled triangular partitions $\mathcal{P}^h$ are also shown, revealing comparable accuracy in both methods, but with less DOF in SBFEM simulations.  Recall that  SBFEM  shape functions are determined  by  the traces over scaled boundary elements,  whilst FE spaces are also populated with  shape functions connected with  triangular DOF other than  the edge ones opposed to the scaling center. One also observe that their error curves approach the  possible optimal  slopes   $-k$ and $-(k+1)$ when  measured by energy or
 $L^2$ norms.
  This experiment illustrates the SBFEM flexibility with respect to mesh generation for numerical simulations  without  convergence deterioration.

\begin{figure}[!htb]
 	\centering
	{\footnotesize
	\begin{tabular}{ccc}
	Quadrilateral $S$-elements \quad & \quad & \quad Quadrilateral  versus polygonal $S$-elements\\
 		\includegraphics[scale=.35]{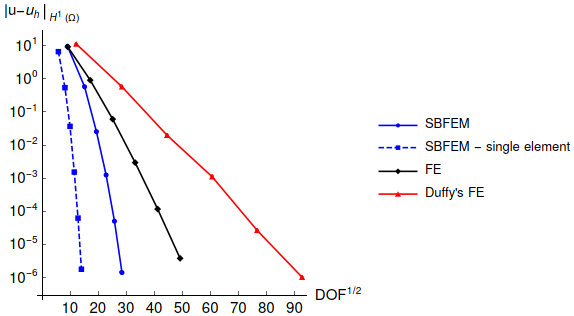} \quad & \quad  & \quad \includegraphics[scale=0.35]{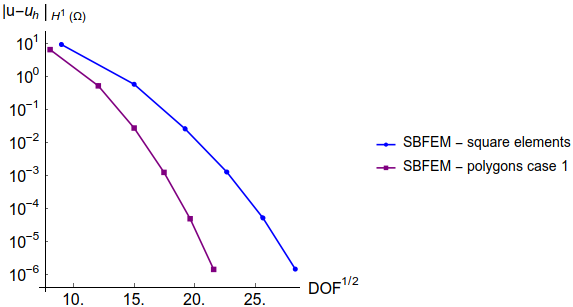} 
 \end{tabular} }
	 	\centering
	\caption{Example \thesection.\theexample -  $k$-convergence histories versus the number of DOF,  for $k=1, \cdots, 6$.  Left: SBFEM interpolation $\Pi^{h,S}_k\, u$  based on the scaled partition of  Figure \ref{fig:example2D},  Galerkin  SBFEM for $\mathbb{S}^{h}_k$,  Duffy's  FE  for $\mathcal{D}^{h}_{k,k}$, both   based on the uniform quadrilateral  partition $\mathcal{T}^h$ of Figure \ref{fig:geometry2D},   and FE method  for  $\mathcal{V}^{h,FE}_k$ based on the scaled triangular partition $\mathcal{P}^h$.  Right: Galerkin SBFEM  solutions in $\mathbb{S}^{h}_k$  using the uniform quadrilaterals and  polygonal meshes  of  case 1 shown in  Figure \ref{fig:geometry2D}.    All cases are for $h=\frac{1}{4}$. \label{fig:p conv harmonic}}
\end{figure}

Plots illustrating  SBFEM $k$-convergence  histories in the energy norm versus DOF are shown in Figure \ref{fig:p conv harmonic}, with  $k=1, \cdots,  6$, and  for $S$-elements with  fixed boundary mesh size $h=\frac{1}{4}$.
The plots on the left  are   for the SBFEM interpolation in the  single $S$-element $\Omega$   (see Figure \ref{fig:example2D}) and for the Galerkin SBFEM experiment for the uniform quadrilateral partition $\mathcal{T}^{h}$ of Figure \ref{fig:geometry2D}. For both cases, the error decay as $k$ increases shows a typical  exponential convergence,  but the interpolation experiment, by just refining the boundary of a single element,   requires  less DOF for a given accuracy threshold. 
For comparison,   $k$-convergence plots  for two  $H^1$-conforming FE methods   are also included:   using $\mathbb{P}_k(K)$ polynomials in the triangles $K$ of the conglomerate partitions $\mathcal{P}^h$ (FE),  and  for  Duffy's spaces $\mathcal{D}_{k,k}^h(S)$  (Duffy's FE).     Errors using usual FE and collapsed FE are comparable, but the  latter has more equations to be solved.  But what is more noticeable on these plots  is that 
SBFEM errors are not only smaller in magnitude  than  the FE errors (as predicted by Theorem \ref{apriori}, since energy FE errors are bounded by FE interpolation errors), but   SBFEM requires less DOF to reach  a given accuracy, the key property expected  to be held for operator adapted methods.

We also compare  in Figure  \ref{fig:p conv harmonic} (right) the $k$-convergence properties of the Galerkin SBFEM    for  spaces based on $\mathcal{T}^h$  of  the uniform quadrilateral and   polygonal $S$-elements (case 1) of  Figure \ref{fig:geometry2D}.    This comparison experiment shows  that  the use of the uniform polygonal mesh of case 1  requires fewer equations to be solved for a given target error.  On the other hand, a bigger eigenvalue system has to be solved for each $S$-element. This kind of polygonal  mesh can be seen as a combination of refining both  the boundary and inside the subdomains. Due to this flexibility, the SBFEM can  generate octree (3D) or quadtree (2D) meshes \cite{Chen2018,Saputra2020}, giving high accuracy, without any additional techniques.

\refstepcounter{example}
\subsection*{Section  \thesection.\theexample: smooth solution in 3D }
The second example refers to approximating Laplace's equation  on a 3D domain $\Omega = [0,1] \times [0,1] \times [0,1]$, with   exact harmonic solution 
\begin{equation*}
	u(x,y,z) = 4 \left(\exp{\left(\frac{\pi x}{4}\right)}\sin{\left(\frac{\pi y}{4}\right)} + \exp{\left(\frac{\pi y}{4}\right)}\sin{\left(\frac{\pi z}{4}\right)} \right). \label{eq:harmonicsolution3}
\end{equation*}
 This problem corresponds to the interpolation Example 2 of Section \ref{Example_Interp3D}. 
 
\begin{figure}[!htb]
 {\footnotesize
\centerline{Hexahedral $S$-elements \hspace{3.8cm} Polygonal $S$-elements - case 1}
\vspace{.3cm}
	\begin{center}	
	\begin{tabular}{cccccc}
	\includegraphics[scale=0.12]{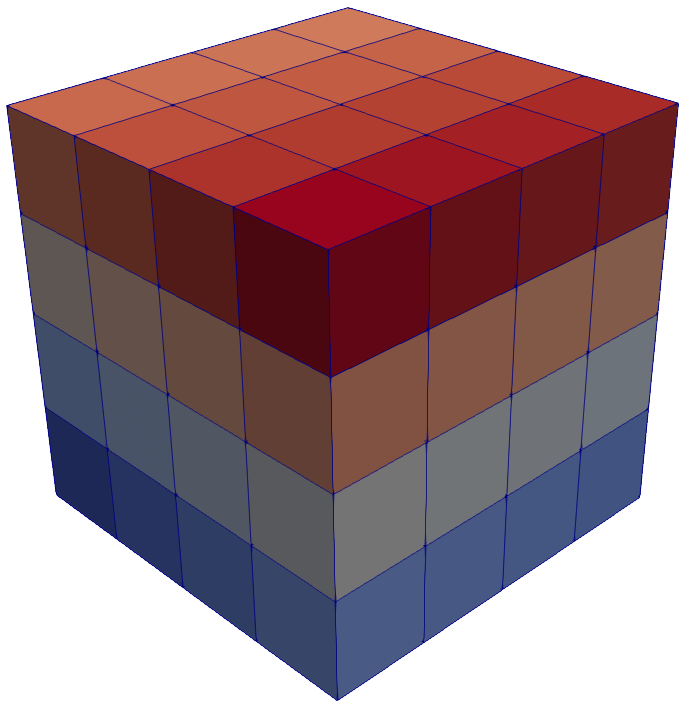} \quad & & \quad \includegraphics[scale=0.12]{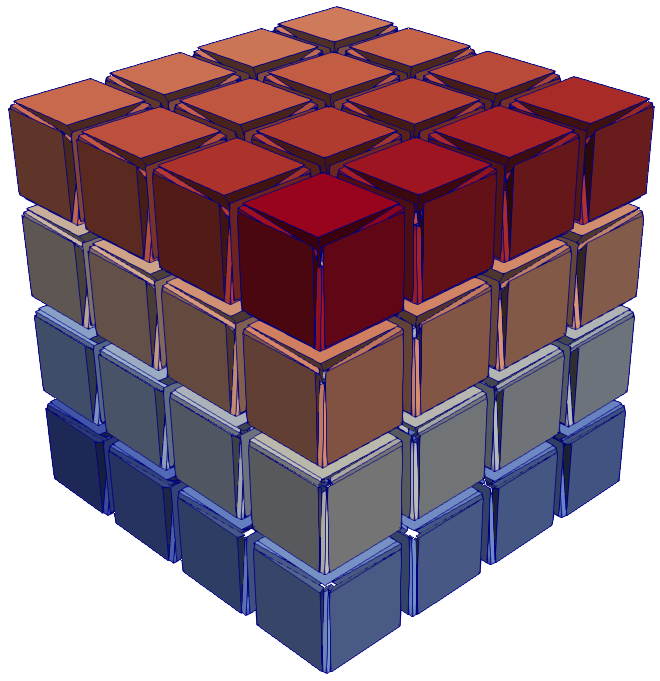} \hspace{1cm} &   \quad \includegraphics[scale=0.12]{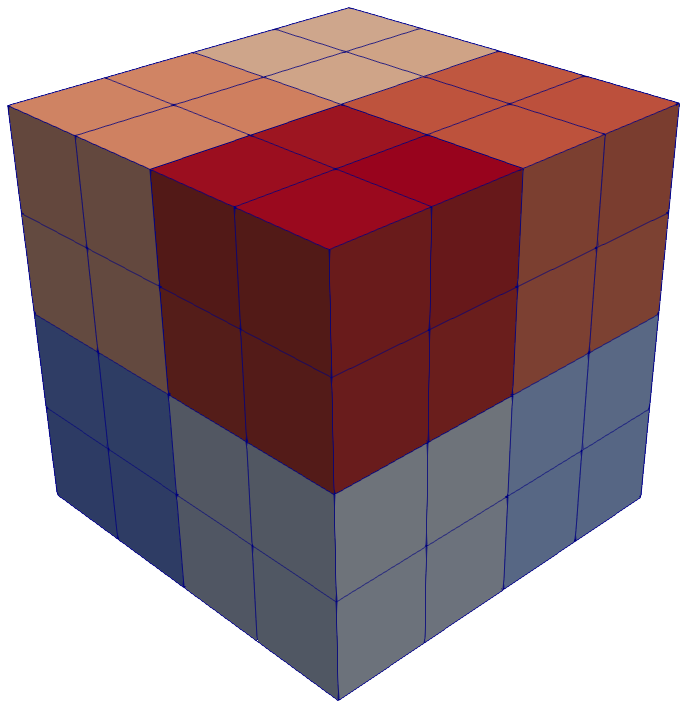} \quad  && \quad \includegraphics[scale=0.12]{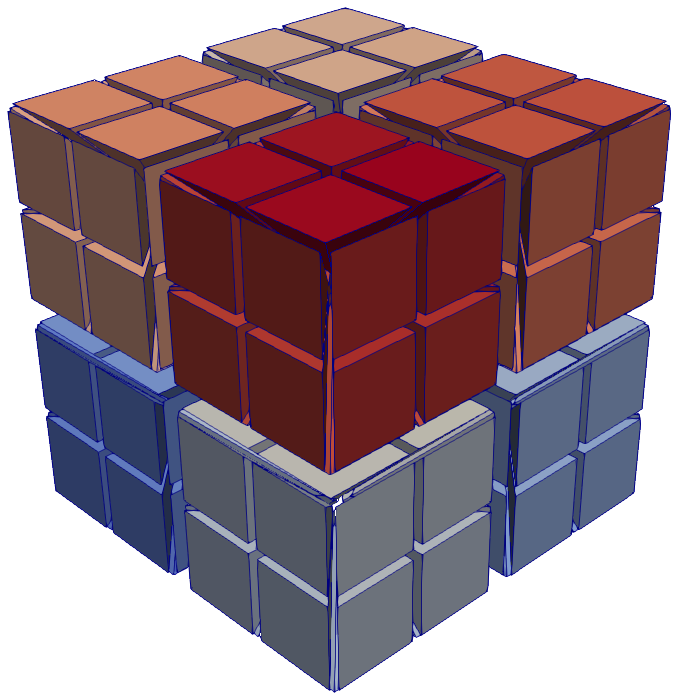}\end{tabular}\end{center}
	\vspace{-.2cm}
\centerline{Polygonal $S$-elements - case 2 }
	\vspace{-.2cm}
\begin{center}
	\begin{tabular}{cccc}
	 \includegraphics[scale=0.12]{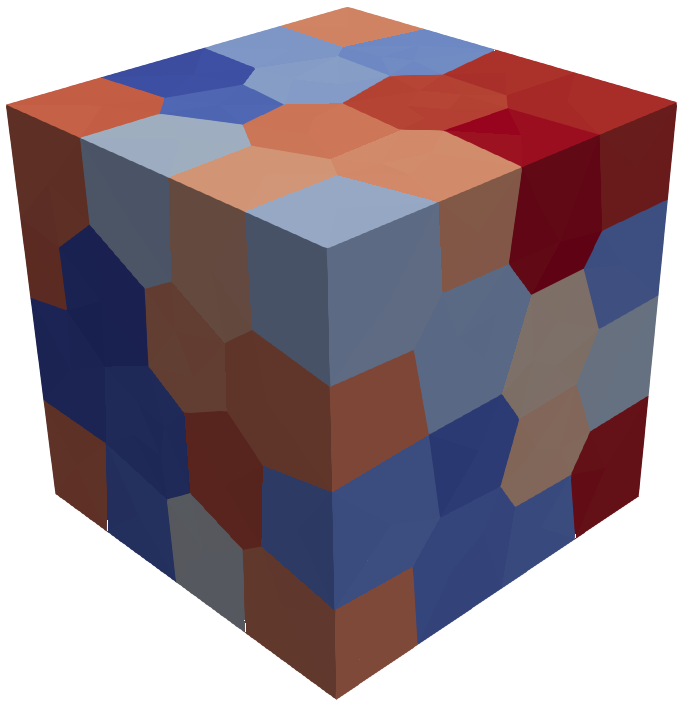}  \quad & & & \quad \includegraphics[scale=0.12]{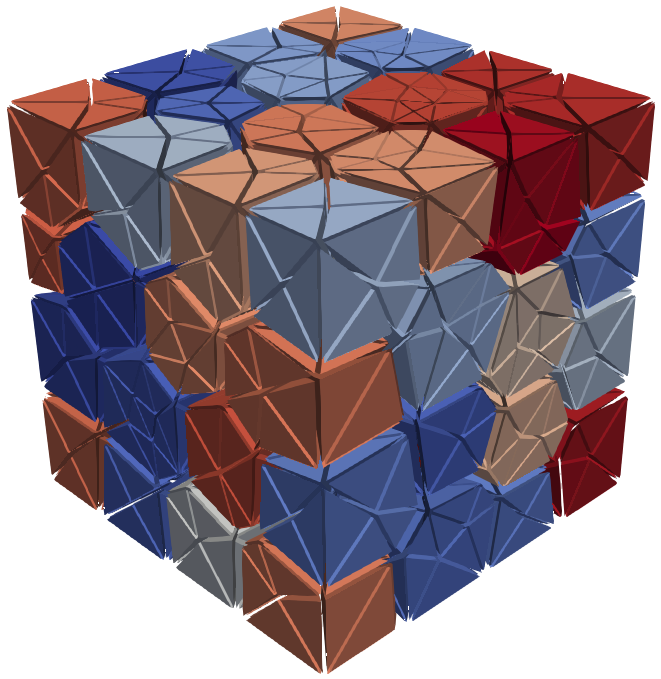}
	\end{tabular}
	\end{center}}
	\caption{Example \thesection.\theexample -  Hexahedral and polyhedral partitions $\mathcal{T}^{h}$, $h=\frac{1}{4}$:   $\mathcal{T}^{h,S}$ composed by scaled pyramids (top), and by  scaled tetrahedra (bottom). \label{fig:mesh3D}}
\end{figure}

 Three types of geometry for $\mathcal{T}^h$ are considered, each one  with   refinement levels $h=2^{-\ell}$,  $\ell=1, 2,$ and $3$.  The illustrations  in Figure \ref{fig:mesh3D} are for $h=\frac{1}{4}$.  One is for   $n\times n \times n $ uniform hexahedral partitions,  $n=2^{\ell}$,  where each $S$-element is decomposed into six pyramids. The second one is composed of polygons (case 1) constructed by subdividing once   each square face of uniform hexahedral partitions into four uniform squares (for this configuration, each $S$-element is a polyhedron with 24 quadrilateral facets, and  composed by $24$  scaled pyramids). More general polyhedral partitions (case 2) are constructed by the software package Neper \cite{Quey2011},   by giving the  number $n$ of $S$-elements in $x$, $y$, and $z$ directions.  Then, for each $S \in  \mathcal{T}^h$, we applied    gmsh \cite{Geuzaine2009}  for the construction of the internal tetrahedral   partitions  $\mathcal{T}^{h,S}$. 
  The average edge  characteristic sizes of the scaled boundary elements of these three irregular partitions resulted to be 
  comparable to the parameter $h$ of the  uniform contexts. 
  The pyramids and  tetrahedra forming $S$ are mapped by  Duffy's  transformations from the reference  hexahedron or  prism, respectively.

\begin{table}[!h] \caption{ Example \thesection.\theexample -  Galerkin  SBFEM errors $E^h_{L^2}=\|u-u^h\|_{L^2(\Omega)}$  and  $E^h_{H^1}=|u-u^h|_{H^1(\Omega)}$   for  uniform   partitions  $\mathcal{T}^h$  of  hexahedral  and   polyhedral   (case 1)  $S$-elements, with  $h=2^{-\ell}$. }
		\label{tab:convcube}
		
\setlength{\tabcolsep}{2pt} 
	\scriptsize
	\begin{center}
			\begin{tabular}{|c|c|c|c||c|c|c||c|c|c||c|c|c|}
			\multicolumn{13}{c}{\footnotesize Uniform hexahedral $S$-elements} \\
			\hline
{\multirow{2}{*}{ $\;\ell\;$} }		&		\multicolumn{3}{c||}{ $k$=1}    & \multicolumn{3}{c|}{ $k$=2}      &		\multicolumn{3}{c||}{{\color[HTML]{000000} $k=3$}}                                                                   & \multicolumn{3}{c|}{{\color[HTML]{000000} $k=4$}}             \\ 
		\cline{2-13}
 &		DOF  & $E^h_{L^2}$ & $E^h_{H^1}$ &  DOF & $E^h_{L^2}$ & $E^h_{H^1}$ &		DOF  & $E^h_{L^2}$ & $E^h_{H^1}$ &  DOF & $E^h_{L^2}$ & $E^h_{H^1}$\\ \hline				
	$1$ &		{\color[HTML]{000000} 27}    & {\color[HTML]{000000} 3.17E-2}           & {\color[HTML]{000000} 3.85E-1}           & {\color[HTML]{000000} 117}   & {\color[HTML]{000000} 1.41E-3}           & {\color[HTML]{000000} 2.40E-2}    &		{\color[HTML]{000000} 279}    & {\color[HTML]{000000} 1.92E-5}           & {\color[HTML]{000000} 5.22E-4}           & {\color[HTML]{000000} 513}   & {\color[HTML]{000000} 6.46E-7}           & {\color[HTML]{000000} 2.09E-5}         \\ \hline
			
	$2$&		{\color[HTML]{000000} 127}   & {\color[HTML]{000000}  7.85E-3}           & {\color[HTML]{000000}  1.87E-1}           & {\color[HTML]{000000} 665}   & {\color[HTML]{000000} 1.93E-4}           & {\color[HTML]{000000} 6.02E-3}     &		{\color[HTML]{000000} 1685}   & {\color[HTML]{000000} 1.25E-6}           & {\color[HTML]{000000} 6.49E-5}            & {\color[HTML]{000000} 3185}   & {\color[HTML]{000000} 2.12E-8}           & {\color[HTML]{000000} 1.31E-6}      \\ \hline
			
	$3$ &		{\color[HTML]{000000} 729}   & {\color[HTML]{000000} 1.93E-3}           & {\color[HTML]{000000} 9.20E-2}           & {\color[HTML]{000000} 4401}   & {\color[HTML]{000000} 2.48E-5}           & {\color[HTML]{000000} 1.51E-3}       &		{\color[HTML]{000000} 11529}   & {\color[HTML]{000000} 7.99E-8}           & {\color[HTML]{000000} 8.07E-6}           & {\color[HTML]{000000} 22113}   & {\color[HTML]{000000} 6.75E-10}           & {\color[HTML]{000000} 8.16E-8}     \\ \hline
			
	& {{\color[HTML]{000000} Rate}}                                                     & {\color[HTML]{000000} 2.03}            & {\color[HTML]{000000} 1.02}            & {\color[HTML]{000000} Rate} & {\color[HTML]{000000} 2.96}            & {\color[HTML]{000000} 2.00}    &     {{\color[HTML]{000000} Rate}}                                                    & {\color[HTML]{000000} 3.97}            & {\color[HTML]{000000} 3.01}            & {\color[HTML]{000000} Rate} & {\color[HTML]{000000} 4.98}            & {\color[HTML]{000000} 4.00}                \\ 
			\hline
			\multicolumn{13}{c}{}\\
\multicolumn{13}{c}{\footnotesize Uniform polyhedral $S$-elements - case 1} \\
\hline
			{\multirow{2}{*}{ $\;\ell\;$} }		&		\multicolumn{3}{c||}{ $k$=1}    & \multicolumn{3}{c|}{ $k$=2}      &		\multicolumn{3}{c||}{{\color[HTML]{000000} $k=3$}}                                                                   & \multicolumn{3}{c|}{{\color[HTML]{000000} $k=4$}}             \\ 
			\cline{2-13}
			&		DOF  & $E^h_{L^2}$ & $E^h_{H^1}$ &  DOF & $E^h_{L^2}$ & $E^h_{H^1}$ &		DOF  & $E^h_{L^2}$ & $E^h_{H^1}$ &  DOF & $E^h_{L^2}$ & $E^h_{H^1}$\\ \hline				
			$1$ &		{\color[HTML]{000000} 26}    & {\color[HTML]{000000} 1.94E-2}           & {\color[HTML]{000000} 2.96E-1}           & {\color[HTML]{000000} 98}   & {\color[HTML]{000000} 1.25E-3}           & {\color[HTML]{000000} 2.21E-2}    &		{\color[HTML]{000000} 218}    & {\color[HTML]{000000} 1.76E-5}           & {\color[HTML]{000000} 5.30E-4}           & {\color[HTML]{000000} 386}   & {\color[HTML]{000000} 5.83E-7}           & {\color[HTML]{000000} 1.92E-5}         \\ \hline
			
			$2$ &		{\color[HTML]{000000} 117}    & {\color[HTML]{000000} 5.24E-3}           & {\color[HTML]{000000} 1.42E-1}           & {\color[HTML]{000000} 513}   & {\color[HTML]{000000} 1.68E-4}           & {\color[HTML]{000000} 5.34E-3}    &		{\color[HTML]{000000} 1197}    & {\color[HTML]{000000} 1.19E-6}           & {\color[HTML]{000000} 6.45E-5}           & {\color[HTML]{000000} 2169}   & {\color[HTML]{000000} 1.91E-8}           & {\color[HTML]{000000} 1.18E-6}         \\ \hline
			
			$3$&		{\color[HTML]{000000} 665}   & {\color[HTML]{000000}  1.35E-3}           & {\color[HTML]{000000}  7.00E-2}           & {\color[HTML]{000000} 3185}   & {\color[HTML]{000000} 2.14E-5}           & {\color[HTML]{000000} 1.33E-3}     &		{\color[HTML]{000000} 7625}   & {\color[HTML]{000000} 7.79E-8}           & {\color[HTML]{000000} 7.92E-6}            & {\color[HTML]{000000} 13985}   & {\color[HTML]{000000} 6.08E-10}           & {\color[HTML]{000000} 7.32E-8}      \\ \hline
			
						& {{\color[HTML]{000000} Rate}}                                                     & {\color[HTML]{000000} 1.96}            & {\color[HTML]{000000} 1.02}            & {\color[HTML]{000000} Rate} & {\color[HTML]{000000} 2.97}            & {\color[HTML]{000000} 2.01}    &     {{\color[HTML]{000000} Rate}}                                                    & {\color[HTML]{000000} 3.93}            & {\color[HTML]{000000} 3.03}            & {\color[HTML]{000000} Rate} & {\color[HTML]{000000} 4.97}            & {\color[HTML]{000000} 4.01}                \\ 
			\hline
		\end{tabular}
	\end{center}
\end{table}

\begin{figure}[!htb]
	\centering
\begin{tabular}{ccc}
	 \multicolumn{3}{c}{\footnotesize{Polyhedral $S$-elements - case 2 }}\\
	\includegraphics[scale=0.4]{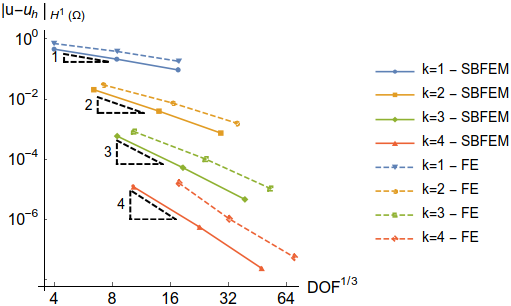} \quad &\quad   & \quad \includegraphics[scale=0.4]{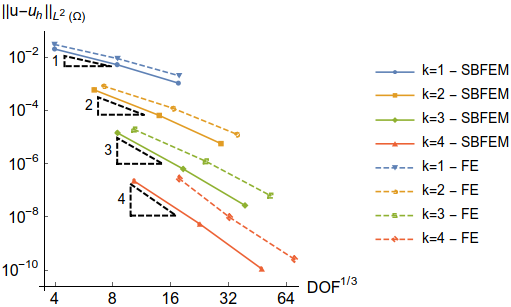}
	\end{tabular}
	\caption{ Example \thesection.\theexample -   Energy  and $L^{2}$  errors   versus  DOF for the Galerkin SBFEM  solution in $\mathbb{S}^h_k$,   for $k=1, \cdots, 4$,  based on the irregular  polyhedral $S$-elements of case 2.  \label{fig:convPolyhdral}}
\end{figure}

  The results for  the Galerkin SBFEM solutions in   $\mathbb{S}^{h}_k$,  with $k=1, \cdots, 4$,  based on the uniform hexahedral
 S-elements and on the   polyhedral $S$-elements of case 1 are   documented  in Table \ref{tab:convcube}.  Optimal  accuracy of order $k$ for energy norm, and $k+1$ for  the $L^2$-norm  occur.   Energy and
 $L^2$ errors  obtained with the polyhedral partitions of case 2 are plotted versus DOF in Figure \ref{fig:convPolyhdral}.  For
 comparison, the Galerkin FE solutions in $\mathcal{V}^{h,FE}_k$  based on the associated scaled tetrahedral 
 partitions $\mathcal{P}^h$ are also shown.   Similar conclusions hold as for the experiment shown in Figure \ref{fig:convSpolygonal}. One can also observe that both Galerkin SBFEM and FE approximation errors have similar  magnitude, but with less DOF in the SBFEM systems. Their error curves   measured with energy and
 $L^2$ norms also  approach the   possible  optimal  slopes  $-k$ and $-(k+1)$, respectively. 

\begin{figure}[!htb]
\begin{center}
{\footnotesize
	\begin{tabular}{cc}
	Hexahedral $S$-elements & {\small Hexahedral vs. Polyhedral $S$-elements} \\
	\includegraphics[scale=0.35]{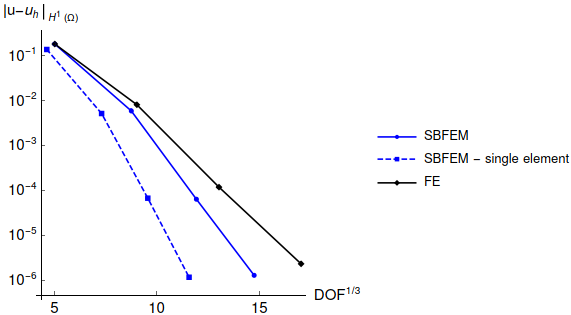} & 	\includegraphics[scale=0.35]{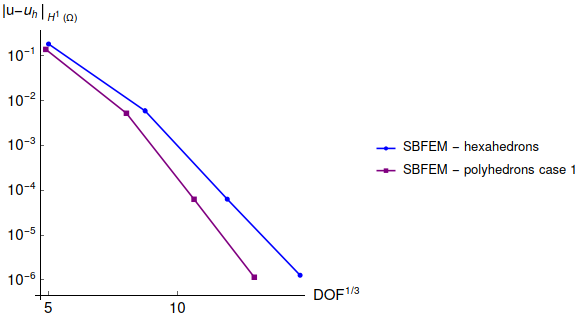}		\end{tabular}}
	\caption{Example \thesection.\theexample -  $k$-convergence histories as function of the number of DOF,  for $k=1, \cdots, 4$: Left: SBFEM interpolation $\Pi^{h,S}_k\, u$  for the scaled partition $\mathcal{T}^{h,S}$ of Figure \ref{fig:mesh3DOne},   Galerkin  SBFEM  for  $ \mathbb{S}^{h}_k$  based on uniform hexahedral  partition $\mathcal{T}^h$,  and  FE method  for  $\mathcal{V}^{h,FE}_k$ based on the  conglomerated scaled  pyramidal partition $\mathcal{P}^h$.   Right:  Galerkin SBFEM  for  $ \mathbb{S}^{h}_k$  based on hexahedral and polyhedral $S$-elements of case 1.  In all the experiments, $h=\frac{1}{4}$. \label{fig:pconvergence3D}}
	\end{center}
\end{figure}
 In the left side of Figure \ref{fig:pconvergence3D}, we compare  the SBFEM $k$-convergence using the fixed  uniform hexahedral  partition at the refinement level $h=\frac{1}{4}$, shown in Figure  \ref{fig:mesh3D},
with  equivalent results for the FE method  using the  spaces $\mathcal{V}^{h,FE}_k \subset H^1(\Omega)$ based on the associated pyramidal partition $\mathcal{P}^h$.  SBFEM approximations lead to lower error values, as predicted by Theorem \ref{apriori}, and the linear systems have a reduced number of equations. The error curve of the interpolation experiment  illustrated in Figure \ref{fig:convonecube} is also included. 

The plots on the right side compare the the $k$-convergence of the two SBFEM solutions  in $\mathbb{S}^h_k$ based on the uniform hexahedral  partition and on the polyhedral partition of case 1  illustrated in   Figure  \ref{fig:mesh3D},  both with $h=\frac{1}{4}$.   Similarly to the comparison experiment of  the previous example, shown in Figure  \ref{fig:p conv harmonic},  these convergence histories  also show  that  the use of polygonal mesh  requires fewer equations to be solved for a given target error, but reminding that  it requires  bigger eigenvalue systems  to be solved for the computatiobn of  SBFEM shape functions in  the $S$-elements.

\refstepcounter{example}
\subsection*{Section  \thesection.\theexample:   coupled FE-SBFEM formulation for a singular problem}
Taking the  singular  harmonic function  interpolated  in  Section \ref{ExamplesInterp}, namely
  \[u=2^{1/4}\sqrt{r}\cos(\frac{\theta}{2})=2^{-1/4}\sqrt{x+\sqrt{x^2+y^2}},\]   
 we enforce  Dirichlet  boundary condition  on  $(x,0)$,  $x<0$,  and  Neumann boundary condition elsewhere.  Due to the lack of regularity of  $u\in H^{\frac{3}{2}-\epsilon}(\Omega)$, the error estimates of Theorem \ref{apriori} in terms of FE interpolant error based on regular partitions are restricted in theory to order $h^{\frac{1}{2}-\epsilon}$. This problem was  considered   in  \cite{DeSiqueira2020} to evaluate   the efficiency of  the mixed FE method  when quarter-point elements are used in the vicinity of the origin  $\mathbf{O}=(0,0)$ (singular point), showing dramatic accuracy improvement. Recall that  the specific 6-noded quarter-point element is also of Duffy's type, obtained by collapsing a reference quadrilateral element on  triangles.  

	\begin{figure}[h]
	\begin{center}
	{\small
	\begin{tabular}{ccc}
	$h=\frac{1}{2}$ \quad & \quad & \quad $h=\frac{1}{16}$ \\
		\centering
		\begin{minipage}{0.32\textwidth}
			\centering
			\includegraphics[scale=0.2]{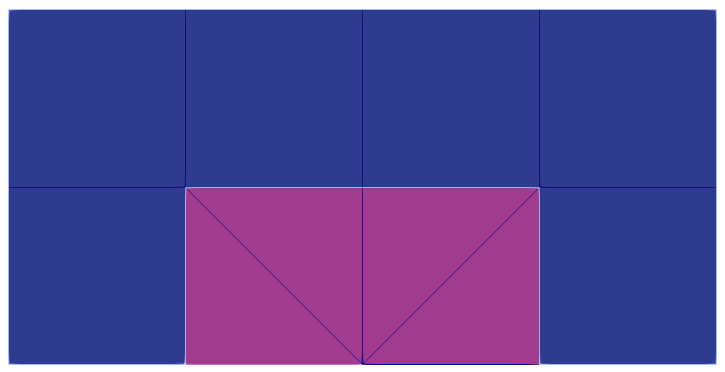}
		\end{minipage}
		 \quad & \quad & \quad 
		\begin{minipage}{0.32\textwidth}
			\centering
			\includegraphics[scale=0.2]{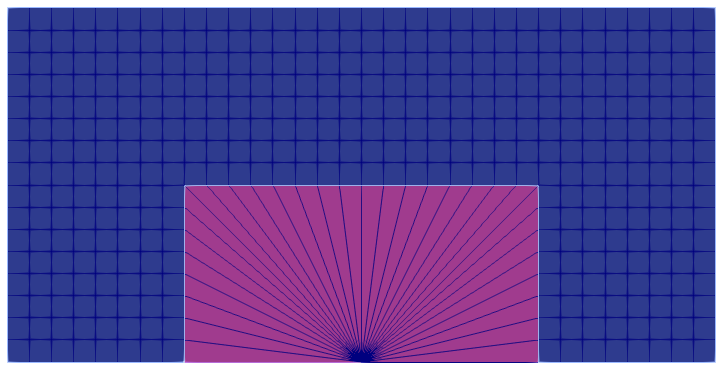}
		\end{minipage}
		\centering
		\end{tabular}}
		\end{center}
		\caption{Example \thesection.\theexample -  Meshes for the coupled FE-SBFEM formulation:  FE (blue) in the smooth region and SBFEM (magenta)  close to the singularity point. \label{fig:singularhref}}
	\end{figure}
	
With this motivation,  we propose  a formulation composing   SBFEM approximations in a  single  element  $S=[-0.5,0.5]\times [0, 0.5]$  and FE  approximations elsewhere, in the region where the solution is smooth.  Similarly to the interpolation experiments  in Section \ref{ExamplesInterp},  the space $\mathbb{S}^{h}_k(S)$ is conceived in such a way that the scaling center is located on the singularity point, which means that an open scaled boundary element is applied. The vertical and top-horizontal edges of $S$  are uniformly subdivided to form an interface partition   $\Gamma^{h,S}$. Elsewhere,  a uniform quadrilateral mesh matching   $\Gamma^{h,S}$ is adopted, as illustrated in Fig \ref{fig:singularhref} for $h = \frac{1}{2}$ and   $h = \frac{1}{16}$.  The coupling between FE and SBFEM approximations is  straightforward  since SBFEM uses compatible FE  spaces at the interface.
Four mesh  sizes  $h = 2^{-\ell}$,  $\ell = 1, \cdots, 4$, and polynomials  of degree $k=1, \cdots, 4$ are performed.  

 The corresponding results  are documented in Table \ref{tab:convsbfemfem}. As for regular problems with smooth solutions, optimal rates of convergence of order $k$ and $k+1$ for energy and $L^2$ errors  hold for this singular problem, without any  adaptivity, i.e.  uniform  degree $k$ is used over the domain and no $h$-adaptivity is applied as well.  
	
	\begin{table}[!h]
	\caption{Example \thesection.\theexample  -  Errors  $E^h_{L^2}=\|u-u^h\|_{L^2(\Omega)}$  and  $E^h_{H^1}=|u-u^h|_{H^1(\Omega)}$, $h=2^{-\ell}$,  for  the combined FE-SBFEM  method.}
	\label{tab:convsbfemfem}
	
	\setlength{\tabcolsep}{2pt} 
	\scriptsize
	\begin{center}
		\begin{tabular}{|c|c|c|c||c|c|c|c|c|c|c|c|c|}
			\hline
			{\multirow{2}{*}{ $\;\ell\;$} }		&		\multicolumn{3}{c||}{ $k$=1}    & \multicolumn{3}{c|}{ $k$=2}    &	\multicolumn{3}{c||}{{\color[HTML]{000000} $k=3$}}                                                                   & \multicolumn{3}{c|}{{\color[HTML]{000000} $k=4$}}                                                                          \\ 
			\cline{2-13}
			&		DOF  & $E^h_{L^2}$ & $E^h_{H^1}$ &  DOF & $E^h_{L^2}$ & $E^h_{H^1}$  & 	DOF  & $E^h_{L^2}$ & $E^h_{H^1}$ &  DOF & $E^h_{L^2}$ & $E^h_{H^1}$\\ \hline	
			
			$1$ &		{\color[HTML]{000000} 14}    & {\color[HTML]{000000} 8.44E-4}           & {\color[HTML]{000000} 1.11E-1}           & {\color[HTML]{000000} 26}    & {\color[HTML]{000000} 7.87E-4}           & {\color[HTML]{000000} 1.94E-2}   &   {\color[HTML]{000000} 76}    & {\color[HTML]{000000} 9.17E-5}           & {\color[HTML]{000000} 2.82E-3}           & {\color[HTML]{000000} 125}    & {\color[HTML]{000000} 1.19E-5}           & {\color[HTML]{000000} 5.25E-4}       \\ \hline
			
			$2$ &		{\color[HTML]{000000} 39}    & {\color[HTML]{000000} 2.02E-3}           & {\color[HTML]{000000} 5.54E-2}           & {\color[HTML]{000000} 117}   & {\color[HTML]{000000} 1.12E-4}           & {\color[HTML]{000000} 4.23E-3}      &		{\color[HTML]{000000} 259}    & {\color[HTML]{000000} 6.35E-6}           & {\color[HTML]{000000} 3.83E-4}           & {\color[HTML]{000000} 441}   & {\color[HTML]{000000} 4.59E-7}           & {\color[HTML]{000000} 3.72E-5}         \\ \hline
			
			$3$&		{\color[HTML]{000000} 125}   & {\color[HTML]{000000}  4.95E-4}           & {\color[HTML]{000000}  2.73E-2}           & {\color[HTML]{000000} 665}   & {\color[HTML]{000000} 1.45E-5}           & {\color[HTML]{000000} 1.06E-3}     &		{\color[HTML]{000000} 949}   & {\color[HTML]{000000} 4.16E-7}           & {\color[HTML]{000000} 4.87E-5}            & {\color[HTML]{000000} 1649}   & {\color[HTML]{000000} 1.54E-8}           & {\color[HTML]{000000} 2.39E-6}      \\ \hline
			
			$4$ &		{\color[HTML]{000000} 441}   & {\color[HTML]{000000} 1.23E-4}           & {\color[HTML]{000000} 1.36E-2}           & {\color[HTML]{000000} 4401}   & {\color[HTML]{000000} 1.84E-6}           & {\color[HTML]{000000} 2.66E-4}       &		{\color[HTML]{000000} 3625}   & {\color[HTML]{000000} 2.67E-8}           & {\color[HTML]{000000} 6.10E-6}           & {\color[HTML]{000000} 6369}   & {\color[HTML]{000000} 4.97E-10}           & {\color[HTML]{000000} 1.50E-7}      \\ \hline
			
		& {{\color[HTML]{000000} Rate}}                                                     & {\color[HTML]{000000} 2.01}            & {\color[HTML]{000000} 1.01}            & {\color[HTML]{000000} Rate} & {\color[HTML]{000000} 2.99}            & {\color[HTML]{000000} 2.00}      &       {{\color[HTML]{000000} Rate}}                                                    & {\color[HTML]{000000} 3.96}            & {\color[HTML]{000000} 3.00}            & {\color[HTML]{000000} Rate} & {\color[HTML]{000000} 4.96}            & {\color[HTML]{000000} 3.99}                 \\ 
			\hline
		\end{tabular}
	\end{center}
\end{table}

\begin{figure}[!htb]
		   \begin{center}
		   {\footnotesize
		\begin{tabular}{c}
		SBFEM interpolation\\
		\includegraphics[scale=0.12]{Figures/oneelementrefskeletonsing3.png}  \\ 
		FE-SBFEM\\
		 \includegraphics[scale=0.12]{Figures/singularmeshnref3.png}  
\end{tabular}   \hspace{1cm} 	\begin{tabular}{c} \\  \includegraphics[scale=0.35]{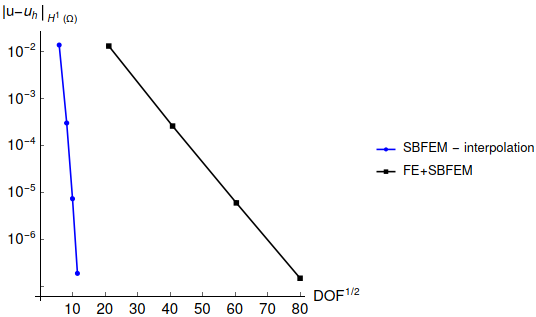}\end{tabular}}
	   \end{center}
		\caption{Example \thesection.\theexample - Partitions and  $k$-convergence histories versus the number of DOF,  with $k=1, \cdots, 4$, for  SBFEM interpolants $\Pi^{h,S}\, u$ of Example 3,  Section \ref{ExamplesInterp}, and    Galerkin FE-SBFEM solutions. 	\label{fig:pconvsbfemfem}}
	\end{figure}

For comparison,  two $k$-convergence histories as function of the number of DOF   are shown  in Figure \ref{fig:pconvsbfemfem}  for fixed partitions of the domain $\Omega$: one  for the SBFEM  interpolation errors  computed in Example 3,  Section \ref{ExamplesInterp},  and the other  for the  combined  Galerkin FE-SBFEM method.  The  partitions used in these experiments are  illustrated in Figure \ref{fig:pconvsbfemfem},  noticing that they coincide within the region $S$ around the singularity, but the FE partition in the smooth region being more refined. Whilst SBFEM interpolation in the single element $\Omega$ requires much less DOF, both experiments reach very close  error values, 
because the error in this problem is governed by the singularity, modeled using SBFEM in both experiments. 
However,  the results of Fig \ref{fig:pconvsbfemfem} could be deceiving. It should be emphasized that global SBFEM interpolation in the whole domain $\Omega$ was feasible in this particular  test problem, but this will not the case in practical  singular problems, for which coupled FE+SBFEM simulations reveal to be  a simple and efficient option. 

\section{Conclusions} 
\label{sec:conclusions}

We provide a priori error estimates in energy norm for  Galerkin SBFEM approximations of harmonic solutions by exploring  two aspects of SBFEM's methodology. 

The SBFEM approximation spaces are  based on star-shaped polytopal subregions ($S$-elements), where the functions are parametrized in  the radial  and surface directions. We show that  they can be  presented in the context of Duffy's approximations based on sub-partitions of the $S$-elements.   Piecewise polynomial discretization is adopted for surface  traces,  which are radially extended to   the interior of $S$ by solving  local harmonic problems using test functions restricted to Duffy's spaces. As a consequence,  shape functions can be derived from analytical solutions  defined by eigenvalue  problems,  whose coefficients are determined by the geometry of the $S$-elements. 

We demonstrate that there is an equivalence between the SBFEM ODE equation and an orthogonality  property of SBFEM spaces, with respect to the gradient inner product for a wide class of  Duffy's approximations vanishing on the  facets of $S$. This orthogonal   property is the key for  the derivation  of the estimation of SBFEM errors  in energy norm. The Galerkin SBFEM approximation error is necessarily smaller than the FE interpolant error  for the FE space   included in the Duffy's space sharing the same interface traces.

We show that SBFEM errors in the approximation of harmonic functions  come from two sources: there is the kind of error caused when the trace of  harmonic functions are discretized over the facets of $S$, occurring in virtual harmonic approximations, and there is the error caused by the deviation of SBFEM approximations of being harmonic inside $S$. The fact that the first source of error is the  dominant one  is favorable for  applications for singular problems, where the singularity  may be isolated,  without interaction  with the $S$-element facets. For this class of problems, the solution away from the
singularity is regular. If the convergence rate is dominated by the approximation
on the boundary of $S$, then this explains regular convergence rates even for
singular problems, as illustrated by the verification  simulations. 

Numerical tests in 2D and 3D problems emphasize the optimal rate of convergence of the scaled boundary approximations, proven theoretically for the energy norm. Although we have considered only harmonic solutions, the demonstration can be extended for more general homogeneous elliptic PDEs, for instance, elasticity problems without body loads.

\section*{Acknowledgements}

 The authors thankfully acknowledge financial support from: FAPESP - S\~ao Paulo Research Foundation, grants  2016/05155-0 (Gomes) and  17/08683-0 (Devloo),    CNPq - Conselho Nacional de Desenvolvimento Cient\'\i fico e Tecnol\'ogico, grants 305823-2017-5 (Devloo) and  306167/2017-4 (Gomes), and   ANP - Brazilian National Agency of Petroleum, Natural  Gas and Biofuels, grant 2014/00090-2 (Coelho, Devloo).

\bibliographystyle{elsarticle-num-names}
\bibliography{CMAMEReferences.bib}

\appendix

\section{Scaled ODE equation}\label{Riccati}
The second-order ODE problem  \eqref{ODE} can be solved using standard   methods through a system of first-order differential equations. Given  $\hat{\vetor{Q}}_{i}(\xi)=\left[  \xi^{d-1}  \tensor{E}_{11} \hat{\vetor{\Phi}}_i^{'}(\xi) +   \xi^{d-2} \tensor{E}_{21} \hat{\vetor{\Phi}}_i(\xi)  \right]$, the ODE \eqref{ODE}  can be expressed by the two equations:
{\small
\begin{align}
    \xi \vetor{\hat{\Phi}}'_i(\xi) &= \left(-\tensor{E}_{11}^{-1}\tensor{E}_{12} + 0.5(d-2)\mathbf{I}\right)\vetor{\hat{\Phi}}_i(\xi) + \tensor{E}_{11}^{-1}\vetor{\hat{Q}}_i(\xi), \label{eq:odedu}\\
    \xi\vetor{\hat{Q}}'_i(\xi) &= \left(-\tensor{E}_{21}\tensor{E}_{11}^{-1}\tensor{E}_{12} + \tensor{E}_{22}\right)\vetor{\hat{\Phi}}_i(\xi)  + \left(\tensor{E}_{11}\tensor{E}_{21}^{-1} - 0.5(d-2)\mathbf{I}\right)\hat{\vetor{Q}}_i(\xi). \label{eq:odedq}
\end{align}
This ODE system can be  grouped in a matrix form as 
\begin{equation}   \xi \tensor{X}'(\xi) =- \tensor{Z}\,\tensor{X}(\xi),  \;\xi \in  [-1,1], \label{eq:eigenvalueproblem} \end{equation}
for  $ \tensor{X}(\xi) = \begin{bmatrix}   \tensor{\hat{\Phi}}(\xi) \\  \tensor{\hat{Q}}(\xi)\end{bmatrix}$, where
 $\tensor{\hat{\Phi}}(\xi)=[ \vetor{\hat{\Phi}}_i(\xi)]$,  and $\tensor{\hat{Q}}(\xi)=[\vetor{\hat{Q}}_i(\xi)]$  are $\mathcal{N}^S\times \mathcal{N}^S$ matrices with columns $ \vetor{\hat{\Phi}}_i(\xi)$ and $\vetor{\hat{Q}}_i(\xi)$,
and  $\tensor{Z}$ is the $2\mathcal{N}^S \times 2\mathcal{N}^S$ matrix
 {\small
$$
    \tensor{Z} = \begin{bmatrix}
    \left(\tensor{E}_{11}^{-1}\tensor{E}_{12} - 0.5(d-2)\mathbf{I}\right) & -\tensor{E}_{11}^{-1}\\
    -\tensor{E}_{22}+ \tensor{E}_{21}\tensor{E}_{11}^{-1}\tensor{E}_{21} & \left(-\tensor{E}_{21}\tensor{E}_{11}^{-1} + 0.5(d-2)\mathbf{I}\right)
    \end{bmatrix}.
$$}
If  $\begin{bmatrix}\tensor{A} \\  \tensor{Q}\end{bmatrix}$   are linearly independent eigenvectors  of the matrix $\tensor{Z}$ corresponding  to  eigenvalues $\vetor{\lambda}$, then the function 
$    \tensor{X}(\xi) =\begin{bmatrix}\tensor{A} \\  \tensor{Q}\end{bmatrix} \xi^{\vetor{\lambda}}$
solves  \eqref{eq:eigenvalueproblem}.
The functions $ \xi^{\vetor{\lambda}}$ corresponding to   eigenvalues  having  negative 
 real parts  are unbounded for 
$\xi\rightarrow 0$, and are unsuited to describe solutions 
at the interior of the
S-element, whilst those of positive  real parts represent solutions
that are zero at the scaling center of $S$.
Thus, the  desired solutions of  the system \eqref{eq:odedu}-\eqref{eq:odedq} are taken  as
$$
 \tensor{\hat{\Phi}}(\xi) = \tensor{A}_+\mbox{diag}(\xi^{\vetor{\lambda}_{+}}),\quad   \hat{\tensor{Q}}(\xi) = \tensor{Q}_+\mbox{diag}(\xi^{\vetor{\lambda}_{+}}),$$
where $\vetor{\lambda}_{+} \in \mathbb{R}^N$ represents the positive real part of $\vetor{\lambda}$,  $\tensor{A}_+=[\vetor{A}_{+i}]$ and  $\tensor{Q}_+=[\vetor{Q}_{+i}]$ are the associated eigenvector components.  For simplicity, the index $+$ is dropped in Section \ref{sec:ODE}.
\end{document}